\documentclass[aap]{imsart}
\makeatletter\def\journal@name{}\makeatother % no "Submitted to ..." line
\RequirePackage{cmap}
\input{glyphtounicode}\pdfgentounicode=1
\RequirePackage{amsthm,amsmath,amssymb,amsfonts}
\RequirePackage[numbers]{natbib}
\RequirePackage{mathtools,booktabs,array,enumitem,adjustbox}
\RequirePackage[colorlinks,citecolor=blue,linkcolor=blue,urlcolor=blue]{hyperref}
\startlocaldefs
\theoremstyle{plain}
\newtheorem{theorem}{Theorem}
\newtheorem{proposition}{Proposition}
\newtheorem{lemma}{Lemma}
\newtheorem{tool}{Theorem}
\makeatletter\newcommand{\tlsec}[1]{{\def\@currentlabel{\ifnum\value{subsection}>0 \thesubsection\else\thesection\fi}\label{ts:#1}}}\makeatother
\newcommand{\TS}[1]{Theorem~\ref{tl:#1} (Section~\ref{ts:#1})}

\theoremstyle{definition}
\newtheorem{definition}{Definition}
\newtheorem{remark}{Remark}

\newcommand{\E}{\mathbb E}\renewcommand{\P}{\mathbb P}\newcommand{\Var}{\operatorname{Var}}
\newcommand{\1}{\mathbf 1}\newcommand{\R}{\mathbb R}\newcommand{\T}{\mathbb T}
\newcommand{\Dbar}{\mu_n}\newcommand{\dK}{d_{\mathrm K}}

\newcommand{\Line}[1]{\[\adjustbox{max width=\linewidth}{$\displaystyle #1$}\]}
\newcommand{\LineE}[2]{\begin{equation}\label{#1}\adjustbox{max width=0.9\linewidth}{$\displaystyle #2$}\end{equation}}
\newcommand{\PtD}[1]{\E_1[A_{ki}A_{il}\mid X_k=#1,X_l=0]}
\newcommand{\PthD}[1]{\E_1[A_{ki}A_{ij}A_{jl}\mid X_k=#1,X_l=0]}
\numberwithin{equation}{section}
\newenvironment{restated}[3]{\par\smallskip\noindent\textbf{#1~\ref{#2}#3.}\ \itshape}{\par\smallskip\upshape}
\renewcommand{\paragraph}[1]{\par\smallskip\noindent\textbf{#1}\ }
\endlocaldefs

\begin{document}

\begin{frontmatter}
\title{Limit laws for the friendship paradox in random geometric graphs}
\runtitle{Limit laws for the friendship paradox}

\begin{aug}
\author[A]{\fnms{Mingao}~\snm{Yuan}\ead[label=e1]{myuan2@utep.edu}}
\address[A]{Department of Mathematical Sciences, The University of Texas at El Paso, El Paso, Texas 79968, USA\printead[presep={,\ }]{e1}}
\end{aug}

\begin{abstract}
The friendship paradox says that, on average, the neighbours of a vertex have more neighbours than the vertex itself. We study the size of this effect, the average over all vertices of the mean degree of the neighbours minus the degree of the vertex, in the random geometric graph on the torus with the maximum-norm distance. We determine its limit law in every regime of the connection radius and quantify the approximation by non-asymptotic error bounds in the Kolmogorov distance or the total variation distance. We also give its exact mean for every $n$ and every radius below one half. For a fixed radius the limit is not Gaussian but a weighted sum of centred exponential variables; when the radius shrinks and the expected degree grows, the statistic is asymptotically Gaussian, with a Berry--Esseen bound; for bounded expected degree it is Gaussian at the classical rate; and when the expected degree tends to zero it is asymptotically Gaussian, Poisson or identically zero, depending on how fast the radius shrinks. When the expected degree is large, the proofs rest on the vanishing of the first Hoeffding projection, which makes a degenerate $U$-statistic of pairs of points the source of the fluctuations; it is analysed by Fourier series for a fixed radius and by Stein's method of exchangeable pairs for a shrinking radius. For bounded expected degree we use normal approximation for local functionals, and for vanishing expected degree the Chen--Stein method.
\end{abstract}

\begin{keyword}[class=MSC]
\kwdgroup[type=primary]{\kwd{60F05}
\kwd{60D05}}
\kwdgroup[type=secondary]{\kwd{05C80}}
\end{keyword}

\begin{keyword}
\kwd{friendship paradox}
\kwd{random geometric graph}
\kwd{Berry--Esseen bound}
\kwd{degenerate $U$-statistic}
\kwd{exchangeable pairs}
\kwd{Chen--Stein method}
\end{keyword}

\end{frontmatter}

\section{Introduction}
The \emph{friendship paradox}, discovered by the sociologist Feld~\cite{Feld91}, is often phrased as ``your friends have more friends than you do''. To make this precise, consider a group of individuals joined by mutual friendships. For each individual, compare the average number of friends of his or her friends with his or her own number of friends, and average this difference over the whole group. The result is never negative, and it is strictly positive unless every two friends have the same number of friends. The effect, although counterintuitive at first glance, is thus a mathematical fact and not a statistical accident.

In mathematical terms, let $\mathcal{G}$ be a graph on $n$ vertices with adjacency matrix $(A_{ij})$ and degrees $d_i=\sum_{j\ne i}A_{ij}$. The local version of the paradox is described by the \emph{friendship index}
\begin{equation}\label{eq:Delta}
\Delta_i=\frac1{d_i}\sum_{j\ne i}A_{ij}d_j-d_i\qquad(\Delta_i=0\text{ if }d_i=0),
\end{equation}
the mean degree of the neighbours of vertex $i$ minus the degree of $i$. The global version of the paradox is described by the average
\[\mathcal{F}_n=\frac1n\sum_{i=1}^n\Delta_i .\]
Every edge is seen from both of its ends, and a short computation (Proposition~\ref{prop:id} in Section~\ref{ts:prop:id}) gives
\begin{equation}\label{eq:identity}
\mathcal{F}_n=\frac1n\sum_{i<j}A_{ij}\frac{(d_i-d_j)^2}{d_id_j}\ \ge0 .
\end{equation}
Hence the paradox holds in every graph, and the size of $\mathcal{F}_n$ measures how unequal the degrees of adjacent vertices are.

The paradox is more than a curiosity: it can be turned into a sampling and monitoring device. Immunising random friends protects a network against epidemics far more efficiently than immunising random individuals (Cohen, Havlin and ben-Avraham~\cite{CHbA03}); friends of random individuals act as ``sensors'' that detect outbreaks earlier (Christakis and Fowler~\cite{CF10}); seeding random neighbours improves influence maximisation (Lattanzi and Singer~\cite{LS15}); sampling friends yields efficient algorithms for polling a social network (Nettasinghe and Krishnamurthy~\cite{NK19}) and for maximum likelihood estimation of power-law degree distributions (Nettasinghe and Krishnamurthy~\cite{NK21}); and the paradox explains systematic biases in perceptions of popularity and social norms, and with them why highly connected agents behave differently from their less connected peers, for instance in adolescent drug use and binge drinking (Hodas, Kooti and Lerman~\cite{HKL13}; Jackson~\cite{Jackson19}). The paradox has also been generalised in several directions: to edge weights and node attributes (Eom and Jo~\cite{EJ14}; Evtushenko and Kleinberg~\cite{EK24}), to triangles (Bhattacharya et al.~\cite{BGHJS26}) and to network centrality indices, where the average centrality of the neighbours of a vertex exceeds the global average (Hazra and Verbitskiy~\cite{HV26}). In all these uses, the size of the gap between the mean degree of the neighbours and the degree of the vertex matters, not only its sign. This size is exactly what $\mathcal{F}_n$ measures.

Despite this interest, little is known about the distribution of $\mathcal{F}_n$ in random graph models. For general graphs, Feld~\cite{Feld91} explained the paradox by the variance of the degrees, and Cao and Ross~\cite{CR16} and Jackson~\cite{Jackson19} give general conditions, bounds on its size and generalisations to other vertex attributes. For random graphs, Hazra, den Hollander and Parvaneh~\cite{HdHP23,HdHP25} study the friendship bias of a \emph{typical vertex} in sparse random graphs that converge locally to a random tree (Erd\H{o}s--R\'enyi graphs, the configuration model, preferential attachment) and identify its limiting distribution. For random geometric graphs, the statistic $\mathcal{F}_n$ itself has been studied by Yuan: \cite{Yuan24} determines its asymptotic distribution for a dense graph,  and a weak law of large numbers for sparse graphs is proven in \cite{Yuan26}. A description of the fluctuations of $\mathcal{F}_n$ that covers all radii is, however, missing, as is a bound on the maximal error of the approximation of the true distribution of $\mathcal{F}_n$ by its limit distribution. It is not a consequence of the typical-vertex results: unlike the friendship bias of a typical vertex, the fluctuations of the \emph{graph-wide average} $\mathcal{F}_n$ are governed not by a single vertex but by the dependence among vertices, and in a geometric graph this dependence is strong, because neighbouring vertices share neighbours.

We address this question for the random geometric graph, introduced by Gilbert~\cite{Gil61} and treated systematically by Penrose~\cite{Penrose03}. It is the standard model of a network embedded in space, used for wireless and ad hoc networks (Haenggi et al.~\cite{Haenggi09}) and for spatial networks in general (Barth\'elemy~\cite{Barth11}); see Duchemin and De~Castro~\cite{DDC22} for recent developments. Let $X_1,\dots,X_n$ be i.i.d.\ uniform points on the torus $\T^m=[0,1)^m$, let $d$ be the $L_\infty$ distance on $\T^m$, and join $i\ne j$ if $d(X_i,X_j)<r_n$ (the $L_\infty$ random geometric graph of Bangachev and Bresler~\cite{BB24}, Yuan and Yu \cite{YY26}). The torus avoids boundary effects, and the $L_\infty$ distance makes neighbourhoods cubes, so that all probabilities of subgraphs reduce to volumes of intersections of cubes. We write $\mathcal{F}_n$ for the statistic of the resulting graph.

The most distinctive feature of a geometric graph is \emph{transitivity}: two neighbours of a vertex are likely to be neighbours of each other, so the degrees of adjacent vertices are strongly positively correlated (Dall and Christensen~\cite{DC02}; Penrose~\cite{Penrose03}; Barth\'elemy~\cite{Barth11}). This is why the paradox is weaker than in tree-like graphs. In a sparse Erd\H{o}s--R\'enyi graph a neighbour of a vertex has on average one more neighbour than the vertex itself, so $\mathcal{F}_n\to1$ (see Hazra, den Hollander and Parvaneh~\cite{HdHP23}), whereas in the random geometric graph we show that
\begin{equation}\label{eq:mean-intro}
\E[\mathcal{F}_n]\longrightarrow\kappa_m:=1-\Big(\frac34\Big)^m<1
\end{equation}
when the expected degree tends to infinity and $r_n\to0$ (the exact formula, for every $n$ and every $r_n<\frac12$, is Theorem~\ref{thm:mean} in Section~\ref{ts:thm:mean}). The natural questions are how $\mathcal{F}_n$ fluctuates around its limit, how fast it becomes Gaussian, and what replaces the Gaussian limit when the graph is sparse or the radius is constant. Because of the transitivity, the usual tools for sums of independent or locally tree-like terms do not apply directly. Normal approximation for sums of local functionals of point processes is well developed (see, e.g., Lachi\`eze-Rey, Schulte and Yukich~\cite{LSY19}), but it requires the fluctuations to be driven by individual points. Here this is not the case: by translation invariance of the torus, the first Hoeffding projection of $\mathcal{F}_n$ vanishes, so the usual source of fluctuations of order $n^{-1/2}$ is missing, and a degenerate $U$-statistic built from pairs of points takes its place.

In this paper, we describe the limit law of $\mathcal{F}_n$ in every regime of the radius $r_n$, from the dense regime down to the very sparse one, and we quantify the approximation by non-asymptotic error bounds in the Kolmogorov distance or the total variation distance. We also obtain an explicit formula for the mean of $\mathcal{F}_n$ and its asymptotic variance in every regime. Table~\ref{tab:regimes} summarizes the picture: for each regime, determined by the size of the radius $r_n$, it lists the limit law of $\mathcal{F}_n$ and the rate of the approximation, where $\dK$ denotes the Kolmogorov distance and $d_{\mathrm{TV}}$ the total variation distance; the precise statements are in Section~\ref{sec:main}.
\begin{table}[h]\centering\small\renewcommand{\arraystretch}{1.7}
\adjustbox{max width=\textwidth}{\begin{tabular}{rlll}\toprule
 & regime & limit law of $\mathcal{F}_n$ & rate\\\midrule
1 & $nr_n^m\asymp n$ ($r_n$ constant) & non-Gaussian ($Q_\infty$) & $\dK\lesssim n^{-1/2}(\log n)^{3/2}$\\
2 & $nr_n^m\to\infty$, $r_n=o(1)$ & Gaussian & $\dK\lesssim r_n^{m/2}+(nr_n^m)^{-1/3}$\\
3 & $nr_n^m\asymp1$ & Gaussian & $\dK\lesssim n^{-1/2}$\\
4 &$nr_n^m=o(1)$, $n^3r_n^{2m}\to\infty$ & Gaussian & $\dK\lesssim(nr_n^m)^{1/3}+(n^3r_n^{2m})^{-1/2}$\\
5 & $n^3r_n^{2m}\asymp1$ & Poisson (for $n\mathcal{F}_n$) & $d_{\mathrm{TV}}\lesssim n^4r_n^{3m}$\\
6 & $n^3r_n^{2m}=o(1)$ & degenerate at $0$ & $\P(\mathcal{F}_n\ne0)\lesssim n^3r_n^{2m}$\\\bottomrule
\end{tabular}}
\caption{Limit laws and approximation error bounds of \(\mathcal{F}_{n}\) from the dense to the sparse regime.}\label{tab:regimes}
\end{table}

\noindent The first result concerns a constant radius (regime~1 of Table~\ref{tab:regimes}), where the fluctuations are of order one and the limit is not Gaussian: $\mathcal{F}_n-\E[\mathcal{F}_n]$ converges, with Kolmogorov rate $n^{-1/2}(\log n)^{3/2}$, to
\begin{equation}\label{eq:Q-intro}
Q_\infty=\sum_{\xi\in\mathbb Z^m_+}\lambda_\xi\big(|Z_\xi|^2-1\big),
\end{equation}
where the $Z_\xi$ are i.i.d.\ standard complex Gaussian variables and the weights $\lambda_\xi\ge0$ are explicit Fourier coefficients of a kernel built from the cube (Theorem~\ref{thm:fixed} in Section~\ref{ts:thm:fixed}). The second result concerns the dense regime where $r_n\to0$ and $nr_n^m\to\infty$ (regime~2 of Table~\ref{tab:regimes}). We show that $\mathcal{F}_n$ is asymptotically Gaussian, with the Berry--Esseen bound $r_n^{m/2}+(nr_n^m)^{-1/3}$ in the Kolmogorov distance (Theorems~\ref{thm:var} and~\ref{thm:be}). When the expected degree stays bounded (regime~3 of Table~\ref{tab:regimes}), $\mathcal{F}_n$ is asymptotically Gaussian with the classical rate $n^{-1/2}$ (Theorem~\ref{thm:sparse} in Section~\ref{ts:thm:sparse}). When the expected degree tends to zero (regimes~4--6 of Table~\ref{tab:regimes}) the graph splits into tiny components, and $n\mathcal{F}_n$ is essentially the number of induced two-paths: depending on whether the expected number of these, of order $n^3r_n^{2m}$, tends to zero, stays bounded or tends to infinity, the statistic vanishes with high probability (regime~6 of Table~\ref{tab:regimes}), is asymptotically Poisson (regime~5 of Table~\ref{tab:regimes}), or is asymptotically Gaussian (regime~4 of Table~\ref{tab:regimes}) (Theorem~\ref{thm:vsparse} in Section~\ref{ts:thm:vsparse}).

Together, these results show that one and the same statistic has a non-Gaussian, a Gaussian, a Poisson or a degenerate limit, and that the transitions are located by two parameters, the expected degree $nr_n^m$ and the expected number of induced two-paths $n^3r_n^{2m}$. The non-Gaussian limit (regime~1 of Table~\ref{tab:regimes}) and the Gaussian limit in the dense regime (regime~2 of Table~\ref{tab:regimes}) are two faces of the same degenerate pair statistic: for a fixed radius its kernel has a fixed discrete Fourier spectrum, whereas for a shrinking radius the spectrum spreads over many modes, and heuristically no single mode dominates.

To prove these results, we combine the exact identity~\eqref{eq:identity} with standard tools. The Efron--Stein inequality controls the remainder that arises when the degrees $d_id_j$ in the denominator of~\eqref{eq:identity} are replaced by the squared expected degree, and the Hoeffding decomposition splits the resulting polynomial in subgraph counts into pieces that depend on one, two, three and four points. The one-point piece vanishes by translation invariance, and the two-point piece, a degenerate $U$-statistic, carries the fluctuations. For a fixed radius, the limit theory of degenerate $U$-statistics together with Fourier analysis on the torus gives the non-Gaussian limit. For a shrinking radius and large expected degree, Stein's method of exchangeable pairs gives the Gaussian limit with a Berry--Esseen bound. For bounded expected degree, the normal approximation for stabilising functionals gives the rate $n^{-1/2}$. For vanishing expected degree, the Chen--Stein method gives the Poisson regime, and, together with the Berry--Esseen theorem, the Gaussian one.

 The rest of the paper is organised as follows. Section~\ref{sec:main} states all results; Section~\ref{sec:guide} explains the proof strategy; Sections~\ref{sec:proofmean}--\ref{sec:vsparse} contain the proofs (exact mean, dense regime with the main machinery, constant radius, bounded expected degree, very sparse regimes); the proofs of the lemmas and propositions are given in the Supplement.

\section{Main results}\label{sec:main}
This section states the results. We first fix conventions (asymptotic notation and the laws used), then define the model and the statistic and introduce the two distances used to compare laws. Next comes the exact mean, which holds for all radii and gives the centring used in every later result. The remaining subsections describe the fluctuations of $\mathcal{F}_n$ from the densest regime to the sparsest, in the order of Table~\ref{tab:regimes}: constant radius, shrinking radius with growing expected degree, bounded expected degree, and vanishing expected degree.

\paragraph{Asymptotic notation.} All limits are as $n\to\infty$. The dimension $m\ge1$ is a fixed positive integer throughout the paper, and the radius $r_n\in(0,\frac12)$ may depend on $n$. For positive sequences $a_n,b_n$ (possibly depending on $r_n$ as well) we use the Bachmann--Landau notation:
\begin{itemize}[leftmargin=*,itemsep=0pt]
\item $a_n=O(b_n)$ (equivalently $a_n\lesssim b_n$) means $a_n\le Cb_n$ for a constant $C$ and all large $n$;
\item $a_n=o(b_n)$ (equivalently $a_n\ll b_n$) means $a_n/b_n\to0$;
\item $a_n\asymp b_n$ means $cb_n\le a_n\le Cb_n$ for constants $0<c\le C$ and all large $n$ (that is, $a_n=O(b_n)$ and $b_n=O(a_n)$);
\item $a_n\sim b_n$ means $a_n/b_n\to1$.
\end{itemize}

\paragraph{Distributions.} For a random variable $X$ and a law $\mathcal L$, $X\sim\mathcal L$ means that $X$ has law $\mathcal L$. We use the following laws: $\mathrm{Bin}(N,p)$, the binomial law with $N$ trials and success probability $p$, so $\P(X=k)=\binom Nk p^k(1-p)^{N-k}$, $k=0,\dots,N$; $\mathrm{Po}(\lambda)$, the Poisson law with mean $\lambda>0$, so $\P(X=k)=e^{-\lambda}\lambda^k/k!$, $k=0,1,\dots$ ($\mathrm{Po}(0)$ is the point mass at $0$); $\mathrm{Exp}(1)$, the exponential law with mean $1$; and $\mathcal{N}(0,1)$, the standard normal law, with $\mathcal{N}$ denoting a random variable having this law. The notation $\mathrm{Po}(\lambda_n)$ is used when the mean depends on $n$.    $\Phi$ denotes the standard normal distribution function. Further, $\Rightarrow$ denotes convergence in distribution.

\subsection{The model and the friendship paradox }
\begin{definition}[Random geometric graph on the torus]\label{def:model}
Let $X_1,\dots,X_n$ be i.i.d.\ uniform on $[0,1)^m$ with opposite sides identified, let $A_{ij}=\1\{d(X_i,X_j)<r_n\}$ for $i\ne j$, where
\Line{d(x,y)=\max_{1\le b\le m}\min(|x_b-y_b|,1-|x_b-y_b|).} The graph $\mathcal{G}_n(r_n)$ on the vertex set $\{1,\dots,n\}$ with adjacency matrix $(A_{ij})$ is the random geometric graph. The radius $r_n$ may depend on $n$.
\end{definition}

\noindent The cube $[0,1)^m$ with opposite sides identified is the $m$-dimensional torus, and $d$ is the $L_\infty$ distance on it.  On the torus every point plays the same role, so the law of the graph does not change under translations; we use this invariance throughout. The $L_\infty$ distance makes the neighbourhood of a vertex a cube, so that subgraph probabilities factorise over the coordinates and can be computed explicitly (see Section~\ref{pf:prop:alpha} of the Supplement for example).  The edges in $\mathcal{G}_n(r_n)$ are correlated, because they depend on the distances between the random positions of the vertices. The radius \(r_{n}\) parameterizes the sparsity of the graph, where a larger \(r_{n}\) yields a higher expected edge density. Because the graph becomes complete for \(r_n \ge \frac{1}{2}\), this regime is excluded from our analysis throughout the paper. Furthermore, the dimension \(m\) is assumed to be a fixed positive integer.

Definition~\ref{def:model} is the torus version of the model of Gilbert~\cite{Gil61}, which is studied systematically in~\cite{Penrose03}. The graph $\mathcal{G}_n(r_n)$ is the $L_\infty$ random geometric graph considered by Bangachev and Bresler~\cite{BB24} and by Yuan and Yu~\cite{YY26}; the same distance on the torus is used by Bringmann, Keusch and Lengler~\cite{BKL19} to define their geometric inhomogeneous random graphs. The model is also closely related to that of Goel, Rai and Krishnamachari~\cite{GRK05}, who place $n$ points independently and uniformly in the unit cube $[0,1]^m$ and join two points if their Euclidean distance is at most a prescribed radius. Other models replace the torus by the sphere and join two points when they are close on the sphere; see Devroye, Gy\"orgy, Lugosi and Udina~\cite{DGLU11} and Bubeck, Ding, Eldan and R\'acz~\cite{BDER16}. For further models and recent developments, see Duchemin and De~Castro~\cite{DDC22}.

\begin{definition}[Friendship paradox statistic]\label{def:stat}
Let $\mathcal{G}$ be a graph on the vertex set $\{1,\dots,n\}$ with adjacency matrix $(A_{ij})$. Let\[\Delta_i=\frac1{d_i}\sum_{j\ne i}A_{ij}d_j-d_i,\]  where $d_i=\sum_{j\ne i}A_{ij}$  (the degree of vertex $i$), and  $\Delta_i=0$ if $d_i=0$. The Friendship paradox $\mathcal{F}(\mathcal{G})$ of the graph $\mathcal{G}$ is defined as
\[\mathcal{F}(\mathcal{G})=\frac1n\sum_{i=1}^n\Delta_i .\]
\end{definition}

The quantity $\Delta_i$, the mean degree of the neighbours of $i$ minus the degree of $i$, is called the friendship bias by Hazra, den Hollander and Parvaneh~\cite{HdHP23} and by Bhattacharya et al.~\cite{BGHJS26}, and it is the basic quantity of the friendship paradox in Cantwell, Kirkley and Newman~\cite{CKN21}. A ratio version of $\Delta_i$, called the friendship index, and its network-wide aggregates are studied by Pal, Yu, Novick, Swami and Bar-Noy~\cite{PYNSB19}. The average $\mathcal{F}(G)$ is non-negative for every graph, and this is the mathematical content of the friendship paradox of Feld~\cite{Feld91} (see also~\cite{CR16,CKN21}). The paradox has been generalised to edge weights and node attributes (Eom and Jo~\cite{EJ14}; Evtushenko and Kleinberg~\cite{EK24}), to triangles (Bhattacharya et al.~\cite{BGHJS26}) and to network centrality indices (Hazra and Verbitskiy~\cite{HV26}).  For random graphs, the friendship bias of a typical vertex is studied in sparse graphs by Hazra, den Hollander and Parvaneh~\cite{HdHP23,HdHP25}, and for random geometric graphs the statistic $\mathcal{F}_n$ itself has been studied by Yuan~\cite{Yuan24,Yuan26}.

 To measure how close the law of $\mathcal{F}_n$ is to a limit law, we use two standard distances. The first is the Kolmogorov distance, the largest gap between two distribution functions: it is the largest error made when a probability $\P(X\le x)$ is replaced by the corresponding probability for $Y$. For real random variables $X,Y$, the
Kolmogorov distance between them is defined as
\Line{\dK(X,Y)=\sup_{x\in\R}\big|\P(X\le x)-\P(Y\le x)\big| .}
The Kolmogorov distance is the standard distance for Berry--Esseen bounds. We use it for the Gaussian limits (regimes 2, 3 and 4) and for the non-Gaussian limit $Q_\infty$ of regime 1.

The second is the total variation distance, the largest gap between the probabilities of arbitrary events:
\Line{d_{\mathrm{TV}}(X,Y)=\sup_{\mathcal{B}\subseteq\R}\big|\P(X\in \mathcal{B})-\P(Y\in \mathcal{B})\big| ,}
where the supremum is over all Borel sets $\mathcal{B}$. Clearly $\dK\le d_{\mathrm{TV}}$, so total variation is the stronger of the two. It is the natural distance when the limit is a discrete law, and we use it in regime 5, where the limit is Poisson: $n\mathcal{F}_n$ is then integer-valued (it counts induced two-paths), and the Chen--Stein method (Theorem~\ref{tl:I} in Section~\ref{ts:I}) controls the probability of every set of integers. 

\subsection{The exact mean}
We begin with the expected value, which is known exactly for every $n$ and every radius below $\frac12$; it also fixes the centring used in all later results.
\begin{theorem}[Exact mean]\label{thm:mean}\tlsec{thm:mean}
Let $n\ge2$. Then the expectation of $\mathcal{F}_n$ is given exactly by
\Line{\E[\mathcal{F}_n]=\kappa(r_n)\Bigg[1-\big(1-(2r_n)^m\big)^{n-1}\Bigg(1+\frac{(n-1)(2r_n)^m}{1-(2r_n)^m}\Bigg)\Bigg],}
where 
\[\kappa(r_n)=1-\frac{\E[A_{12}A_{13}A_{23}]}{(2r_n)^{2m}}=1-\rho(r_n)^m\] and, for $t\in(0,\frac12)$,
\Line{\rho(t)=\begin{cases}\dfrac34,&0<t\le\dfrac13,\\[2mm]\dfrac{12t^2-6t+1}{4t^2},&\dfrac13<t<\dfrac12.\end{cases}}
\end{theorem}
The exact formula of Theorem~\ref{thm:mean}  shows that $\E[\mathcal{F}_n]$ approaches $\kappa(r_n)$ exponentially fast as the expected degree $(n-1)(2r_n)^m$ grows:
\Line{0\le\kappa(r_n)-\E[\mathcal{F}_n]\le\kappa(r_n)\Big(1+\frac{(n-1)(2r_n)^m}{1-(2r_n)^m}\Big)e^{-(n-1)(2r_n)^m}.}
The limit of the mean in each regime (the expected degree is of order $nr_n^m$) is therefore
\Line{\E[\mathcal{F}_n]\longrightarrow\begin{cases}1-\Big(\dfrac{12r^2-6r+1}{4r^2}\Big)^m,&r_n=r\in(\frac13,\frac12)\text{ constant},\\[3mm]1-(\tfrac34)^m,&nr_n^m\to\infty\text{ and }r_n\le\frac13,\\[2mm]\big(1-(\tfrac34)^m\big)\Big[1-\big(1+2^mc\big)e^{-2^mc}\Big],&nr_n^m\to c\in(0,\infty),\\[2mm]0,&nr_n^m\to0.\end{cases}}

\begin{remark}[Comparison with the literature]\label{rem:meancomp}
For $m=1$, Yuan~\cite{Yuan26} obtains $\E[\mathcal{F}_n]=\frac14+o(1)$  when $r_n\to0$ and $nr_n\to\infty$; Theorem~\ref{thm:mean} holds for every $n$, every radius below $\frac12$ and every $m$, and has no remainder term. For sparse Erdős--Rényi graphs with mean degree $\lambda$, Hazra, den Hollander and Parvaneh~\cite{HdHP23} show that the friendship bias of a typical vertex has mean $1-(1+\lambda)e^{-\lambda}$; the third case of the display above is this value multiplied by $\kappa_m$. In fact, for the Erdős--Rényi graph $G(n,p)$ a direct computation gives 
\[\E[\mathcal{F}_n]=(1-p)\big[1-(1-p)^{n-2}(1+(n-2)p)\big],\]
which is the formula of Theorem~\ref{thm:mean} with $\kappa(r_n)$ replaced by $1-p$, the probability that two neighbours of a vertex are not adjacent. Thus the transitivity of the geometric graph enters the mean only through the factor $\kappa$.
\end{remark}

\subsection{Regime 1: constant radius $r_n=r\in(0,\frac12)$}
\label{sec:fixed}
If the radius does not shrink, every vertex has of order $n$ neighbours, and one might hope that $\mathcal{F}_n$ becomes deterministic or Gaussian. Neither happens. The vertices are strongly dependent through their shared neighbours, and $\mathcal{F}_n$ keeps random fluctuations of order one, which converge to the non-Gaussian law $Q_\infty$ described below.

\paragraph{The limit law.} Let $r_n=r$ be constant. The limit is built from the following numbers, indexed by integer vectors $\xi=(\xi_1,\dots,\xi_m)\in\mathbb Z^m$. For $\xi\ne0$ put
\begin{equation}\label{eq:ahat}
\lambda_\xi=\frac{2 \eta_\xi^2\big((2r)^m-\eta_\xi\big)}{(2r)^{2m}},\qquad \eta_\xi=\prod_{b=1}^m\Big(\frac{\sin(2\pi r\xi_b)}{\pi\xi_b} \1\{\xi_b\ne0\}+2r \1\{\xi_b=0\}\Big)
\end{equation}
where \Line{\eta_\xi=\int_{\T^m}\1\{\|x\|_\infty<r\}e^{-2\pi i\langle\xi,x\rangle}dx} is the Fourier coefficient of the indicator $\beta(x)=\1\{d(0,x)<r\}$ of the cube of side $2r$ (Section~\ref{sec:fourier}). One checks $|\eta_\xi|<(2r)^m$ for $\xi\ne0$, so $\lambda_\xi\ge0$, and $\lambda_\xi=\lambda_{-\xi}$. Let $\mathbb Z^m_+$ contain exactly one element of each pair $\{\xi,-\xi\}$, $\xi\ne0$, and let $(Z_\xi)_{\xi\in\mathbb Z^m_+}$ be i.i.d.\ standard complex Gaussians ($\E[|Z_\xi|^2]=1$, so $|Z_\xi|^2\sim\mathrm{Exp}(1)$). Put
\begin{equation}\label{eq:Qinf}
Q_\infty=\sum_{\xi\in\mathbb Z^m_+}\lambda_\xi\big(|Z_\xi|^2-1\big),\qquad \sigma^2=\sum_{\xi\in\mathbb Z^m_+}\lambda_\xi^2.
\end{equation}
These are explicit functions of $m$ and $r$ alone.

\paragraph{The constants $\alpha_1(r)$, $\alpha_2(r)$, $\alpha_3(r)$.} The variance $\sigma^2$ of the limit has a closed form in terms of three piecewise polynomials of $r\in(0,\frac12)$ alone:
\LineE{eq:alpha1}{\alpha_1(r)=\begin{cases}\frac{16}3r^3,&0<r\le\frac14,\\[1mm]-\frac13+4r-16r^2+\frac{80}3r^3,&\frac14\le r<\frac12,\end{cases}}
\LineE{eq:alpha2}{\alpha_2(r)=\begin{cases}\frac{115}{12}r^4,&0<r\le\frac15,\\[1mm]\frac1{12}-\frac53r+\frac{25}2r^2-\frac{125}3r^3+\frac{185}3r^4,&\frac15\le r\le\frac13,\\[1mm]-\frac13+\frac{10}3r-10r^2+\frac{10}3r^3+\frac{335}{12}r^4,&\frac13\le r\le\frac25,\\[1mm]1-10r+40r^2-80r^3+80r^4,&\frac25\le r<\frac12,\end{cases}}
\LineE{eq:alpha3}{\alpha_3(r)=\begin{cases}\frac{88}5r^5,&0<r\le\frac16,\\[1mm]-\frac1{60}+\frac r2-6r^2+36r^3-108r^4+\frac{736}5r^5,&\frac16\le r\le\frac14,\\[1mm]\frac1{12}-\frac32r+10r^2-28r^3+20r^4+\frac{224}5r^5,&\frac14\le r\le\frac13,\\[1mm]-\frac9{20}+\frac{13}2r-38r^2+116r^3-196r^4+\frac{872}5r^5,&\frac13\le r<\frac12.\end{cases}}

\begin{theorem}[Constant radii]\label{thm:fixed}\tlsec{thm:fixed}
Let $r_n=r$ be constant, and let $Q_\infty,\sigma^2$, $\alpha_1(r)$, $\alpha_2(r)$, $\alpha_3(r)$ be as in~\eqref{eq:ahat}--\eqref{eq:alpha3}. Then:
\begin{enumerate}[label=(\roman*),leftmargin=*]
\item (\emph{The limit})  $\sum_{\xi\in\mathbb Z^m_+}\lambda_\xi=\kappa(r)$ is finite, the series for $Q_\infty$ converges almost surely and in $L^2$, and $Q_\infty$ has a bounded continuous density. Moreover,
\Line{\E [Q_\infty]=0,\qquad\Var [Q_\infty]=\sigma^2\in(0,\infty),\qquad Q_\infty\ge-\kappa(r).}
 For every $r\in(0,\frac12)$,
\LineE{eq:sigmaformula}{\sigma^2=\frac{2\big[(2r)^{2m}\alpha_1(r)^m-2(2r)^m\alpha_2(r)^m+\alpha_3(r)^m\big]}{(2r)^{4m}},}
which, for \(r < \frac{1}{6}\), yields the simplified form
\Line{\sigma^2=2\big[(\tfrac23)^m-2(\tfrac{115}{192})^m+(\tfrac{11}{20})^m\big](2r)^m.}
\item (\emph{Variance}) The variance of $\mathcal{F}_n$ converges to $\sigma^2$ at the rate
\Line{\big|\Var[\mathcal{F}_n]-\sigma^2\big|\lesssim n^{-1/2}.}
\item (\emph{Limit law and rate}) The Kolmogorov distance between the centred statistic and $Q_\infty$ satisfies
\Line{\dK\big(\mathcal{F}_n-\E[\mathcal{F}_n],Q_\infty\big)\lesssim n^{-1/2}(\log n)^{3/2}.}
\end{enumerate}
\end{theorem}

Part~(i) describes the limit: $Q_\infty$ is a weighted sum of independent centred exponential variables with explicit weights (Fourier coefficients of a kernel built from the cube); it is skewed to the right, hence not Gaussian, and has a smooth density. Part~(ii) says that the variance of $\mathcal{F}_n$ converges to that of the limit. Part~(iii) states the convergence in distribution with an explicit polynomial rate. The rate $n^{-1/2}(\log n)^{3/2}$ is the size of the three- and four-point components of $\mathcal{F}_n$ (which are not part of the limit), up to a power of the logarithm.

\begin{remark}[Comparison with Yuan~\cite{Yuan24}]\label{rem:fixedcomp}
For the random geometric graph on the sphere with the Euclidean distance and a constant radius, Yuan~\cite{Yuan24} shows that $2\mu^2(\mathcal{F}_n-\E[\mathcal{F}_n])$ converges in distribution to $\sum_j\lambda_j(Z_j^2-1)$, where $\mu$ is the edge probability, the $Z_j$ are i.i.d.\ standard normal and the $\lambda_j$ are the eigenvalues of the pair kernel $g$. Theorem~\ref{thm:fixed}  has the same structure. On the torus the kernel is a convolution, so its eigenvalues are explicit Fourier coefficients. This gives a closed form for $\sigma^2$, a proof that $Q_\infty$ is not Gaussian, and the Kolmogorov bound of part~(iii); none of these is available in~\cite{Yuan24}.
\end{remark}

\subsection{Regime 2: $r_n=o(1)$ and $nr_n^m\to\infty$ (the theorems hold for $nr_n^m\ge1$)}
When the radius shrinks but the expected degree still grows, the fluctuations of $\mathcal{F}_n$ become small, and we want their exact size and their limit law. In Theorems~\ref{thm:var} and~\ref{thm:be} below the only conditions on the radius are $r_n<\frac16$ and $nr_n^m\ge1$. The regime of interest is $r_n\to0$ and $nr_n^m\rightarrow\infty$.
\begin{theorem}[Variance]\label{thm:var}\tlsec{thm:var}
Fix $r_{\max}\in(0,\frac12)$ and let $r_n\le r_{\max}$ and $n r_n^m\ge1$. Then
\Line{\Bigg|\frac{\Var[\mathcal{F}_n]}{\sigma^2(r_n)}-1\Bigg|\lesssim(n r_n^m)^{-1/2},}
where $\sigma^2(r)$ is the function defined in~\eqref{eq:sigmaformula} of Theorem~\ref{thm:fixed}(i).
\end{theorem}
The variance of $\mathcal{F}_n$ is, up to a relative error that vanishes with the expected degree, an explicit constant times $(2r_n)^m$ (when $r_n<\frac{1}{6}$), that is, the standard deviation is of order $r_n^{m/2}$. This is much larger than the order $n^{-1/2}$ of an average of $n$ independent terms as soon as $nr_n^m\to\infty$, and it comes from the pair statistic discussed in Section~\ref{sec:guide}.

\begin{theorem}[Berry--Esseen bound]\label{thm:be}\tlsec{thm:be}
Let $r_n<\frac16$ and $n r_n^m\ge1$. Then the standardised statistic satisfies
\Line{\dK\Big(\frac{\mathcal{F}_n-\E[\mathcal{F}_n]}{\sqrt{\Var[\mathcal{F}_n]}},\mathcal{N}\Big)\lesssim r_n^{m/2}+(n r_n^m)^{-1/3}.}
\end{theorem}
After centring and standardising, $\mathcal{F}_n$ is close to a standard normal variable in Kolmogorov distance, provided that $r_n\to0$ and $nr_n^m\to\infty$. The bound has two sources of error: $r_n^{m/2}$, caused by the (small) skewness of the leading pair statistic, and $(nr_n^m)^{-1/3}$, caused by the remaining lower-order terms of $\mathcal{F}_n$.

\subsection{Regime 3: $nr_n^m\asymp1$}
When $nr_n^m\asymp1$ each vertex has a bounded number of neighbours, the dependence between vertices is purely local, and $\mathcal{F}_n$ behaves like an average of $n$ weakly dependent terms.
\begin{theorem}[Sparse regime]\label{thm:sparse}\tlsec{thm:sparse}
Fix $0<c_{\min}\le c_{\max}<\infty$ and let $c_{\min}\le nr_n^m\le c_{\max}$.  
\begin{enumerate}[label=(\roman*),leftmargin=*]
\item (\emph{Variance}) The variance of $\mathcal{F}_n$ is of order $1/n$:
\Line{\Var[\mathcal{F}_n]\asymp\frac1n.}
\item (\emph{Berry--Esseen bound}) The standardised statistic satisfies
\Line{\dK\Big(\dfrac{\mathcal{F}_n-\E[\mathcal{F}_n]}{\sqrt{\Var[\mathcal{F}_n]}},\mathcal{N}\Big)\lesssim n^{-1/2}.}
\end{enumerate}
\end{theorem}
The variance has the usual order $1/n$, as for a sum of $n$ independent terms, and the standardised $\mathcal{F}_n$ is asymptotically Gaussian with the classical rate $n^{-1/2}$.
\subsection{Regimes 4--6: $nr_n^m\to0$ (very sparse)}
Here the graph consists mostly of isolated vertices and isolated edges, which contribute nothing to $n\mathcal{F}_n$; the first contribution comes from \emph{induced two-paths}: triples $\{i,j,k\}$ of vertices with exactly two edges among the three pairs. Each induced two-path contributes exactly $1$ to $n\mathcal{F}_n$, and larger components are negligible (Lemmas~\ref{lem:comp} and~\ref{lem:csize}), so $n\mathcal{F}_n$ is essentially the number $N_\vee$ of induced two-paths. Section~\ref{sec:vsparse} shows
\[
\lambda_n:=\E [N_\vee]=\tfrac12\big(1-(\tfrac34)^m\big) n(n-1)(n-2) (2r_n)^{2m}.
\]

\begin{theorem}[Very sparse regime]\label{thm:vsparse}\tlsec{thm:vsparse}
Suppose $n(2r_n)^m\le\frac1{2e}$.
\begin{enumerate}[label=(\roman*),leftmargin=*]
\item (\emph{Too sparse}) The statistic is non-zero with probability at most
\Line{\P(\mathcal{F}_n\ne0)\lesssim n^3r_n^{2m}.}
\item (\emph{Poisson regime}) The total variation distance between $n\mathcal{F}_n$ and the Poisson law satisfies
\Line{d_{\mathrm{TV}}\big(n\mathcal{F}_n,\mathrm{Po}(\lambda_n)\big)\lesssim n^4r_n^{3m}.}
\item (\emph{Variance and normal regime}) The variance of $n\mathcal{F}_n$ satisfies
\Line{\Big|\frac{\Var[n\mathcal{F}_n]}{\lambda_n}-1\Big|\lesssim(nr_n^m)^{1/2},}
and the standardised statistic satisfies 
\Line{\dK\Bigg(\frac{n\mathcal{F}_n-\E[n\mathcal{F}_n]}{\sqrt{\lambda_n}},\mathcal{N}\Bigg)\lesssim(nr_n^m)^{1/3}+(n^3r_n^{2m})^{-1/2}.}
\end{enumerate}
\end{theorem}
The three cases of the theorem are the regimes 4, 5, 6 of Table~\ref{tab:regimes}: (iii) is the case $n^3r_n^{2m}\to\infty$, (ii) the case $n^3r_n^{2m}\asymp1$ and (i) the case $n^3r_n^{2m}\to0$. In terms of the radius, the consequences are as follows. If $r_n^m\ll n^{-3/2}$, that is, $n^3r_n^{2m}\to0$, then part~(i) gives
\Line{\P(\mathcal{F}_n=0)\to1.}
If $r_n^m\asymp n^{-3/2}$, that is, $n^3r_n^{2m}\to\ell\in(0,\infty)$, then part~(ii) gives that $n\mathcal{F}_n$ is asymptotically Poisson:
\Line{n\mathcal{F}_n\Rightarrow\mathrm{Po}\Big(\tfrac12\big(1-(\tfrac34)^m\big)4^m\ell\Big).}
If $n^{-3/2}\ll r_n^m\ll n^{-1}$, that is, $n^3r_n^{2m}\to\infty$ and $nr_n^m\to0$, then part~(iii) gives that $\mathcal{F}_n$ is asymptotically Gaussian with standard deviation
\Line{\sqrt{\lambda_n}/n\sim(2r_n)^m\sqrt{\tfrac12\big(1-(\tfrac34)^m\big)n}.}

\section{Proof strategy}\label{sec:guide}
\noindent This section outlines the proof strategies for each of the main theorems presented in Section~\ref{sec:main}.

\paragraph{Two basic facts.} \emph{An exact identity.} Adding up the friendship index over all vertices, every edge is seen from both of its ends, and one obtains
\[
n\mathcal{F}_n=\sum_{\{i,j\}\text{ edge}}\frac{(d_i-d_j)^2}{d_id_j}
\]
(Proposition~\ref{prop:id} in Section~\ref{ts:prop:id}). Each edge contributes a non-negative amount, which is zero exactly when its two endpoints have the same degree: the friendship paradox measures how unequal the degrees of adjacent vertices are. \emph{Translation invariance.} On the torus the law of the graph does not change when everything is shifted, and the neighbourhood of a point is always a cube of the same volume. So every other point is a neighbour with the same probability, every degree is binomial, and one-point computations do not depend on where the point is.

\paragraph{Theorem~\ref{thm:mean} (exact mean).} By symmetry it is enough to look at one vertex. Place it at the origin and condition on who its neighbours are: given this, the other points are independent, so the expected degree of a neighbour is a sum of probabilities that two cubes overlap. These probabilities are volumes of intersections of cubes and can be computed exactly (Lemma~\ref{lem:kappa} in Section~\ref{ts:lem:kappa}). No limit theorem is needed (Section~\ref{sec:proofmean}).

\paragraph{The common scheme for Theorems~\ref{thm:fixed}--\ref{thm:be} (dense and intermediate regimes).} Here the expected degree is large and the three theorems are proved in the same three steps (Section~\ref{sec:proofs}; Steps~1--5 hold for every radius $r_n\le r_{\max}<\frac12$, so they also cover a constant radius, Section~\ref{sec:transfer}).
\begin{enumerate}[leftmargin=*,itemsep=2pt]
\item \emph{Replace the random denominators by a constant.} All degrees are close to the expected degree, so in the identity above the products of degrees can be replaced by the square of the expected degree. What is left is a polynomial in the degrees, that is, a combination of subgraph counts, plus a small remainder. The remainder is controlled by the Efron--Stein inequality (\TS{C}): moving one point changes it only slightly (Step~5).
\item \emph{Split the polynomial by the number of points involved (Hoeffding decomposition, \TS{D}).} The statistic is written as a sum of pieces that depend on one, two, three or four of the random points, and the pieces are uncorrelated, so their variances simply add. By translation invariance the one-point piece is constant and disappears after centring. The usual source of fluctuations of order $n^{-1/2}$ is therefore missing, and the two-point piece takes over (Steps~2--3).
\item \emph{The pair statistic dominates.} 
 The two-point piece is a sum over pairs of points of a kernel with mean zero, a degenerate $U$-statistic: the randomness comes entirely from the interaction between pairs. Its variance can be computed in closed form from overlaps of neighbourhoods, while the pieces with three and four points, the degree-and-triangle part and the remainder all have variance of order $1/n$ (Steps~2--5). These three bounds use three different methods. \emph{Three- and four-point pieces:} these components are orthogonal to each other and to the pair piece (\TS{D}), so their variances simply add up over all sets of points. A component can be non-zero only when the points involved are close enough to interact, which makes its second moment small; summing over all sets of points and dividing by the normalisation gives a variance of order $1/n$ (Step~3). \emph{Degree-and-triangle part:} here the Efron--Stein inequality (\TS{C}) is used. Re-sampling one point changes this part by at most a constant times the square of the number of points close to its old or new position, and the moments of this number are controlled by the Rosenthal moment inequality for sums of independent indicators (Lemma~\ref{lem:mom} in Section~\ref{ts:lem:mom}, Step~4). \emph{Denominator error:} the Efron--Stein inequality is used again, but moving one point changes the error only slightly if the degrees involved are close to the expected degree. A Chernoff bound (\TS{A}) shows that a vertex with an atypical degree is exponentially rare. When all affected vertices are typical, the change is small; the rare opposite case is handled by a crude polynomial bound together with the Cauchy--Schwarz inequality and the tiny probability of that case (Step~5). Hence the fluctuation of $\mathcal{F}_n$ is a constant multiple of the pair statistic plus a remainder that is small in every moment.
\end{enumerate}
What remains is to determine the size and the law of the pair statistic, and this is where the three theorems differ.

\paragraph{Theorem~\ref{thm:fixed} (constant radius).} The pair statistic has variance of constant order, so the remainder does not matter. Because the kernel has a fixed range, the pair statistic is not approximately normal. The classical limit theorem for degenerate $U$-statistics (Gregory 1977 \cite{Gregory77}, Serfling \cite{Serfling80}, see \TS{F}) shows that it converges to a weighted sum of centred squared normal variables, with weights given by the eigenvalues of the kernel. On the torus the kernel depends only on the difference of the two points, so its eigenfunctions are the Fourier waves and its eigenvalues are Fourier coefficients, which can be computed explicitly (Section~\ref{sec:eig}); Parseval's identity then gives the variance of the limit. The limit is not Gaussian because it is bounded below. The rate $n^{-1/2}(\log n)^{3/2}$ comes from a Berry--Esseen theorem for degenerate $U$-statistics (Bentkus and G\"otze   \cite{BG99}, G\"otze and Tikhomirov  \cite{GT99}, see \TS{G}) and from a moment bound of every order for the remainder, which allows moments of order about $\log n$ (Lemma~\ref{lem:degmom} in Section~\ref{ts:lem:degmom}, \TS{L}).

\paragraph{Theorems~\ref{thm:var} and~\ref{thm:be} (small radius, expected degree at least one).} \emph{Variance.} The variance of the pair statistic is explicit, and the remainder has variance of order $1/n$. The remainder is negligible exactly when the expected degree tends to infinity, and this gives Theorem~\ref{thm:var} (Section~\ref{ts:thm:var}) with a relative error of the order of the inverse square root of the expected degree. \emph{Berry--Esseen bound.} Resample one uniformly chosen point: this gives a second copy of the pair statistic such that the pair is exchangeable, and, because the pair statistic is a sum over pairs, the average change is \emph{exactly} proportional to the statistic itself. The Shao--Zhang bound for exchangeable pairs (Shao and Zhang   \cite{SZ19}, see \TS{K}) then controls the distance to the normal law by the size of the skewness of the pair statistic (Step~6). Adding the small remainder and fixing the normalisation costs a further term, which is handled by a perturbation lemma (Lemma~\ref{lem:pert} in Section~\ref{ts:lem:pert}). The general theorem for local functionals (Lachi\`eze-Rey--Schulte--Yukich  \cite{LSY19}, see \TS{H}) cannot be used here, because its constants are not uniform in the expected degree (Remark~\ref{rem:noG}).

\paragraph{Theorem~\ref{thm:sparse} (bounded expected degree).} Every vertex has only boundedly many neighbours and the dependence is purely local. \emph{Variance, upper bound:} moving one point changes the statistic by a bounded random amount, so the Efron--Stein inequality gives a variance of order $n$. \emph{Variance, lower bound:}  take of order $n$ cubes of side $\asymp r_n$ that are far apart. Given which cube each point lies in and all points outside the cubes, the cube contributions are conditionally independent, so by the law of total variance their conditional variances add. In each cube, with probability bounded below, exactly three isolated points form either an induced two-path (contribution $1$) or an empty graph (contribution $0$), each with probability bounded below. Each cube thus adds a constant, so $\Var[n\mathcal F_n]\ge cn$. \emph{Berry--Esseen bound:} the contribution of each point depends only on points nearby, moments are bounded and the variance has the order of the volume, so the general theorem for local functionals of point processes (Lachi\`eze-Rey--Schulte--Yukich  \cite{LSY19}, see \TS{H}) applies and gives the rate $n^{-1/2}$ (Section~\ref{sec:sparseproofs}).

\paragraph{Theorem~\ref{thm:vsparse} (very sparse regimes).} When the expected degree tends to zero the graph is nearly empty and splits into tiny components. The statistic is a sum over components. A single edge or a triangle contributes nothing, while a path with three vertices contributes exactly one. So, except on the rare event that some component has at least four vertices (Lemmas~\ref{lem:comp} and~\ref{lem:csize}), the statistic equals the number of induced two-paths. This number is a sum of rare indicators with only local dependence and exactly known mean, so the Chen--Stein method (Arratia, Goldstein and Gordon   \cite{AGG89}, see \TS{I}) compares it with a Poisson variable, with an explicit total variation error (Lemmas~\ref{lem:b} and~\ref{lem:meanNvee}). The three sub-regimes correspond to the mean of this Poisson variable tending to zero (the statistic vanishes with high probability), to a positive constant (Poisson limit), or to infinity (Gaussian limit, using the Berry--Esseen bound for the Poisson law, Lemma~\ref{lem:poisBE} in Section~\ref{ts:lem:poisBE}, and a bound for the contribution of larger components, Lemma~\ref{lem:VarD} in Section~\ref{ts:lem:VarD}), see Section~\ref{sec:vsparse}.

\section{The exact mean}\label{sec:proofmean}
\subsection{Edge probabilities, degrees and triangles}\label{sec:notation}\label{sec:identity}
Recall that $A_{ij}=\1\{d(X_i,X_j)<r_n\}$ on the torus $\T^m=[0,1)^m$ with opposite faces identified. The neighbourhood 
\[\mathcal A(x)=\{y:d(x,y)<r_n\}, \ \ \ \ x\in\T^m\]
is a cube of side $2r_n$, and by translation invariance of the torus its volume does not depend on $x$:
\Line{\mathrm{vol} \mathcal A(x)=(2r_n)^m\qquad(r_n\le\tfrac12).}

\paragraph{Edges and degrees.} Given $X_i$, the point $X_j$ is uniform on $\T^m$, so
\[\P(A_{ij}=1\mid X_i)=\mathrm{vol} \mathcal A(X_i)=(2r_n)^m.\] The indicators $A_{ij}$, $j\ne i$, depend on distinct independent points, so given $X_i$ they are independent Bernoulli$((2r_n)^m)$ variables, and
\LineE{eq:degbin}{d_i\sim\mathrm{Bin}(n-1,(2r_n)^m),\qquad\mu_n=\E [d_i]=(n-1)(2r_n)^m .}
The next lemma computes the triangle probability, and hence the constant $\kappa(r_n)$ of Theorem~\ref{thm:mean} (Section~\ref{ts:thm:mean}), for every $r_n<\frac12$.

\begin{lemma}[The triangle probability and $\kappa(r_n)$ for all $r_n<\frac12$]\label{lem:kappa}\tlsec{lem:kappa}
For $0<r_n<\frac12$, 
\Line{\E[A_{12}A_{13}A_{23}]=\begin{cases} (3r_n^2)^m,&0<r_n\le\dfrac13,\\[2mm](12r_n^2-6r_n+1)^m ,&\dfrac13<r_n<\dfrac12.\end{cases}}
Then $\E[A_{12}A_{13}A_{23}]=(2r_n)^{2m}\rho(r_n)^m$,  
with $\rho(t)$ as in Theorem~\ref{thm:mean} (Section~\ref{ts:thm:mean}), so that $\kappa(r_n)=1-\rho(r_n)^m$. 
\end{lemma}
\noindent The proof of Lemma~\ref{lem:kappa} is given in Section~\ref{pf:lem:kappa} of the Supplement.

\subsection{Proof of the exact mean (Theorem~\ref{thm:mean})}

\noindent By Lemma~\ref{lem:kappa} (Section~\ref{ts:lem:kappa}), $\kappa(r_n)=1-\rho(r_n)^m$ is the function of Theorem~\ref{thm:mean} (Section~\ref{ts:thm:mean}).

\smallskip\noindent\emph{Step 1: reductions.} The vertices are exchangeable, so $\E[\Delta_i]$ is the same for all $i$ and $\E[\mathcal{F}_n]=\E[\Delta_1]$.  

\smallskip\noindent\emph{Step 2: the neighbour patterns.} Given $X_1=0$, $A_{1j}$ for $j\ge2$ depend on the different independent points $X_j$, so they are independent Bernoulli$((2r_n)^m)$ variables. For $k\in\{0,\dots,n-1\}$ let
\[S_k=\{A_{12}=\dots=A_{1,k+1}=1,\ A_{1,k+2}=\dots=A_{1,n}=0\}.\]
Then
\[ 
\P(S_k)=(2r_n)^{mk}(1-(2r_n)^m)^{n-1-k}.\]
The event $\{d_1=k\}$ is the disjoint union of the $\binom{n-1}{k}$ events obtained from $S_k$ by permuting the indices $2,\dots,n$; by exchangeability of $X_2,\dots,X_n$ each of them has the same probability and the same value of $\E[\Delta_1\mid\cdot ]$ as $S_k$. Since $\Delta_1=0$ when $k=0$,
\LineE{eq:avgS}{\E[\mathcal{F}_n]=\E[\Delta_1]=\sum_{k=1}^{n-1}\binom{n-1}{k}(2r_n)^{mk}(1-(2r_n)^m)^{n-1-k} \E[\Delta_1\mid S_k].}

\smallskip\noindent\emph{Step 3: the expected degree of a neighbour.} Let $k\ge1$. Since $S_k$ is a product event in the independent points $X_2,\dots,X_n$, given $S_k$ these points remain independent. The neighbours $2,\dots,k+1$ of vertex $1$ play symmetric roles given $S_k$, so 
\[\E[\Delta_1\mid S_k]=\E[d_2\mid S_k]-k.\] Now \Line{d_2=A_{21}+\sum_{l=3}^{k+1}A_{2l}+\sum_{l=k+2}^nA_{2l}} with $A_{21}=1$ on $S_k$, hence
\Line{\E[d_2\mid S_k]=1+(k-1) \E[A_{23}\mid S_k]+(n-1-k) \E[A_{2,k+2}\mid S_k].}
\emph{Reduction to two or three adjacencies.} The variable $A_{23}$ is a function of $(X_2,X_3)$ only, and $A_{2,k+2}$ of $(X_2,X_{k+2})$ only. The part of $S_k$ that involves the other points, \Line{\{A_{1j}=1,\ 4\le j\le k+1\}\cap\{A_{1j}=0,\ j\ge k+3\}} (resp.\ with $j\notin\{2,k+2\}$), is a function of the remaining points and is therefore independent of $A_{23}$ (resp.\ $A_{2,k+2}$) and of the other conditions. Dropping an independent event from the conditioning does not change a conditional expectation, so
\Line{\E[A_{23}\mid S_k]=\P(A_{23}=1\mid A_{12}A_{13}=1),\qquad\E[A_{2,k+2}\mid S_k]=\P(A_{23}=1\mid A_{12}=1,\ A_{13}=0).}
We compute both by the definition of conditional probability.

\emph{Joint probabilities.} For the triangle, Lemma~\ref{lem:kappa} (Section~\ref{ts:lem:kappa}) gives
\[\P(A_{12}A_{13}A_{23}=1)=\E[A_{12}A_{13}A_{23}]=(2r_n)^{2m}\big(1-\kappa(r_n)\big).\]
For the path, conditioning on $X_2$ first,
\Line{\P(A_{12}A_{23}=1)=\E\big[A_{12} \E[A_{23}\mid X_2]\big]=(2r_n)^m \E[A_{12}]=(2r_n)^{2m}.}

\emph{Two neighbours of vertex $1$.}
\Line{\P(A_{23}=1\mid A_{12}A_{13}=1)=\frac{\P(A_{12}A_{13}A_{23}=1)}{\P(A_{12}A_{13}=1)}=\frac{\E[A_{12}A_{13}A_{23}]}{(2r_n)^{2m}}=1-\kappa(r_n) .}
\emph{A neighbour and a non-neighbour.} Since $\mathbf 1\{A_{13}=0\}=1-A_{13}$,
\Line{\P(A_{12}=1,A_{13}=0,A_{23}=1)=\P(A_{12}A_{23}=1)-\P(A_{12}A_{13}A_{23}=1)=(2r_n)^{2m}\kappa(r_n),}
\[\begin{aligned}\P(A_{23}=1\mid A_{12}=1,A_{13}=0)&=\frac{(2r_n)^{2m}\kappa(r_n)}{\P(A_{12}=1) \P(A_{13}=0)}\\&=\frac{(2r_n)^{2m}\kappa(r_n)}{(2r_n)^m(1-(2r_n)^m)}=\frac{(2r_n)^m\kappa(r_n)}{1-(2r_n)^m}.\end{aligned}\]
Therefore
\[\E[d_2\mid S_k]=1+(k-1)(1-\kappa(r_n))+(n-1-k)\frac{(2r_n)^m\kappa(r_n)}{1-(2r_n)^m},\]
and, subtracting $k$ and using $(n-1-k)(2r_n)^m=\Dbar-k(2r_n)^m$,
\begin{equation}\label{eq:condDelta}\begin{aligned}\E[\Delta_1\mid S_k]=\kappa(r_n)\Big[1-k+\frac{\Dbar-k(2r_n)^m}{1-(2r_n)^m}\Big]=\kappa(r_n)\Big[1+\frac{\Dbar-k}{1-(2r_n)^m}\Big].\end{aligned}\end{equation}

\smallskip\noindent\emph{Step 4: average over the pattern.} Insert~\eqref{eq:condDelta} in~\eqref{eq:avgS}. The weights \Line{\binom{n-1}{k}(2r_n)^{mk}(1-(2r_n)^m)^{n-1-k}=\P(d_1=k)} are the Bin$(n-1,(2r_n)^m)$ probabilities, with \Line{p_0=\P(d_1=0)=(1-(2r_n)^m)^{n-1}} and $\sum_kk \P(d_1=k)=\Dbar$. With $x=\Dbar/(1-(2r_n)^m)$,
\[\begin{aligned}\E[\mathcal{F}_n]&=\kappa(r_n)\sum_{k\ge1}\P(d_1=k)\Big[1+x-\frac{k}{1-(2r_n)^m}\Big]\\&=\kappa(r_n)\Big\{(1-p_0)(1+x)-x\Big\}\\&=\kappa(r_n)\Big\{1-(1-(2r_n)^m)^{n-1}\Big(1+\frac{\Dbar}{1-(2r_n)^m}\Big)\Big\}.\end{aligned}\]
\hfill$\square$
\section{Proofs for the dense and intermediate regimes}\label{sec:proofs}
This section proves Theorems~\ref{thm:var} and~\ref{thm:be}; its Steps~1--5 are also used for Theorem~\ref{thm:fixed} in Section~\ref{sec:fixedproofs}. The proofs follow the same scheme in four stages:
\begin{enumerate}[leftmargin=*]
\item[(1)] \emph{An exact formula.} $\mathcal{F}_n$ equals a sum over edges of $(d_i-d_j)^2/(d_id_j)\ge0$ (Proposition~\ref{prop:id} in Section~\ref{ts:prop:id}). If all degrees were exactly $\Dbar$ we could replace $d_id_j$ by $\Dbar^2$ and obtain the simple polynomial statistic $\Psi_n/(n\Dbar^2)$, $\Psi_n=\sum_{i<j}A_{ij}(d_i-d_j)^2$; the error is shown to be small in Step~5.
\item[(2)] \emph{Counting.} $\Psi_n$ is a combination of counts of small subgraphs: a part $\mathcal{S}$ built from degrees and triangles, and a 4-point part $\Upsilon_n$.
\item[(3)] \emph{Orthogonal decomposition.} Split $\Upsilon_n$ into pieces that depend on two, three and four of the random points (Hoeffding decomposition). Because the graph is translation invariant, the pieces depending on a single point vanish. The two-point piece is a sum $\sum_{k<l}g(X_k,X_l)$ with $\E[g(X_k,X_l)]=0$: a \emph{degenerate $U$-statistic}. It is the only piece that matters; all others have variance $O(1/n)$, much smaller than $r_n^m$.
\item[(4)] \emph{Normal approximation} of the two-point piece by an exchangeable-pair argument (Step~6) with a sharp rate, and a perturbation argument to add back the small pieces.
\end{enumerate}
\noindent Steps are numbered so that the proofs of the theorems can refer to them. Throughout Steps~1--5 and the proof of Theorem~\ref{thm:var} (Section~\ref{ts:thm:var}) we fix $r_{\max}\in(0,\frac12)$ and assume $0<r_n\le r_{\max}$ and $nr_n^m\ge1$; constants may depend on $m$ and $r_{\max}$. This covers the radii of Theorems~\ref{thm:var} and~\ref{thm:be} ($r_n\to0$) as well as the constant radius of Theorem~\ref{thm:fixed} (Section~\ref{ts:thm:fixed}); in Step~6 and in the proof of Theorem~\ref{thm:be} (Section~\ref{sec:proofbe}) we assume in addition $r_n<\frac16$ (small radii). Hence \Line{\Dbar=(n-1)(2r_n)^m\ge\frac n2(2r_n)^m\ge2^{m-1}\ge1,} so Lemma~\ref{lem:mom} (Section~\ref{ts:lem:mom}) applies with $c_*=1$, and $n(2r_n)^m\ge c$ with $c=2^m$. 

\subsection{Notation for the proofs}\label{sec:notation2}
\emph{Norms.} For $x\in\T^m$, $\|x\|_\infty=d(x,0)$. For a bounded function $f$ on $\T^m$ or $\T^m\times\T^m$ (such as the pair kernel $g$), $\|f\|_\infty=\sup|f|$; for a random variable $Y$ and $p\ge1$, $\|Y\|_{L^p}=(\E|Y|^p)^{1/p}$. We write $\sigma_m^2=2\big[(\frac23)^m-2(\frac{115}{192})^m+(\frac{11}{20})^m\big]$, $\kappa_m=1-(\frac34)^m$, and $\mu_n=(n-1)(2r_n)^m$.

\paragraph{Edges along a tree.} The indicators $A_{ij}$ of different pairs are not independent in general (because of triangles), but they are independent along any tree:
\begin{lemma}[Trees]\label{lem:tree}\tlsec{lem:tree}
Let $\mathcal T$ be a tree whose vertex set is a set $V\subseteq\{1,\dots,n\}$ of $v$ vertices and whose $v-1$ edges are pairs $\{u,w\}\subset V$. Then, exactly,
\Line{\P\big(A_{uw}=1\ \text{for all edges }\{u,w\}\text{ of }\mathcal T\big)=(2r_n)^{m(v-1)}.}
\end{lemma}
\noindent The proof of Lemma~\ref{lem:tree} is given in Section~\ref{pf:lem:tree} of the Supplement.

\paragraph{The pair kernel.} For distinct $i,j,k,l$ and $x,y\in\T^m$ put
\Line{P_2(x,y)=\E[A_{ki}A_{il}\mid X_k=x,X_l=y],\qquad P_3(x,y)=\E[A_{ki}A_{ij}A_{jl}\mid X_k=x,X_l=y],}
the probability that a further point is adjacent to both $x$ and $y$, respectively that two further points form a path with three edges from $x$ to $y$ (they do not depend on the labels $i,j$). Since $\E[A_{ij}]=(2r_n)^m$ and $A_{ij}$ is independent of the other factors for a new point $j$, we have \[(2r_n)^mP_2(X_k,X_l)=\E[A_{ki}A_{il}A_{ij}\mid X_k,X_l].\]
Define
\LineE{eq:gdef}{g(x,y)=2\big((2r_n)^mP_2(x,y)-P_3(x,y)\big)=2 \E[A_{ki}A_{ij}(A_{il}-A_{jl})\mid X_k=x,X_l=y]}
and
\Line{G(x,y)=\E[g(x,X_i)g(X_i,y)].}
The pair statistic and its variance are
\LineE{eq:Wdef}{W_n=\sum_{1\le k<l\le n}g(X_k,X_l),\qquad\varsigma_n^2=\Var [W_n].}

\paragraph{Translation invariance.} Let $\beta(x)=\1\{d(0,x)<r_n\}$ be the indicator of $\mathcal A(0)$, so that $A_{uv}=\beta(X_u-X_v)$. Since the uniform law on $\T^m$ is invariant under shifts, $P_2,P_3$ and $g$ are invariant under common shifts of $(x,y)$, so they depend on $(x,y)$ only through $x-y$. We write $P_2(x)=P_2(x,0)$, $P_3(x)=P_3(x,0)$ and $g(x)=g(x,0)$, and similarly for $G$ and $h$. Explicitly, with $U,U_1,U_2$ independent and uniform on $\T^m$ and $\beta$ even,
\Line{\begin{gathered}P_2(x)=\E[\beta(x-U_1)\beta(U_1)]=|\mathcal A(0)\cap\mathcal A(x)|,\\ P_3(x)=\E[\beta(x-U_1)\beta(U_1-U_2)\beta(U_2)],\end{gathered}}
and $\E[A_{12}A_{13}A_{23}]=\E[\beta(U) P_2(U)]$. By Lemma~\ref{lem:tree} (Section~\ref{ts:lem:tree}) (the paths $k-i-l$ and $k-i-j-l$ are trees),
\Line{\E[P_2(X_k,X_l)]=(2r_n)^{2m},\qquad\E[P_3(X_k,X_l)]=(2r_n)^{3m},\qquad\E[g(X_k,X_l)]=0,}
and for all $x$,
\Line{0\le P_2(x)\le (2r_n)^m,\qquad0\le P_3(x)\le (2r_n)^{2m},\qquad P_2(0)=(2r_n)^m.}
Moreover, $P_2(x,y)=P_2(x-y,0)$ and $P_3(x,y)=P_3(x-y,0)$. Therefore, $g(x,y)=g(x-y,0)$.

\subsection{Step 1: decomposition of $\mathcal{F}_n$}\label{decomf}

\begin{proposition}[Identity]\label{prop:id}\tlsec{prop:id}
In every graph (with $\Delta_i=0$ for isolated $i$),\Line{\mathcal{F}_n=\frac1n\sum_{i<j}A_{ij}\frac{(d_i-d_j)^2}{d_id_j}\ge0 .}
\end{proposition}
\noindent The proof of Proposition~\ref{prop:id} is given in Section~\ref{pf:prop:id} of the Supplement.

\paragraph{Replacing the random denominator by a constant.} Put
\LineE{eq:Ndef}{\Psi_n=\sum_{i<j}A_{ij}(d_i-d_j)^2,\qquad \widehat{\mathcal{F}}_n=\frac{\Psi_n}{n\Dbar^2},\qquad R_n=\mathcal{F}_n-\widehat{\mathcal{F}}_n .}
So
\Line{R_n=\frac1n\sum_{i<j}A_{ij}(d_i-d_j)^2\big(\frac1{d_id_j}-\frac1{\Dbar^2}\big);}
it is small because $d_id_j\approx\Dbar^2$.

\paragraph{Expanding the square.} Fix an edge $\{i,j\}$ ($A_{ij}=1$). Since $d_i=A_{ij}+\sum_{k\ne i,j}A_{ik}$ and likewise for $d_j$, and the two copies of $A_{ij}$ cancel,
\Line{d_i-d_j=\sum_{k\ne i,j}(A_{ik}-A_{jk}),\qquad A_{ik}-A_{jk}\in\{-1,0,1\}.}
Squaring, \Line{(d_i-d_j)^2=\sum_{k\notin\{i,j\}}(A_{ik}-A_{jk})^2+\sum_{k\ne l, k,l\notin\{i,j\}}(A_{ik}-A_{jk})(A_{il}-A_{jl}).}
Hence
\Line{\Psi_n=\mathcal{S}+\Upsilon_n}
with
\LineE{eq:SC}{\mathcal{S}=\sum_{i<j}\sum_{k\ne i,j}A_{ij}(A_{ik}-A_{jk})^2,\qquad \Upsilon_n=\sum_{i<j}\sum_{\substack{(k,l):\ k\ne l\\ k,l\notin\{i,j\}}}A_{ij}(A_{ik}-A_{jk})(A_{il}-A_{jl}).}
\begin{lemma}\label{lem:S}\tlsec{lem:S}
$\mathcal{S}=\sum_id_i(d_i-1)-6\mathcal T$, where $\mathcal T$ is the number of triangles of the graph.
\end{lemma}
\noindent The proof of Lemma~\ref{lem:S} is given in Section~\ref{pf:lem:S} of the Supplement.

\begin{lemma}\label{lem:ECn}\tlsec{lem:ECn}
$\E [\Upsilon_n]=0$, and every summand \Line{H_{ij;kl}=A_{ij}(A_{ik}-A_{jk})(A_{il}-A_{jl})} of $\Upsilon_n$ satisfies $\E[H_{ij;kl}\mid X_U]=0$ whenever $k\notin U$ or $l\notin U$, $U\subseteq\{i,j,k,l\}$.
\end{lemma}
\noindent The proof of Lemma~\ref{lem:ECn} is given in Section~\ref{pf:lem:ECn} of the Supplement.

 The Hoeffding decomposition writes a function of independent random variables as a sum of mutually orthogonal pieces, each depending on a subset of the variables, so that variances add.  We need it to separate the two-point component of $\Upsilon_n$, which carries the main fluctuations of $\mathcal{F}_n$, from the three- and four-point components, whose variances are small and bounded separately; by translation invariance the one-point component vanishes. 
\begin{tool}[Hoeffding decomposition \cite{Hoeffding48}, see \cite{Lee90}]\label{tl:D}\tlsec{D}
Let $X_1,\dots,X_n$ be independent and $f=\sum_{\mathbf i}h_{\mathbf i}(X_{\mathbf i})$ a finite sum of square-integrable terms, each depending on the coordinates in an index set $\mathbf i\subseteq\{1,\dots,n\}$. For every $S\subseteq\{1,\dots,n\}$ define
\Line{\pi_S f=\sum_{U\subseteq S}(-1)^{|S|-|U|} \E\big[f \big| X_U\big]\qquad(\E[f|X_\emptyset]=\E [f]).}
Then (a) $f=\sum_{S}\pi_Sf$; (b) $\pi_Sf$ is a function of $X_S$ and the $\pi_Sf$ are pairwise orthogonal in $L^2$, hence $\Var [f]=\sum_{S\ne\emptyset}\E[(\pi_Sf)^2]$; (c) $\E[\pi_Sf\mid X_U]=0$ unless $S\subseteq U$; (d) $\pi_Sf=\pi_S\big(\E[f|X_S]\big)$ is an orthogonal projection of $\E[f|X_S]$, so $\|\pi_Sf\|_{L^2}\le\|\E[f|X_S]\|_{L^2}$ (for $f$ with $\E[f|X_U]=0$ for all $U\subsetneq S$ one has $\pi_Sf=\E[f|X_S]$).
\end{tool}

\noindent\textbf{Which components exist.} Write $\Upsilon_n=\sum H_{ij;kl}$, where each summand is \Line{H_{ij;kl}=A_{ij}(A_{ik}-A_{jk})(A_{il}-A_{jl})} with $A_{uv}=\beta(X_u-X_v)$. There are $\binom n2$ choices of the unordered pair $\{i,j\}$ and $(n-2)(n-3)$ choices of the ordered pair $(k,l)$, so $12\binom n4$ summands; each is attached to the $4$-set $V=\{i,j,k,l\}$, and each $4$-set carries exactly $12$ summands ($6$ choices of $\{i,j\}$ times $2$ orders of $(k,l)$). We apply the Hoeffding decomposition in Theorem~\ref{tl:D}  to one summand $H$ and then add up over summands.

Since $H_{ij;kl}$ depends only on $X_V=\{X_l,l\in V\}$ with $V=\{i,j,k,l\}$, $\E[H_{ij;kl}\mid X_U]=\E[H_{ij;kl}\mid X_{U\cap V}]$ for every $U$; if $a\in S$ and $a\notin V$, pairing $U\subseteq S\setminus\{a\}$ with $U\cup\{a\}$ cancels all terms of $\pi_SH_{ij;kl}$, so $\pi_SH_{ij;kl}=0$ unless $S\subseteq V$. By Lemma~\ref{lem:ECn}, $\E[H_{ij;kl}\mid X_U]=0$ whenever $U$ misses $k$ or $l$, so $\pi_SH_{ij;kl}=0$ unless $\{k,l\}\subseteq S$. Hence the only non-zero components are those on $\{k,l\}$, $\{i,k,l\}$, $\{j,k,l\}$, $\{i,j,k,l\}$, namely
\[\begin{aligned}
\pi_{\{k,l\}}H_{ij;kl}&=\E[H\mid X_k,X_l],\\
\pi_{\{i,k,l\}}H_{ij;kl}&=\E[H_{ij;kl}\mid X_i,X_k,X_l]-\E[H_{ij;kl}\mid X_k,X_l],\\
\pi_{\{j,k,l\}}H_{ij;kl}&=\E[H_{ij;kl}\mid X_j,X_k,X_l]-\E[H_{ij;kl}\mid X_k,X_l],\\
\pi_{\{i,j,k,l\}}H_{ij;kl}&=H_{ij;kl}-\E[H_{ij;kl}\mid X_i,X_k,X_l]-\E[H_{ij;kl}\mid X_j,X_k,X_l]+\E[H_{ij;kl}\mid X_k,X_l].
\end{aligned}\]
In particular there are no one-point components: $\Upsilon_n$ has no linear part (translation invariance). Summing over all summands according to the size of the index set $S\subset\{1,2,\dots,n\}$ gives the three components
\[\Upsilon_n^{(2)}=\sum_{\substack{|S|=2}}\pi_{S}\Upsilon_n,\ \ \ \  \ \Upsilon_n^{(3)}=\sum_{\substack{|S|=3}}\pi_{S}\Upsilon_n,\ \ \ \  \ \Upsilon_n^{(4)}=\sum_{\substack{|S|=4}}\pi_{S}\Upsilon_n,\]
so that \Line{\Upsilon_n=\Upsilon_n^{(2)}+\Upsilon_n^{(3)}+\Upsilon_n^{(4)}.}  Then
\LineE{eq:split}{\mathcal{F}_n-\E[\mathcal{F}_n]=\underbrace{\tfrac{1}{n\Dbar^2}\Upsilon_n^{(2)}}_{M_n}+\underbrace{\tfrac{1}{n\Dbar^2}\big(\Upsilon_n-\Upsilon_n^{(2)}\big)}_{\mathrm{Rem}^{(1)}}+\underbrace{\tfrac{1}{n\Dbar^2}(\mathcal{S}-\E [\mathcal{S}])}_{\mathrm{Rem}^{(2)}}+\underbrace{(R_n-\E [R_n])}_{\mathrm{Rem}^{(3)}},}
where we used $\mathcal{F}_n=\frac{\mathcal{S}+\Upsilon_n}{n\Dbar^2}+R_n$ and $\E [\Upsilon_n]=0$.  $M_n$ 
is the main term, built from the two-point component only, and $\mathrm{Rem}^{(1)}$ (the three- and four-point components of $\Upsilon_n$), $\mathrm{Rem}^{(2)}$ (the centred degree-and-triangle part) and $\mathrm{Rem}^{(3)}$ (the centred denominator error) form the remainder \Line{\mathrm{Rem}=\mathrm{Rem}^{(1)}+\mathrm{Rem}^{(2)}+\mathrm{Rem}^{(3)},} so that 
\[\mathcal{F}_n-\E[\mathcal{F}_n]=M_n+\mathrm{Rem}.\]

\subsection{Step 2: the second Hoeffding component of $\Upsilon_n$}\label{sec:step2}

\begin{proposition}[The constants $\alpha_1,\alpha_2,\alpha_3$]\label{prop:alpha}\tlsec{prop:alpha}
Let $0<r<\frac12$, $r_n=r$, and let $X_k,X_l$ be independent and uniform on $\T^m$, with $P_2=P_2(X_k,X_l)$ and $P_3=P_3(X_k,X_l)$. Let $\alpha_1,\alpha_2,\alpha_3$ be the piecewise polynomials~\eqref{eq:alpha1}--\eqref{eq:alpha3}. Then
\[\E[P_2^2]=\alpha_1(r)^m,\ \ \ \ \ \ \ \ \E[P_2P_3]=\alpha_2(r)^m,\ \ \ \ \ \ \ \ \ \E[P_3^2]=\alpha_3(r)^m,\]
and
\begin{align*}
\E[g(X_k,X_l)^2]=4\big[(2r)^{2m}\alpha_1(r)^m-2(2r)^m\alpha_2(r)^m+\alpha_3(r)^m\big].
\end{align*}
\end{proposition}
\noindent The proof of Proposition~\ref{prop:alpha} is given in Section~\ref{pf:prop:alpha} of the Supplement.

\medskip\noindent\textbf{(a) The pair component.} 
Expand the product  in $H_{ij;kl}$ and take conditional expectations given $X_k,X_l$:
\[\begin{aligned}\E[H_{ij;kl}\mid X_k,X_l]&=\E[A_{ij}A_{ik}A_{il}\mid X_k,X_l]-\E[A_{ij}A_{ik}A_{jl}\mid X_k,X_l]\\
&\quad-\E[A_{ij}A_{jk}A_{il}\mid X_k,X_l]+\E[A_{ij}A_{jk}A_{jl}\mid X_k,X_l].\end{aligned}\]
In the first term, condition also on $X_i$: the edge $A_{ij}$ involves the new point $X_j$ only, and $\E[A_{ij}\mid X_i]=(2r_n)^m$, so the term equals \Line{(2r_n)^m \E[A_{ik}A_{il}\mid X_k,X_l]=(2r_n)^m P_2(X_k,X_l);} the fourth term is the same by symmetry of $i$ and $j$. In the second term the edges $\{k,i\},\{i,j\},\{j,l\}$ form a path from $X_k$ to $X_l$, so the term is $P_3(X_k,X_l)$; the third term is the path $k,j,i,l$ and equals $P_3(X_k,X_l)$ as well. Therefore
\begin{equation}\begin{aligned}
\E[H_{ij;kl}\mid X_k,X_l]&=2 \big((2r_n)^m P_2(X_k,X_l)-P_3(X_k,X_l)\big)\\&=2 \E[A_{ki}A_{ij}(A_{il}-A_{jl})\mid X_k,X_l]\\
&=g(X_k,X_l).\end{aligned}\end{equation}
A given unordered pair $\{k,l\}$ occurs as the roles $(k,l)$ or $(l,k)$ (two choices), and the pair $\{i,j\}$ can be any of the $\binom{n-2}2$ pairs among the remaining $n-2$ points; so each $\{k,l\}$ occurs in $2\binom{n-2}2=(n-2)(n-3)$ summands, each with the same pair component $g(X_k,X_l)$. Hence, with $W_n$ as in~\eqref{eq:Wdef}, \Line{\Upsilon_n^{(2)}=(n-2)(n-3) W_n,} the sum of the two-point components of $\Upsilon_n$.

\medskip\noindent\textbf{(b) Properties of $g$.}
Write $P_2=P_2(X_k,X_l)$, $P_3=P_3(X_k,X_l)$, so that \[g=2\big((2r_n)^mP_2-P_3\big).\]

\emph{The exact value of $\E[g(X_k,X_l)^2]$.}  Note that
\[\E[g(X_k,X_l)^2]=4\big((2r_n)^{2m}\E[P_2^2]-2(2r_n)^m\E[P_2P_3]+\E[P_3^2]\big).\]
 By Proposition~\ref{prop:alpha}   (valid for every $r_n<\frac12$), 
\LineE{eq:aexact}{\omega:=\E[g(X_k,X_l)^2]=4\big[(2r_n)^{2m}\alpha_1^m-2(2r_n)^m\alpha_2^m+\alpha_3^m\big]=2\sigma^2(r_n) (2r_n)^{4m}>0 .}
Here $\sigma^2(r_n)>0$ (defined in~\eqref{eq:sigmaformula}) for every $r_n<\frac12$ , and for $r_n<\frac16$,  $\sigma^2(r_n)=\sigma_m^2(2r_n)^m$, i.e.\ $\omega=2\sigma_m^2(2r_n)^{5m}$.

\medskip\noindent\textbf{(c) Variance of the pair component.} For two different pairs $\{k,l\}\ne\{k',l'\}$ the variables $g(X_k,X_l)$ and $g(X_{k'},X_{l'})$ are uncorrelated: if the pairs are disjoint they are independent with mean $0$; if they share one index, say $k=k'$, then
\Line{\E[g(X_k,X_l)g(X_k,X_{l'})]=\E\big[\E[g(X_k,X_l)|X_k] \E[g(X_k,X_{l'})|X_k]\big]=0}
because \Line{\E[g(X_k,X_l)\mid X_k]=\E\big[\E[H_{ij;kl}\mid X_k,X_l]\mid X_k\big]=\E[H_{ij;kl}\mid X_k]=0} by Lemma~\ref{lem:ECn} (Section~\ref{ts:lem:ECn}). Hence $\Var [W_n]=\binom n2 \E[g(X_k,X_l)^2]$ and, with $M_n=\frac{(n-2)(n-3)}{n\Dbar^2}W_n$ and $\Dbar=(n-1)(2r_n)^m$,
\begin{equation}\label{eq:varM}\begin{aligned}\Var [M_n]&=\frac{(n-2)^2(n-3)^2}{n^2(n-1)^4(2r_n)^{4m}}\binom n2 \omega=\frac{(n-2)^2(n-3)^2}{n(n-1)^3}\cdot\frac{\omega}{2(2r_n)^{4m}}\\
&=\frac{(n-2)^2(n-3)^2}{n(n-1)^3} \sigma^2(r_n)=\sigma^2(r_n) \Big(1+O\big(\tfrac1n\big)\Big).\end{aligned}\end{equation}

\subsection{Step 3: the three- and four-point components of $\Upsilon_n$}

\smallskip\noindent\emph{Triples.} Fix $T=\{i,k,l\}$ (the case $\{j,k,l\}$ is identical by symmetry). The summand is $H_{ij;kl}$,  and $X_T=(X_i,X_k,X_l)$ determines everything in it except $X_j$. We bound $|H_{ij;kl}|$ pointwise by a product of an indicator that is measurable with respect to $X_T$ and the single factor $A_{ij}$, and then take the conditional expectation, exactly as for $g(X_k,X_l)=\E[H_{ij;kl}\mid X_k,X_l]$.

\emph{Pointwise bound.} Both factors $A_{ik}-A_{jk}$ and $A_{il}-A_{jl}$ lie in $\{-1,0,1\}$. Suppose $A_{ij}=1$, i.e.\ $d(X_i,X_j)<r_n$. If $A_{ik}-A_{jk}\ne0$, then $A_{ik}=1$ or $A_{jk}=1$, so $d(X_k,X_i)<r_n$ or $d(X_k,X_j)<r_n$; in the second case \Line{d(X_k,X_i)\le d(X_k,X_j)+d(X_j,X_i)<2r_n} by the triangle inequality. In both cases $d(X_k,X_i)<2r_n$, and the same argument gives $d(X_l,X_i)<2r_n$ when $A_{il}-A_{jl}\ne0$. Hence
\[|H_{ij;kl}|\le A_{ij} \1\{d(X_k,X_i)<2r_n,\ d(X_l,X_i)<2r_n\}.\]

\emph{Conditional expectation.} The indicator is a function of $X_T$, so it can be taken out of the conditional expectation, and \Line{\E[A_{ij}\mid X_T]=\E[A_{ij}\mid X_i]=(2r_n)^m} because $X_j$ is independent of $(X_k,X_l)$ given $X_i$. Therefore
\[\begin{aligned}\big|\E[H_{ij;kl}\mid X_T]\big|&\le\E\big[|H_{ij;kl}|\ \big|\ X_T\big]\\&\le\E[A_{ij}\mid X_T]\;\1\{d(X_k,X_i)<2r_n,\ d(X_l,X_i)<2r_n\}\\&=(2r_n)^m \1\{d(X_k,X_i)<2r_n,\ d(X_l,X_i)<2r_n\}.\end{aligned}\]
Taking squares and using independence of $X_k,X_l$ given $X_i$,
\Line{\E\big[\big(\E[H_{ij;kl}|X_T]\big)^2\big]\le (2r_n)^{2m} \big((4r_n)^m\big)^2\le C(2r_n)^{4m} .}
By Theorem~\ref{tl:D}(d) (Section~\ref{ts:D}), \Line{\E[(\pi_TH_{ij;kl})^2]\le\E[\E[H_{ij;kl}|X_T]^2]\le C(2r_n)^{4m}.} A fixed triple $T$ lies in $12(n-3)$ of the $4$-sets' summands (choose the fourth index; $12$ summands per $4$-set; those whose role pattern does not contain $\{k,l\}\subset T$ have zero component and may be ignored, which only helps). By the Cauchy--Schwarz inequality for a sum of $12(n-3)$ terms, \Line{\E[(\pi_T\Upsilon_n)^2]\le(12(n-3))^2C(2r_n)^{4m}\le Cn^2(2r_n)^{4m},} and by orthogonality of components on different sets,
\Line{\Var\big[\Upsilon_n^{(3)}\big]=\sum_{|T|=3}\E[(\pi_T\Upsilon_n)^2]\le\binom n3 Cn^2(2r_n)^{4m}\le Cn^5(2r_n)^{4m} .}
\emph{Quadruples.} The argument has four steps: (1) only the summands living on $U$ contribute to $\pi_U\Upsilon_n$; (2) the projection does not increase the norm; (3) the second moment of one summand is $\le4(2r_n)^{3m}$; (4) different $U$ are orthogonal, so the variance is a sum over $U$.

\smallskip\noindent(1) Fix a $4$-set $U$. By Section \ref{decomf}, every summand $H_{ij;kl}$ with $\{i,j,k,l\}\ne U$ has $\pi_UH_{ij;kl}=0$. The summands with $\{i,j,k,l\}=U$ are the $12$ summands (listed in Section \ref{decomf}): $6$ choices of $\{i,j\}\subset U$ times $2$ orders of the remaining pair $(k,l)$. Hence $\pi_U\Upsilon_n=\pi_U\Sigma_U$ with $\Sigma_U=\sum_{12\text{ summands}}H_{ij;kl}$.

\smallskip\noindent(2) By Theorem~\ref{tl:D}(d) (Section~\ref{ts:D}), $\pi_U$ is an orthogonal projection in $L^2$, so by Minkowski's inequality
\[\E[(\pi_U\Upsilon_n)^2]\le\E[\Sigma_U^2]\le12\sum_{12\text{ terms}}\E[H_{ij;kl}^2]=144 \E[H_{ij;kl}^2],\]
because the $12$ summands have the same law (the points are i.i.d.\ and the summands differ by relabelling).

\smallskip\noindent(3) With $H_{ij;kl}$ as in Step~2(a) of Section~\ref{sec:proofs}, since $A_{ij}\in\{0,1\}$, \Line{H_{ij;kl}^2=A_{ij}(A_{ik}-A_{jk})^2(A_{il}-A_{jl})^2.} For $\{0,1\}$-valued variables \Line{(A_{ik}-A_{jk})^2=A_{ik}+A_{jk}-2A_{ik}A_{jk}\le A_{ik}+A_{jk}.} Given $(X_i,X_j)$, the points $X_k,X_l$ are independent of each other and of $(X_i,X_j)$, so the two factors are conditionally independent, and $\E[A_{ik}\mid X_i,X_j]=(2r_n)^m$. Therefore
\[\E[H_{ij;kl}^2\mid X_i,X_j]\le A_{ij} \big(2(2r_n)^m\big)\big(2(2r_n)^m\big)=4(2r_n)^{2m}A_{ij},\]
and taking expectations gives $\E[H_{ij;kl}^2]\le4(2r_n)^{3m}$.

\smallskip\noindent(4) By Theorem~\ref{tl:D} (Section~\ref{ts:D}), the components $\pi_U\Upsilon_n$ for different sets $U$ are orthogonal, and $\Upsilon_n^{(4)}=\sum_{|U|=4}\pi_U\Upsilon_n$. Combining (2) and (3), and $\E[\pi_U\Upsilon_n]=0$,
\Line{\Var[\Upsilon_n^{(4)}]=\sum_{|U|=4}\E[(\pi_U\Upsilon_n)^2]\le\binom n4\cdot12^2\cdot4(2r_n)^{3m}\le Cn^4(2r_n)^{3m}\le Cn^5(2r_n)^{4m},}
because $n(2r_n)^m\ge c$. Hence, since $\Upsilon_n^{(3)}$ and $\Upsilon_n^{(4)}$ are orthogonal,
\LineE{eq:rem1}{\Var\big[\mathrm{Rem}^{(1)}\big]=\frac{\Var\big[\Upsilon_n^{(3)}\big]+\Var\big[\Upsilon_n^{(4)}\big]}{n^2\Dbar^4}\le\frac{Cn^5(2r_n)^{4m}}{n^2 (n-1)^4(2r_n)^{4m}}\le\frac Cn .}

\subsection{Step 4: the degree-and-triangle part $\mathcal{S}$}

This step shows that $\mathcal{S}$, the degree-and-triangle part of $\Psi_n$ in~\eqref{eq:SC}, is negligible: $\Var[\mathrm{Rem}^{(2)}]\le C/n$ (see~\eqref{eq:rem2}). It uses one lemma and one standard tool, in this order:
\begin{itemize}[leftmargin=*]
\item Lemma~\ref{lem:mom}  bounds the moments of degrees and of neighbourhood sizes, and is the only input about the graph that the argument needs;
\item Theorem~\ref{tl:C} (the Efron--Stein inequality) bounds $\Var[\mathcal{S}]$ by the effect of re-sampling one point at a time. It is applied, together with Lemma~\ref{lem:mom}, in the proof of~\eqref{eq:rem2} that follows.
\end{itemize}

\begin{lemma}[Moment bounds]\label{lem:mom}\tlsec{lem:mom}
Let \Line{\nu_k=\#\{j\ne k:d(X_j,X_k)<2r_n\}} be the number of points within distance $2r_n$ of $X_k$. Assume $\Dbar\ge c_*$ for a constant $c_*>0$. For every $p\ge1$ and a constant $C_p$ dependent on $p$:
\begin{enumerate}[label=(\roman*),leftmargin=*]
\item $\|d_i-\Dbar\|_{L^p}\le C_p\Dbar^{1/2}$ and $\|d_i\|_{L^p}\le C_p\Dbar$;
\item $\|\nu_k\|_{L^p}\le C_p\Dbar$, and the same bound holds for \Line{\tilde\nu_k=\#\{j\ne k:d(X_j,X_k)<3r_n\};}
\item $d_k\le \nu_k$, and every neighbour $l$ of $k$ has $d_l\le \nu_k+1$; hence $\sum_{l:A_{kl}=1}d_l\le \nu_k(\nu_k+1)$;
\item the bounds in (i)--(ii) remain true, with the same constants, for the conditional laws of $d_v$ and $\nu_k$ given the positions of any fixed finite set of other points, provided the number of conditioned points is bounded by a constant (the bounds then hold with $C_p$ depending also on that number).
\end{enumerate}
\end{lemma}
\noindent The proof of Lemma~\ref{lem:mom}, which rests on the Rosenthal moment inequality for sums of independent variables \cite{Rosenthal70}, is given in Section~\ref{pf:lem:mom} of the Supplement.

The Efron--Stein inequality controls the variance of a function of independent random elements by the average squared change of that function when a single coordinate is re-sampled. It is therefore the natural tool whenever moving one point changes the statistic only slightly. We use it for $\mathcal{S}$ in this step, in Step~5, in the sparse upper bound (Section~\ref{sec:sparseproofs}) and in Lemma~\ref{lem:VarD} (Section~\ref{ts:lem:VarD}).
\begin{tool}[Efron--Stein inequality \cite{EfronStein81,Steele86}]\label{tl:C}\tlsec{C}
Let $f$ be a square-integrable function of independent random elements $X_1,\dots,X_n$, let $X_k'$ be an independent copy of $X_k$ and \Line{f^{(k)}=f(X_1,\dots,X_{k-1},X_k',X_{k+1},\dots,X_n).} Then
\Line{\Var [f]\le\frac12\sum_{k=1}^n\E\big[\big(f-f^{(k)}\big)^2\big] .}
\end{tool}

\medskip\noindent\textbf{Proof of~\eqref{eq:rem2}.} The proof applies Theorem~\ref{tl:C}  to $\mathcal{S}$ and bounds the effect of moving one point by Lemma~\ref{lem:mom}.

\smallskip\noindent\textbf{Setup.} By Lemma~\ref{lem:S} (Section~\ref{ts:lem:S}), $\mathcal{S}=\sum_{i}d_i(d_i-1)-6\mathcal T$, a function of the independent points $X_1,\dots,X_n$. Fix $k$ and let $X_k'$ be an independent copy of $X_k$. Write $\mathcal{S}^{(k)}$ for $\mathcal{S}$ computed with $X_k$ replaced by $X_k'$, and use primes for the quantities of the new graph: $d_l'$ is the degree of $l$ in it, $N(k)$ and $N'(k)$ are the neighbourhoods of $k$ in the old and the new graph, $d_k=|N(k)|$, $d_k'=|N'(k)|$, and $\nu_k$, $\nu_k'$ count the other points within distance $2r_n$ of $X_k$, respectively $X_k'$. Put $\nu=\nu_k+\nu_k'$. The goal is the pointwise bound 
\[|\mathcal{S}-\mathcal{S}^{(k)}|\le C(\nu+1)^2.\]

\medskip\noindent\textbf{(a) What changes.} The edges not containing $k$ are the same in both graphs, because they depend only on the other points. So (i) the degree of $k$ changes from $d_k$ to $d_k'$; (ii) for $l\ne k$, $d_l'=d_l-A_{kl}+A_{kl}'$, so $d_l$ changes (by exactly $\pm1$) only if $l\in N(k) \triangle N'(k)$, and the number of such vertices is at most $d_k+d_k'$; (iii) only the triangles through $k$ can be created or destroyed.

\medskip\noindent\textbf{(b) The vertex $k$.} Since $0\le d(d-1)\le d^2$ for nonnegative integer $d$, the term $d_k(d_k-1)$ is replaced by $d_k'(d_k'-1)$, a change of at most $d_k^2+d_k'^2$.

\medskip\noindent\textbf{(c) The other vertices.} Let $f(d)=d(d-1)$, so that 
\[f(d+1)-f(d)=2d,\ \ \ \  \ \ \ f(d)-f(d-1)=2(d-1).\]
A vertex $l\ne k$ whose degree changes by $\pm1$ therefore changes $f(d_l)$ by at most $2\max(d_l,d_l')$. We bound this by $2(\nu+1)$. If $l\in N(k)$, Lemma~\ref{lem:mom}(iii)  gives $d_l\le\nu_k+1$, and $d_l'\le d_l$ because a vertex of $N(k)$ loses the edge $kl$ and gains at most the edge $kl$ back, so $d_l'\le d_l$. If $l\in N'(k)\setminus N(k)$, Lemma~\ref{lem:mom}(iii) (Section~\ref{ts:lem:mom}) applied to the new graph gives $d_l'\le\nu_k'+1$, and $d_l=d_l'-1<d_l'$. In both cases $\max(d_l,d_l')\le\nu+1$. Since at most $d_k+d_k'\le\nu_k+\nu_k'=\nu$ vertices change (we used $d_k\le\nu_k$, $d_k'\le\nu_k'$ from Lemma~\ref{lem:mom}(iii)), the total change over $l\ne k$ is at most $2\nu(\nu+1)$.

\medskip\noindent\textbf{(d) The triangles.} Triangles through $k$ are determined by pairs of neighbours of $k$, so there are at most $\binom{d_k}2$ of them in the old graph and $\binom{d_k'}2$ in the new one, and the change in $\mathcal T$ is at most $\max\{\binom{d_k}2,\binom{d_k'}2\}\le\frac12\nu^2$. Hence $6 |\mathcal T-\mathcal T^{(k)}|\le3\nu^2$.

\medskip\noindent\textbf{(e) Adding up.} By (b)--(d), using $d_k^2+d_k'^2\le\nu_k^2+\nu_k'^2\le\nu^2$,
\Line{|\mathcal{S}-\mathcal{S}^{(k)}|\le\nu^2+2\nu(\nu+1)+3\nu^2\le 8(\nu+1)^2 .}

\medskip\noindent\textbf{(f) Moments and Efron--Stein.} By Minkowski and Lemma~\ref{lem:mom}(ii) with $p=4$ (both $\nu_k$ and $\nu_k'$ are distributed as $\nu_k$), and $\Dbar\ge c_*$,
\[\E[(\nu+1)^4]\le\big(2C_4\Dbar+1\big)^4\le C\Dbar^4,\]
hence \Line{\E[(\mathcal{S}-\mathcal{S}^{(k)})^2]\le C \E[(\nu+1)^4]\le C\Dbar^4} for every $k$. By the Efron--Stein inequality (Theorem~\ref{tl:C}), summing over the $n$ points,
\[\Var[\mathcal{S}]\le\frac12\sum_{k=1}^n\E\big[(\mathcal{S}-\mathcal{S}^{(k)})^2\big]\le\frac12 n C^2\Dbar^4 .\]
Dividing by the normalisation $n^2\Dbar^4$ of $\mathrm{Rem}^{(2)}$,
\LineE{eq:rem2}{\Var\big[\mathrm{Rem}^{(2)}\big]=\frac{\Var [\mathcal{S}]}{n^2\Dbar^4}\le\frac{Cn\Dbar^4}{n^2\Dbar^4}=\frac Cn .}

\subsection{Step 5: the remainder $R_n$ from replacing $d_id_j$ by $\Dbar^2$}

\noindent The remainder $R_n$ is a sum over edges of the terms $\psi(d_i,d_j)$ defined below, which compare $1/(d_id_j)$ with $1/\Dbar^2$. Such a term is small only when both degrees are comparable to $\Dbar$; for a vertex of unusually small degree it is not. The proof therefore needs to know that atypical degrees are rare. Degrees are binomial~\eqref{eq:degbin}, and the Chernoff bound shows that a binomial variable is within a constant factor of its mean except with exponentially small probability. It is used below to show that all degrees are ``typical'' when $\Dbar$ is large, and, together with Lemma~\ref{lem:mom} (Section~\ref{ts:lem:mom}) for the conditional case, that atypical vertices contribute negligibly to the variance.
\begin{tool}[Chernoff bound; Chernoff 1952, see e.g.\ \cite{MU05}]\label{tl:A}\tlsec{A}
If $B\sim\mathrm{Bin}(N,p)$ and $\mu=Np$, then for $0<\varepsilon\le1$,
\Line{\P(|B-\mu|\ge\varepsilon\mu)\le2e^{-\varepsilon^2\mu/3}.}
\end{tool}

\paragraph{Typical vertices.} Assume $\Dbar\ge16$. Call a vertex $v$ \emph{typical} if $|d_v-\Dbar|\le\Dbar/4$. By Theorem~\ref{tl:A} (applied to the binomial~\eqref{eq:degbin} with $\varepsilon=\frac14$), $\P(v\text{ is atypical})\le2e^{-\Dbar/48}$. When some other points are held fixed (``conditioned''), the degree is still a sum of independent Bernoulli variables whose parameters change by at most $O((2r_n)^m)$, and the same bound holds up to constants (this is made precise in Lemma~\ref{lem:mom}(iv) in Section~\ref{ts:lem:mom}). We shall use the following consequence: for $\Dbar\ge16$,
\Line{\P(v\text{ atypical})\le Ce^{-\Dbar/60}.}

For integers $x,y\ge1$ let
\Line{\psi(x,y)=\frac{(x-y)^2}{xy}-\frac{(x-y)^2}{\Dbar^2}.}
By Proposition~\ref{prop:id} (Section~\ref{ts:prop:id}),
\Line{R_n=\frac1n\sum_{i<j}A_{ij} \psi(d_i,d_j).}
We show $\Var [R_n]\le C/n$ when $r_n\le r_{\max}$ and $nr_n^m\ge1$, assuming nothing on $\Dbar$ beyond $\Dbar\ge c_*$. By \TS{C} it suffices to prove
\Line{\E\big[(R_n-R_n^{(k)})^2\big]\le \frac{C}{n^2},}
where $R_n^{(k)}$ is obtained by replacing $X_k$ by an independent copy $X_k'$.

\medskip\noindent\textbf{Set-up.} Write $\mathcal{G},\mathcal{G}'$ for the graphs before and after the replacement, $d_i,d_i'$ for the degrees, and $N(k),N'(k)$ for the neighbourhoods of $k$. Let $\nu_k,\nu_k'$ ($\tilde\nu_k,\tilde\nu_k'$) count the points $j\ne k$ within distance $2r_n$ ($3r_n$) of $X_k$, $X_k'$, and put $\nu=\nu_k+\nu_k'\le\tilde\nu=\tilde\nu_k+\tilde\nu_k'$. Let \Line{\mathcal D_k=\{v\ne k:\ d(X_v,X_k)<2r_n\text{ or }d(X_v,X_k')<2r_n\},} so $|\mathcal D_k|\le\nu$. The triangle inequality gives:
\begin{enumerate}[label=(F\arabic*),leftmargin=*,itemsep=1pt]
\item An edge $\{i,j\}\not\ni k$ is present in $\mathcal{G}$ iff it is present in $\mathcal{G}'$. For $i\ne k$, \Line{d_i'-d_i=A_{ik}'-A_{ik}\in\{-1,0,1\},} and it is non-zero only for $i\in N(k)\triangle N'(k)\subseteq\mathcal D_k$.
\item Every neighbour $u\ne k$ of a vertex of $N(k)\cup N'(k)$ belongs to $\mathcal D_k$; and the degrees $d_k,d_k'$ and $d_v,d_v'$, $v\in\mathcal D_k$, are at most $\tilde\nu+1$ (the neighbours of $v\in\mathcal D_k$ other than $k$ lie within $3r_n$ of $X_k$ or of $X_k'$). In particular each vertex of $\mathcal D_k$ lies on at most $\tilde\nu+1$ edges.
\end{enumerate}
By (F1), splitting the edges according to whether they contain $k$,
\[n(R_n-R_n^{(k)})=\underbrace{\sum_{i\in N(k)}\psi(d_k,d_i)-\sum_{i\in N'(k)}\psi(d_k',d_i')}_{\text{(a)}}+\underbrace{\sum_{\{i,j\}\in\mathcal B_k}\big[\psi(d_i,d_j)-\psi(d_i',d_j')\big]}_{\text{(b)}},\]
where $\mathcal B_k$ is the set of edges $\{i,j\}\not\ni k$ with an endpoint in $N(k)\triangle N'(k)$ (the other terms of (b) vanish). By (F2) both endpoints of an edge of $\mathcal B_k$ lie in $\mathcal D_k$, so $|\mathcal B_k|\le\tilde\nu(\tilde\nu+1)$; term (a) has $d_k+d_k'\le\nu$ summands.

\medskip\noindent\textbf{A crude bound.} We show, for every $\Dbar\ge c_*$,
\LineE{eq:crude}{|n(R_n-R_n^{(k)})|\le C(\tilde\nu+2)^3 .}
For integers $x,y\ge1$, \[\frac{(x-y)^2}{xy}=\frac xy+\frac yx-2\le x+y,\]
so \Line{|\psi(x,y)|\le x+y+\frac{(x+y)^2}{c_*^2}\le C(x+y)^2.} Also $\partial_x\psi=\frac1y-\frac y{x^2}-\frac{2(x-y)}{\Dbar^2}$, so by the mean value theorem \Line{|\psi(x\pm1,y)-\psi(x,y)|\le C(1+x+y)} if $x,x\pm1\ge1$, and the same in $y$. Each summand of (a) is $\psi$ at integers in $[1,\tilde\nu+1]$ by (F2), hence at most $C(\tilde\nu+2)^2$, and $|\text{(a)}|\le C\tilde\nu(\tilde\nu+2)^2$. For $\{i,j\}\in\mathcal B_k$ the four degrees $d_i,d_j,d_i',d_j'$ are in $[1,\tilde\nu+1]$ and $(d_i',d_j')$ differs from $(d_i,d_j)$ by at most one unit in each coordinate, so \Line{|\psi(d_i,d_j)-\psi(d_i',d_j')|\le C(\tilde\nu+2)} (two unit steps), and $|\text{(b)}|\le C\tilde\nu(\tilde\nu+1)(\tilde\nu+2)$. This gives~\eqref{eq:crude}. 

If $\Dbar\le16$, then~\eqref{eq:crude} and Lemma~\ref{lem:mom}(ii) (Section~\ref{ts:lem:mom}) give 
\[\E[(n(R_n-R_n^{(k)}))^2]\le C.\]

\medskip\noindent\textbf{Case $\Dbar>16$: the good event.} Let $\mathcal E_k$ be the event that all vertices of $\mathcal D_k$ are typical and $|d_k-\Dbar|,|d_k'-\Dbar|\le\Dbar/4$. On $\mathcal E_k$ every degree that occurs (in $\mathcal{G}$ or $\mathcal{G}'$, also shifted by one unit) lies in $[\Dbar/2,2\Dbar]$. For $x,y\in[\Dbar/2,2\Dbar]$ (and $x\pm1$ in the same range) put \[Y=|x-\Dbar|+|y-\Dbar|\le2\Dbar.\]
Then
\LineE{eq:psi}{|\psi(x,y)|\le\frac{12Y^3}{\Dbar^3},\qquad|\psi(x\pm1,y)-\psi(x,y)|\le\frac{32(Y+1)^2}{\Dbar^3}.}
Indeed $|x-y|\le Y$ and \Line{|\Dbar^2-xy|=|(x-\Dbar)(y-\Dbar)+\Dbar(x-\Dbar)+\Dbar(y-\Dbar)|\le3\Dbar Y,} so with $xy\ge\Dbar^2/4$ we get $|\frac1{xy}-\frac1{\Dbar^2}|\le\frac{12Y}{\Dbar^3}$ and \Line{|\psi(x,y)|=(x-y)^2\left|\frac1{xy}-\frac1{\Dbar^2}\right|\le\frac{12Y^3}{\Dbar^3}.} For the second bound, \Line{|\partial_x\psi|=\big|2(x-y)(\frac1{xy}-\frac1{\Dbar^2})-\frac{(x-y)^2}{x^2y}\big|\le\frac{24Y^2+8Y^2}{\Dbar^3}} (as $x^2y\ge\Dbar^3/8$), and $Y$ changes by at most $1$ along the segment from $x$ to $x\pm1$.

Applying~\eqref{eq:psi} (all arguments are in range by the above), with $u_v:=(d_v-\Dbar)^2+1$ and $Y_{kl}=|d_k-\Dbar|+|d_l-\Dbar|$, $Y_{kl}'$ the same in $\mathcal{G}'$, we get on $\mathcal E_k$
\[|\text{(a)}|\le\frac C{\Dbar^3}\Big[\sum_{l\in N(k)}Y_{kl}^3+\sum_{l\in N'(k)}(Y_{kl}')^3\Big],\qquad|\text{(b)}|\le\frac C{\Dbar^3}(\tilde\nu+1)\sum_{v\in\mathcal D_k}u_v .\]
For (b), each edge of $\mathcal B_k$ contributes at most $\frac{C}{\Dbar^3}(u_i+u_j)$ (two unit steps, $(Y+1)^2\le C(u_i+u_j)$), both endpoints lie in $\mathcal D_k$ and each vertex is on at most $\tilde\nu+1$ edges. We bound second moments, using $\P(v\in\mathcal D_k)\le C(2r_n)^m$, $n(2r_n)^m\le C\Dbar$, and the conditional moment bounds of Lemma~\ref{lem:mom}(iv) (Section~\ref{ts:lem:mom}) given $X_k,X_k',X_v$. For (b), by Cauchy--Schwarz $(\sum_{v\in\mathcal D_k}u_v)^2\le\nu\sum_{v\in\mathcal D_k}u_v^2$ and
\Line{\E\Big[(\tilde\nu+2)^2\Big(\sum_{v\in\mathcal D_k}u_v\Big)^2\Big]\le\sum_{v=1}^n\E\big[\1\{v\in\mathcal D_k\}(\tilde\nu+2)^3u_v^2\big]\le n C(2r_n)^m C\Dbar^3 C\Dbar^2\le C\Dbar^6 .}
For (a), \Line{Y_{kl}^3\le4(|d_k-\Dbar|^3+|d_l-\Dbar|^3)} and $|N(k)|\le\nu$, so \Line{\sum_{l\in N(k)}Y_{kl}^3\le4\nu|d_k-\Dbar|^3+4\sum_{l\in\mathcal D_k}|d_l-\Dbar|^3,} and by Hölder and Lemma~\ref{lem:mom}(i),(ii) (Section~\ref{ts:lem:mom})
\[\begin{aligned}\E[\nu^2|d_k-\Dbar|^6]&\le(\E\nu^4)^{1/2}(\E|d_k-\Dbar|^{12})^{1/2}\le C\Dbar^5,\\
\E\Big[\Big(\sum_{l\in\mathcal D_k}|d_l-\Dbar|^3\Big)^2\Big]&\le n C(2r_n)^m C\Dbar C\Dbar^3\le C\Dbar^5,\end{aligned}\]
and the same for $\mathcal{G}'$. Hence \Line{\E[(n(R_n-R_n^{(k)}))^2\1_{\mathcal E_k}]\le C\Dbar^{-6}(\Dbar^5+\Dbar^6)\le C.}

\medskip\noindent\textbf{The bad event.} On $\mathcal E_k^c$ use~\eqref{eq:crude} and Cauchy--Schwarz:
\Line{\E\big[(n(R_n-R_n^{(k)}))^2\1_{\mathcal E_k^c}\big]\le C\big(\E[(\tilde\nu+2)^{12}]\big)^{1/2}\P(\mathcal E_k^c)^{1/2}\le C\Dbar^6\big(\Dbar e^{-\Dbar/60}\big)^{1/2}\le C,}
since $\E[(\tilde\nu+2)^{12}]\le C\Dbar^{12}$ and $\P(\mathcal E_k^c)\le C\Dbar e^{-\Dbar/60}$. For the last bound, \Line{\mathcal E_k^c\subseteq\{\exists v\in\mathcal D_k\text{ atypical}\}\cup\{|d_k-\Dbar|>\Dbar/4\}\cup\{|d_k'-\Dbar|>\Dbar/4\}.} 
Note that $|\mathcal D_k|=\sum_{l\neq k}\1\{l\in\mathcal D_k\}$.
By the union bound, \Line{\P(\exists v\in\mathcal D_k\text{ atypical})\le\sum_{v\ne k}\E[\1\{v\in\mathcal D_k\}\P(v\text{ atypical}\mid X_k,X_k',X_v)]\le Ce^{-\Dbar/60}\E|\mathcal D_k|\le C\Dbar e^{-\Dbar/60},} because given $(X_k,X_k',X_v)$ the degree of $v$ is a fixed term plus independent Bernoulli variables (Lemma~\ref{lem:mom}(iv) in Section~\ref{ts:lem:mom} and the typical-vertex bound) and $\E|\mathcal D_k|\le\E\nu\le2^{m+1}\Dbar$. Both $d_k$ and $d_k'$ are $\mathrm{Bin}(n-1,(2r_n)^m)$, so by Theorem~\ref{tl:A}  each deviates by more than $\Dbar/4$ with probability at most $2e^{-\Dbar/48}\le Ce^{-\Dbar/60}$. Finally $\Dbar^{13/2}e^{-\Dbar/120}$ is bounded for $\Dbar\ge16$.

In all cases $\E[(R_n-R_n^{(k)})^2]\le C/n^2$, and Theorem~\ref{tl:C} (Section~\ref{ts:C}) gives
\LineE{eq:rem3}{\Var\big[\mathrm{Rem}^{(3)}\big]=\Var [R_n]\le\tfrac12 n\cdot\frac C{n^2}=\frac C{2n}.}

\subsection{Proof of Theorem~\ref{thm:var} (variance)}

Let   $r_n\le r_{\max}<\frac12$ and $nr_n^m\ge1$. By~\eqref{eq:rem1},~\eqref{eq:rem2},~\eqref{eq:rem3} and the elementary inequality \Line{\Var[U+V+W]\le3(\Var [U]+\Var [V]+\Var [W]),} 
we have
\[\Var[\mathrm{Rem}]\le \frac{C}{n}.\] By~\eqref{eq:split}, 
\[\mathcal{F}_n-\E[\mathcal{F}_n]=M_n+\mathrm{Rem}\] with \Line{\Var [M_n]=\sigma^2(r_n)(1+O(1/n))} from~\eqref{eq:varM}. Hence
\[\begin{aligned}\big|\Var[\mathcal{F}_n]-\sigma^2(r_n)\big|&\le\big|\Var [M_n]-\sigma^2(r_n)\big|+2\sqrt{\Var [M_n] \Var[\mathrm{Rem}]}+\Var[\mathrm{Rem}]\\
&\le C\sigma^2(r_n)\Big(\frac1n+\frac1{\sqrt{n(2r_n)^m}}+\frac1{n(2r_n)^m}\Big),\end{aligned}\]
because $\sigma^2(r_n)\asymp (2r_n)^m$.
Since $n(2r_n)^m\asymp\Dbar\ge1$, we have  
\Line{\Big|\frac{\Var[\mathcal{F}_n]}{\sigma^2(r_n)}-1\Big|\le C(n(2r_n)^m)^{-1/2}\le C(nr_n^m)^{-1/2}.}
For $r_n<\frac16$ we have $\sigma^2(r_n)=\sigma_m^2(2r_n)^m$ (Proposition~\ref{prop:alpha} in Section~\ref{sec:step2}), which is Theorem~\ref{thm:var} (Section~\ref{ts:thm:var}); for constant $r$ the same bound with $\sigma^2(r)$ is used in the proof of Theorem~\ref{thm:fixed}(ii) (Section~\ref{ts:thm:fixed}). \hfill$\square$

\subsection{Step 6: the pair statistic $W_n$ is asymptotically Gaussian at the rate $\sqrt{(2r_n)^m}$}\label{sec:step6}
This step proves that the pair statistic $W_n$ is asymptotically normal, with a Berry--Esseen rate, so that the result can be transferred to $\mathcal F_n$ in the proof of Theorem~\ref{thm:be} (Section~\ref{sec:proofbe}). Throughout this step and that proof we assume $r_n<\frac16$; for a constant radius the pair statistic is not asymptotically Gaussian (Section~\ref{sec:fixedproofs}). It rests on four ingredients, presented in the order in which they are used:
\begin{itemize}[leftmargin=*]
\item Theorem~\ref{tl:K} (the Shao--Zhang bound for exchangeable pairs) turns such bounds into a Kolmogorov distance to the normal law;
\item Proposition~\ref{prop:kolm} applies Theorem~\ref{tl:K} to $W_n$ and gives the Berry--Esseen bound for the pair statistic;
\item Lemma~\ref{lem:pert} (the perturbation lemma) is the device that carries this bound over to $\mathcal F_n=M_n+\mathrm{Rem}$ in Section~\ref{sec:proofbe}: it shows that a rescaling by a constant near $1$ and an additive remainder with small second moment change the Kolmogorov distance only slightly.
\end{itemize}

\medskip\noindent\textbf{Normal approximation for exchangeable pairs.} For an exchangeable pair, the Kolmogorov distance to the normal law is bounded by the error in the conditional variance of the increment, the error $R$ in the linear regression, and a third-moment term (Theorem \ref{tl:K} below). We use it for the pair statistic.

\begin{tool}[Kolmogorov bound for exchangeable pairs; Shao--Zhang 2019 \cite{SZ19}, Cor.~2.2 of Thm.~2.1]\label{tl:K}\tlsec{K}
Let $(W,W')$ be an exchangeable pair, $\Delta=W-W'$, and suppose $\E[\Delta\mid W]=\lambda(W+R)$ for a constant $\lambda\in(0,1)$ and a random variable $R$. If $\E[|W|]\le2$, then
\Line{\sup_{z\in\R}\big|\P(W\le z)-\Phi(z)\big|\le\E\big[\Big|1-\frac1{2\lambda}\E\big[\Delta^2\mid W\big]\Big|\big]+\E[|R|]+2\sqrt{\frac{\E [|\E[\Delta^3\mid W]|]}{\lambda}} .}
\end{tool}
\noindent (The source states a more general bound with a free symmetric $\Delta^*\ge|\Delta|$ and the term $\lambda^{-1}\E[|\E[\Delta\Delta^*\mid W]|]$ in its Thm.~2.1; the displayed form is its Cor.~2.2, obtained with $\Delta^*=s/2+\Delta^2/(2s)$.)

\medskip\noindent\textbf{Berry--Esseen bound for the pair statistic.} 
\begin{proposition}[Kolmogorov bound for the pair statistic]\label{prop:kolm}\tlsec{prop:kolm}
If $r_n<\frac16$ and $nr_n^m\ge1$, then
\Line{\dK\big(W_n/\varsigma_n,\mathcal{N}\big)\le C ((2r_n)^{m/2}+n^{-1/2})\le C r_n^{m/2},}
where $\varsigma_n^2=\Var [W_n]$.
\end{proposition}
\noindent The proof of Proposition~\ref{prop:kolm} is based on the Berry--Esseen bound for exchangeable pairs of Shao and Zhang 2019 \cite{SZ19} (Theorem \ref{tl:K}), and is given in Section~\ref{pf:prop:kolm} of the Supplement.

\medskip\noindent\textbf{Perturbation of the Kolmogorov distance.} The next lemma is elementary; it is the  tool used to pass from $W_n$ to $\mathcal F_n$.
\begin{lemma}[Perturbation]\label{lem:pert}\tlsec{lem:pert}
Let $X,Y$ be random variables and $\varepsilon>0$. Then
\Line{\dK(X+Y,\mathcal{N})\le\dK(X,\mathcal{N})+\frac{\varepsilon}{\sqrt{2\pi}}+\P(|Y|>\varepsilon).}
In particular \Line{\dK(X+Y,\mathcal{N})\le\dK(X,\mathcal{N})+C(\E [Y^2])^{1/3}.} Moreover for $c\in[\frac12,2]$: \Line{\dK(cX,\mathcal{N})\le\dK(X,\mathcal{N})+C|c-1|,} and for every $x_0\in\R$: \Line{\dK(X+x_0,\mathcal{N})\le\dK(X,\mathcal{N})+|x_0|/\sqrt{2\pi}.}
\end{lemma}
\noindent The proof of Lemma~\ref{lem:pert} is given in Section~\ref{pf:lem:pert} of the Supplement.

\subsection{Proof of Theorem~\ref{thm:be} (Berry--Esseen bound)}\label{sec:proofbe}

Throughout, $r_n<\frac16$ and $nr_n^m\ge1$, as in the theorem.

\medskip\noindent\emph{Step 1: the standardised statistic.} By~\eqref{eq:split}, \[\mathcal{F}_n-\E[\mathcal{F}_n]=M_n+\mathrm{Rem}\]
and $M_n=\frac{(n-2)(n-3)}{n\Dbar^2}W_n$, so
\Line{\frac{\mathcal{F}_n-\E[\mathcal{F}_n]}{\sqrt{\Var[\mathcal{F}_n]}}=c_n \frac{W_n}{\varsigma_n}+\frac{\mathrm{Rem}}{\sqrt{\Var[\mathcal{F}_n]}},}
where
\[c_n=\frac{(n-2)(n-3)\varsigma_n}{n\Dbar^2\sqrt{\Var[\mathcal{F}_n]}}=\frac{\sqrt{\Var [M_n]}}{\sqrt{\Var[\mathcal{F}_n]}}.\]

\medskip\noindent\emph{Step 2: the constant $c_n$ and the remainder.} By~\eqref{eq:varM} and Theorem~\ref{thm:var} (Section~\ref{ts:thm:var}),
\[\Var [M_n]=\sigma_m^2(2r_n)^m(1+O(1/n)),\qquad\Var[\mathcal{F}_n]=\sigma_m^2(2r_n)^m(1+O(\Dbar^{-1/2})),\]
therefore, 
\[|c_n-1|\le C\Dbar^{-1/2}.\]
Moreover, by~\eqref{eq:rem1}--\eqref{eq:rem3} $\Var[\mathrm{Rem}]\le C/n$, and $\Var[\mathcal{F}_n]\ge c (2r_n)^m$, so
\Line{\E\Big[\Big(\frac{\mathrm{Rem}}{\sqrt{\Var[\mathcal{F}_n]}}\Big)^2\Big]\le\frac{C/n}{c (2r_n)^m}=\frac C{n(2r_n)^m}.}

\medskip\noindent\emph{Step 3: perturbation.} Lemma~\ref{lem:pert}  gives, for $c_n\in[\frac12,2]$,
\[\dK(c_nX+Y,\mathcal{N})\le\dK(X,\mathcal{N})+C|c_n-1|+C(\E [Y^2])^{1/3},\]
where $X=W_n/\varsigma_n$ and $Y=\mathrm{Rem}/\sqrt{\Var[\mathcal{F}_n]}$. Apply Proposition~\ref{prop:kolm} (Section~\ref{ts:prop:kolm}) for $\dK(X,\mathcal{N})$, Step~2 for $|c_n-1|$ and $\E[Y^2]$:
\Line{\dK\Big(\frac{\mathcal{F}_n-\E[\mathcal{F}_n]}{\sqrt{\Var[\mathcal{F}_n]}},\mathcal{N}\Big)\le C\Big((2r_n)^{m/2}+\Dbar^{-1/2}+(n(2r_n)^m)^{-1/3}\Big)\le C\Big(r_n^{m/2}+(nr_n^m)^{-1/3}\Big).}

\medskip\noindent\emph{Step 4: the range of validity.} Step~3 requires $c_n\in[\frac12,2]$. If $\Dbar$ is below a large constant $D_0$, then $(n(2r_n)^m)^{-1/3}\ge c(D_0)>0$ and the claim is trivial (take $C\ge1/c(D_0)$, as $\dK\le1$). If $\Dbar\ge D_0$, then $|c_n-1|\le C\Dbar^{-1/2}\le\frac12$ and Step~3 applies. \hfill$\blacksquare$

\begin{remark}[Why Theorem~\ref{tl:H} does not give Theorem~\ref{thm:be}]\label{rem:noG}
Theorem~\ref{tl:H} (Section~\ref{ts:H}) is used in regime~3 but not in regime~2, where $n(2r_n)^m\to\infty$: its constant depends on the stabilisation and moment constants, which are not uniform in $n(2r_n)^m$ (the score of a point depends on a ball of radius $2r_n$ and $\E[|\xi_n|^{4+p}]\asymp(n(2r_n)^m)^{4+p}$), and a uniform constant would give a bound smaller than the skewness $r_n^{m/2}$ of the leading pair statistic. Theorem~\ref{thm:be} (Section~\ref{ts:thm:be}) is therefore proved directly by the exchangeable-pair argument of Step~6.
\end{remark}

\section{Proofs for constant radii}\label{sec:fixedproofs}
\noindent Throughout, $C$ depends on $(m,r)$; $r_n=r\in(0,1)$ is a constant, $\Dbar=(n-1)(2r_n)^m\asymp n$.

\subsection{Fourier series on the torus (a reminder)}\label{sec:fourier}

 \noindent  Throughout, $f$ and $h$ denote functions on $\T^m$ ($h$ is a generic integrable function, the ``kernel''; in the applications below it is $\beta$ or the normalised kernel $g/(2r_n)^{2m}$), and all integrals are over $\T^m$ with respect to Lebesgue measure. The Fourier transform preserves the $L^2$ norm. We use this on the torus to express $\sigma^2$ and $\sum_\xi\lambda_\xi$ through the Fourier coefficients of the cube indicator.

\begin{enumerate}[label=(\roman*),leftmargin=*,itemsep=3pt]
\item \emph{Fourier series and Parseval.} Every $f\in L^2(\T^m)$ has the Fourier series
\[
f(x)=\sum_{\xi\in\mathbb Z^m}\hat f(\xi) e^{2\pi i\langle\xi,x\rangle},\qquad \hat f(\xi)=\int_{\T^m}f(x) e^{-2\pi i\langle\xi,x\rangle} dx .
\]
The functions $e_\xi(x)=e^{2\pi i\langle\xi,x\rangle}$ are orthonormal and complete, so Parseval's identity \[\int|f|^2=\sum_\xi|\hat f(\xi)|^2\]
holds. If moreover $f$ is continuous and $\sum_\xi|\hat f(\xi)|<\infty$, the series converges uniformly to $f$ (its sum is continuous and has the same Fourier coefficients as $f$, and a continuous function is determined by them).
\item \emph{Product rule.} Let $F(x)=\int_{\T^m}f(u) h(x-u) du$ (the convolution of $f$ and $h$ on $\T^m$). Then $\hat F(\xi)=\hat f(\xi) \hat h(\xi)$. 
\item Suppose $h$ is real and even. By Parseval's identity, 
\[h(x-y)=\sum_{\xi\in\mathbb Z^m}\hat h(\xi)e^{2\pi i\langle\xi,x-y\rangle}\] in $L^2(\T^m)$. For $\xi\ne0$, $h$ real and even gives $\hat h(\xi)=\hat h(-\xi)\in\mathbb R$. Hence the terms for $\xi$ and $-\xi$ combine:
\[
\hat h(\xi)\Big[e^{2\pi i\langle\xi,x-y\rangle}+e^{-2\pi i\langle\xi,x-y\rangle}\Big]=2\hat h(\xi)\cos\big(2\pi\langle\xi,x-y\rangle\big).
\]
Summing over one representative of each pair, that is over $\xi\in\mathbb Z^m_+$,
\[
h(x-y)=\hat h(0)+\sum_{\xi\in\mathbb Z^m_+}2\hat h(\xi)\cos\big(2\pi\langle\xi,x-y\rangle\big).
\]

\item \emph{The cube.} For $\beta(x)=\mathbf 1\{\|x\|_\infty<r\}$, $r<\frac12$, the Fourier coefficient is the number $\eta_\xi$ of~\eqref{eq:ahat},
\[
\eta_\xi=\int_{\T^m}\beta(w)e^{-2\pi i\langle\xi,w\rangle}dw=\prod_{b=1}^m\Big(\frac{\sin(2\pi r\xi_b)}{\pi\xi_b} \1\{\xi_b\ne0\}+2r \1\{\xi_b=0\}\Big),
\]
because the integral factorises over the coordinates. Clearly,  $\eta_\xi$ is real and $\eta_\xi=\eta_{-\xi}$; in particular, $\eta_0=(2r_n)^m$.
\item \emph{The functions $P_2,P_3$.} Recall that \[P_2(x)=\E[\beta(x-U_1)\beta(U_1)],\ \ \ \ \ \ \  \ P_3(x)=\E[\beta(x-U_1)\beta(U_1-U_2)\beta(U_2)]\]
with $U_1,U_2$ independent and uniform. By the product rule (ii), applied once, resp.\ twice, their Fourier coefficients are $\eta_\xi^2$ and $\eta_\xi^3$.
\end{enumerate}

\subsection{The remainder for a constant radius}\label{sec:transfer}
Steps~1--5 and the variance computation of Section~\ref{sec:proofs} were proved there for every radius $r_n\le r_{\max}<\frac12$ with $nr_n^m\ge1$, with $\kappa(r_n)$ and $\sigma^2(r_n)$ in place of $\kappa_m$ and $\sigma_m^2(2r_n)^m$; in particular they apply to $r_n=r$ constant (then $\Dbar=(n-1)(2r_n)^m\asymp n$). Consequently, exactly as in the dense regime,
\LineE{eq:remfixed}{\Var [\mathrm{Rem}]\le \frac{C}{n},\qquad\mathrm{Rem}:=\mathcal{F}_n-\E[\mathcal{F}_n]-M_n,\qquad M_n=\tfrac{(n-2)(n-3)}{n\Dbar^2}W_n .}
(The decomposition $\mathcal{F}_n-\E[\mathcal{F}_n]=M_n+\mathrm{Rem}$ is~\eqref{eq:split}.) What is not available for a constant radius is Step~6: the pair statistic $W_n$ is not asymptotically Gaussian.

The variance bound~\eqref{eq:remfixed} is enough for the weak limit, but the rate in Theorem~\ref{thm:fixed}(iii) (Section~\ref{ts:thm:fixed}) needs the remainder to be small with overwhelming probability, and this requires bounds on its high moments $\|\mathrm{Rem}\|_{L^q}$ with $q\asymp\log n$. The aim of this section is Lemma~\ref{lem:remmom}, which gives such a bound; it is obtained from two moment inequalities, applied to the three parts of \[\mathrm{Rem}=\mathrm{Rem}^{(1)}+\mathrm{Rem}^{(2)}+\mathrm{Rem}^{(3)}\] 
in~\eqref{eq:split}:

\begin{itemize}[leftmargin=*]
\item Lemma~\ref{lem:degmom} bounds the $L^q$ norm of a sum of completely degenerate terms. It handles the three- and four-point Hoeffding components of $\Upsilon_n$ ($\mathrm{Rem}^{(1)}$) and the two- and three-point components of the degree-and-triangle part $S$ ($\mathrm{Rem}^{(2)}$);
\item Theorem~\ref{tl:L} (the moment inequality of Boucheron--Bousquet--Lugosi--Massart) bounds the $L^q$ norm of a fluctuation by the $L^q$ norms of one-point re-sampling effects. It handles the centred denominator error $\mathrm{Rem}^{(3)}=R_n-\E R_n$, for which no degenerate structure is available;
\item Lemma~\ref{lem:remmom} combines the two and gives $\|\mathrm{Rem}\|_{L^q}\le Cq^{3/2}n^{-1/2}$ for $q\le n^{1/3}$.
\end{itemize}

\medskip\noindent\textbf{Degenerate sums.} \begin{lemma}[Moments of degenerate sums]\label{lem:degmom}\tlsec{lem:degmom}
Let $X_1,\dots,X_n$ be independent, $t\ge1$ and $s\ge1$ integers with $n\ge 2st$, and for each $t$-element set $T\subseteq\{1,\dots,n\}$ let $\chi_T(X_T)$ be a random variable with $|\chi_T|\le B$ and $\E[\chi_T\mid X_{T\setminus\{i\}}]=0$ for every $i\in T$ (complete degeneracy). Then, with $q=2s$,
\Line{\Big\|\sum_{|T|=t}\chi_T(X_T)\Big\|_{L^q}\le C_t q^{t/2} n^{t/2}B .}
\end{lemma}
\noindent The proof of Lemma~\ref{lem:degmom} is given in Section~\ref{pf:lem:degmom} of the Supplement.

\medskip\noindent\textbf{One-point re-sampling.} The $L^q$ norm of a fluctuation is bounded by the $L^q$ norms of one-point re-sampling effects, with growth $\sqrt q$. We use it with $q\asymp\log n$ for the remainder (Lemma~\ref{lem:remmom} in Section~\ref{ts:lem:remmom} and Step~5 of Theorem~\ref{thm:fixed}(iii) in Section~\ref{ts:thm:fixed}).
\begin{tool}[moment inequality; Boucheron--Bousquet--Lugosi--Massart 2005 \cite{BBLM05}, Thm.~2]\label{tl:L}\tlsec{L}
Let $X_1,\dots,X_n$ be independent, $Z=f(X_1,\dots,X_n)$ square-integrable, $X_i'$ an independent copy of $X_i$ and $Z_i'=f(X_1,\dots,X_i',\dots,X_n)$. Put \Line{V^+=\E\big[\sum_i((Z-Z_i')_+)^2 \big| X_1,\dots,X_n\big].} Then for every $q\ge2$, \Line{\|(Z-\E Z)_+\|_{L^q}\le\sqrt{\varkappa q \|V^+\|_{L^{q/2}}}} with an absolute constant $\varkappa<1.3$.
\end{tool}
\noindent Applying this to $Z$ and to $-Z$, and using \Line{\|V^+\|_{L^{q/2}}\le\sum_i\|\E[(Z-Z_i')^2|X]\|_{L^{q/2}}\le\sum_i\|Z-Z_i'\|_{L^q}^2} (Minkowski and the conditional Jensen inequality), we obtain the form used below:
\LineE{eq:BBLM}{\|Z-\E Z\|_{L^q}\le2\sqrt{\varkappa q} \Big(\sum_{i=1}^n\|Z-Z_i'\|_{L^q}^2\Big)^{1/2}\qquad(q\ge2).}

\medskip\noindent\textbf{Moments of the remainder.} The two inequalities in Lemma~\ref{lem:degmom} and Theorem~\ref{tl:L} combine into the following moment bound for the remainder.  The detailed proof is in  Section \ref{pf:lem:remmom}.

\begin{lemma}\label{lem:remmom}\tlsec{lem:remmom}
For constant $r\in(0,\frac12)$ there are $C,n_0$ such that for $n\ge n_0$ and every even integer $q$ with $2\le q\le n^{1/3}$, \[\|\mathrm{Rem}\|_{L^q}\le C q^{3/2} n^{-1/2}.\]
\end{lemma}

\subsection{The eigenvalues $\lambda_\xi$ and the limit $Q_\infty$ (Theorem~\ref{thm:fixed}(i))}\label{sec:eig}

Let $r<\frac12$ and let $g(x,y)$ be as in~\eqref{eq:gdef}. We define $h=h(x)$ on $\T^m$ by \[h(x):=\frac{g(x,0)}{(2r)^{2m}}=\frac{2((2r)^{m}P_2(x,0)-P_3(x,0))}{(2r)^{2m}}.\]

\paragraph{(A) Many non-zero eigenvalues.} Take $\xi=(k,0,\dots,0)$. Then 
\[\eta_\xi=(2r)^{m-1}\frac{\sin(2\pi rk)}{(\pi k)}\] 
vanishes iff $2rk\in\mathbb Z$. If $2r$ is irrational this never happens; if $2r=p/d$ in lowest terms ($d\ge2$), then $\eta_\xi\ne0$ for all $k\not\equiv0\pmod d$. So $\eta_\xi\ne0$ for infinitely many $\xi\in\mathbb Z^m_+$.

\paragraph{Fourier coefficients of $h$.} The Fourier coefficients of $P_2(x):=P_2(x,0)$ and $P_3(x):=P_3(x,0)$ are $\eta_\xi^2$ and $\eta_\xi^3$ (Section~\ref{sec:fourier} (v) ), then
\[\hat h(\xi)=\frac{2\eta_\xi^2\big((2r)^m-\eta_\xi\big)}{(2r)^{2m}}=\lambda_\xi\quad(\xi\ne0),\qquad\hat h(0)=0 .\]
For $\xi\ne0$, $|\eta_\xi|<(2r)^m$. Hence $0<(2r)^m-\eta_\xi\le2(2r)^m$, so  $0\leq\lambda_\xi\le4\eta_\xi^2/(2r)^m$, and $\lambda_\xi>0$ iff $\eta_\xi\ne0$; by (A), $\lambda_\xi>0$ for infinitely many $\xi$.

\paragraph{Summability and the value of $\kappa(r)$.} By Parseval's identity,
\[\sum_\xi\eta_\xi^2=\int\beta^2=(2r)^m,\] so $\sum_{\xi\ne0}\lambda_\xi\le4$. Thus $h$ is continuous with absolutely summable Fourier series, which converges to $h$ everywhere. Note that $P_2(0)=(2r)^m$,
\[\sum_{\xi\ne0}\lambda_\xi=h(0)=\frac{2\big((2r_n)^{2m}-\E[A_{12}A_{13}A_{23}]\big)}{(2r)^{2m}}=2\kappa(r),\]
and $\lambda_\xi=\lambda_{-\xi}$. Then
\begin{equation}\label{hfourier}
\sum_{\xi\in\mathbb Z^m_+}\lambda_\xi=\kappa(r),\ \ \ \ \  \ \ h(x)=\sum_{\xi\in\mathbb Z^m_+}2\lambda_\xi\cos(2\pi\langle\xi,x\rangle).
\end{equation}

\paragraph{The series $Q_\infty$.} The variables $\lambda_\xi(|Z_\xi|^2-1)$ are independent and centred with variance $\lambda_\xi^2$, and \[0<\sigma^2=\sum_{\xi\in\mathbb Z^m_+}\lambda_\xi^2\le\left(\sum_{\xi\in\mathbb Z^m_+}\lambda_\xi\right)^2<\infty.\]
By Kolmogorov's theorem the series converges a.s.\ and in $L^2$, $\E[Q_\infty]=0$, $\Var[Q_\infty]=\sigma^2$. Since $\sum_{\xi\in\mathbb Z^m_+}\lambda_\xi|Z_\xi|^2$ has non-negative terms and mean $\kappa(r)<\infty$, it converges a.s., and \Line{Q_\infty=\sum_{\xi\in\mathbb Z^m_+}\lambda_\xi|Z_\xi|^2-\kappa(r)\ge-\kappa(r).} 

\paragraph{Density.} Let \Line{Y:=Q_\infty+\kappa(r)=\sum_{\xi\in\mathbb Z^m_+}\lambda_\xi|Z_\xi|^2,} the a.s.\ limit of $Y_N=\sum_{\xi\in F_N}\lambda_\xi|Z_\xi|^2$ over finite sets $F_N\uparrow\mathbb Z^m_+$. For $E\sim\mathrm{Exp}(1)$, $\E[e^{iu\lambda E}]=(1-i\lambda u)^{-1}$, so by independence and dominated convergence \Line{\varphi_Y(u)=\lim_N\prod_{\xi\in F_N}(1-i\lambda_\xi u)^{-1}.} Each factor has modulus $(1+\lambda^2u^2)^{-1/2}\le1$; choosing two positive eigenvalues $\lambda_1,\lambda_2$ (by (A)),
\[|\varphi_Y(u)|\le(1+\lambda_1^2u^2)^{-1/2}(1+\lambda_2^2u^2)^{-1/2}\le\frac1{\lambda_1\lambda_2u^2},\]
so $\varphi_Y\in L^1(\R)$. By Fourier inversion $Y$ has a bounded continuous density \[f_Y(y)=\frac1{2\pi}\int e^{-iuy}\varphi_Y(u)du,\]
hence $Q_\infty$ has the bounded continuous density $f_\infty(x)=f_Y(x+\kappa(r))$; this is used in the proof of Theorem~\ref{thm:fixed}(iii) (Section~\ref{ts:thm:fixed}).

\paragraph{The value of $\sigma^2$.} Since  $\hat h(\xi)=\hat h(-\xi)$, then $\sum_{\xi\ne0}\hat h(\xi)^2=2\sigma^2$. By Parseval's identity and $\hat h(0)=0$, \Line{\sum_{\xi\ne0}\hat h(\xi)^2=\int_{\T^m} h(x)^2dx=\E[h(X_k-X_l)^2]} ($X_k-X_l$ is uniform). Hence, with $P_2=P_2(X_k,X_l)$, $P_3=P_3(X_k,X_l)$ and~\eqref{eq:gdef},
\[\sigma^2=\frac12 \E[h(X_k-X_l)^2]=\frac{\E[g(X_k,X_l)^2]}{2(2r)^{4m}}=\frac{2 \E\big[\big((2r)^mP_2-P_3\big)^2\big]}{(2r)^{4m}}.\]
By Proposition~\ref{prop:alpha}(a) (Section~\ref{ts:prop:alpha}), $\E[P_2^2]=\alpha_1^m$, $\E[P_2P_3]=\alpha_2^m$, $\E[P_3^2]=\alpha_3^m$ with $\alpha_i$ the piecewise polynomials~\eqref{eq:alpha1}--\eqref{eq:alpha3}, hence
\Line{\sigma^2=\frac{2\big[(2r)^{2m}\alpha_1^m-2(2r)^m\alpha_2^m+\alpha_3^m\big]}{(2r)^{4m}}.}

\subsection{Proof of Theorem~\ref{thm:fixed}(ii)--(iii)}\label{sec:fixedproof}
Part~(ii) is a direct computation from the variance of $W_n$ and the remainder bound~\eqref{eq:remfixed}. For part~(iii), the remainder is discarded using Lemma~\ref{lem:remmom} (Section~\ref{ts:lem:remmom}), and what is left is $M_n=c_nT_n$ with $T_n$ a degenerate $U$-statistic with a bounded kernel on the torus. Two theorems on such statistics are needed, and they play different roles: Theorem~\ref{tl:F} (Section~\ref{ts:F}) identifies the weak limit of $T_n$ as the weighted sum $Q_\infty$ of centred $\chi^2$ variables whose weights are the eigenvalues of the kernel, so it gives the limit law; Theorem~\ref{tl:G} (Section~\ref{ts:G}) gives the rate $n^{-1}\log n$ at which $T_n$ approaches $Q_\infty$, which is what turns the limit law into the rate stated in (iii). The proof is organised accordingly: part~(ii), the two tools, and then Steps~1--5 of part~(iii).

\medskip\noindent\textbf{Part (ii): variance.} Recall that $\Dbar=(n-1)(2r_n)^m$, 
\[M_n=\frac{(n-2)(n-3)}{n\Dbar^2}W_n,\ \ \ \ \Var [W_n]=\binom n2\E[g(X_k,X_l)^2],  \ \ \ \ \ \ \sigma^2=\frac{\E[g(X_k,X_l)^2]}{2(2r)^{4m}}.\]
Then
{\relpenalty=0 \binoppenalty=700 \Line{\Var [M_n]=\frac{(n-2)^2(n-3)^2}{n(n-1)^3} \sigma^2=\sigma^2\Big(1+O\big(\tfrac1n\big)\Big).}
By~\eqref{eq:remfixed},
\[\begin{aligned}|\Var[\mathcal{F}_n]-\sigma^2|&\le|\Var [M_n]-\sigma^2|+2\sqrt{\Var [M_n]\Var[\mathrm{Rem}]}+\Var[\mathrm{Rem}]\\&\le C/n+C/\sqrt n\le Cn^{-1/2}.\end{aligned}\]
}
\medskip\noindent\textbf{Tools for part (iii).} A degenerate $U$-statistic converges to a weighted sum of centred $\chi^2_1$ variables, with the kernel's eigenvalues as weights. This identifies $Q_\infty$ (Step~2 below).
\begin{tool}[limit of degenerate $U$-statistics; Gregory 1977 \cite{Gregory77}, see Serfling \cite{Serfling80}, \S5.5.2]\label{tl:F}\tlsec{F}
Let $X_1,X_2,\dots$ be i.i.d.\ and $h$ a symmetric kernel with $\E [h(X_1,X_2)^2]<\infty$ and $\E[h(x,X_1)]=0$ for a.e.\ $x$. Let $h(x,y)=\sum_{j\ge1}\lambda_j\varphi_j(x)\varphi_j(y)$ be its spectral expansion in $L^2$ (orthonormal $\varphi_j$; $\sum\lambda_j^2<\infty$) and $T_n=\frac1n\sum_{1\le k<l\le n}h(X_k,X_l)$. Then $T_n\Rightarrow\frac12\sum_j\lambda_j(\zeta_j^2-1)$ with $\zeta_j$ i.i.d.\ $\mathcal{N}(0,1)$.
\end{tool}

\noindent For a bounded degenerate kernel with enough non-zero eigenvalues, the convergence in Theorem~\ref{tl:F} (Section~\ref{ts:F}) holds at the rate $n^{-1}\log n$. This gives the rate for a constant radius (Step~3 below).
\begin{tool}[rate for degenerate $U$-statistics; Bentkus--G\"otze 1999 \cite{BG99}, G\"otze--Tikhomirov 1999 \cite{GT99}]\label{tl:G}\tlsec{G}
In the setting of the previous theorem assume in addition that $h$ is bounded and has at least nine non-zero eigenvalues. Then, with $Q=\frac12\sum_j\lambda_j(\zeta_j^2-1)$, \Line{\sup_x|\P(T_n\le x)-\P(Q\le x)|\le C_h n^{-1}\log n} (the original statements are for the normalisation $\frac{2n}{n-1}T_n$; the logarithm absorbs the change of scale).
\end{tool}

\medskip\noindent\textbf{Part (iii): limit law and rate.} We use Theorems~\ref{tl:F} and~\ref{tl:G}.

\medskip\noindent\textbf{Step 1 (spectral expansion).} Recall that 
\[P_2(x,y)=P_2(x-y,0),\ \ \ \ \ \ \ \ P_3(x,y)=P_3(x-y,0),\]
which implies $g(x,y)=g(x-y,0)$.
Then $M_n=c_nT_n$ with $c_n=\frac{(n-2)(n-3)}{(n-1)^2}\to1$ and \[T_n=\frac{\sum_{k<l}g(X_k,X_l)}{n(2r)^{2m}}=\frac{\sum_{k<l}g(X_k-X_l,0)}{n(2r)^{2m}}=\frac1n\sum_{k<l}h(X_k-X_l).\]
By  \eqref{hfourier}, we get
\Line{h(x-y)=\sum_{\xi\in\mathbb Z^m_+}\lambda_\xi\big(2\cos(2\pi\langle\xi,x\rangle)\cos(2\pi\langle\xi,y\rangle)+2\sin(2\pi\langle\xi,x\rangle)\sin(2\pi\langle\xi,y\rangle)\big).}

\medskip\noindent\textbf{Step 2 (weak limit of $T_n$).} Theorem~\ref{tl:F} (Section~\ref{ts:F}) gives
\Line{T_n\Rightarrow\tfrac12\sum_{\xi\in\mathbb Z^m_+}\lambda_\xi\big(\zeta_{\xi,1}^2+\zeta_{\xi,2}^2-2\big)\overset d=\sum_{\xi\in\mathbb Z^m_+}\lambda_\xi\big(|Z_\xi|^2-1\big)=Q_\infty,}
where the $\zeta_{\xi,\iota}$ are i.i.d.\ $\mathcal{N}(0,1)$: they are the limits of $n^{-1/2}\sum_k\varphi_{\xi,\iota}(X_k)$ for the orthonormal centred eigenfunctions $\varphi_{\xi,1}=\sqrt2\cos(2\pi\langle\xi,\cdot\rangle)$, $\varphi_{\xi,2}=\sqrt2\sin(2\pi\langle\xi,\cdot\rangle)$, so by the multivariate CLT they are jointly Gaussian with identity covariance. Finally, $\frac12(\zeta_{\xi,1}^2+\zeta_{\xi,2}^2)\sim\mathrm{Exp}(1)$, i.e.\ equals $|Z_\xi|^2$ in law, independently over $\xi$.

\medskip\noindent\textbf{Step 3 (rate for $T_n$).} Theorem~\ref{tl:G} (Section~\ref{ts:G}) requires a symmetric degenerate bounded kernel with at least nine non-zero eigenvalues. Clearly, the kernel is bounded, and part~(A) of Section~\ref{sec:eig} shows that the kernel has infinitely many non-zero $\lambda_\xi$. The limit in Theorem~\ref{tl:G} (Section~\ref{ts:G}) is $Q_\infty$ by Step~2, hence \[\dK(T_n,Q_\infty)\le C n^{-1}\log n.\]

\medskip\noindent\textbf{Step 4 (from $T_n$ to $M_n$).}
(a) For $n\ge8$, $c_n\in[\frac12,1)$ and \[|c_n^{-1}-1|=\frac{1-c_n}{c_n}\le\frac6n.\]
(b) $Q_\infty$ has a bounded density, $Q_\infty\ge-\kappa(r)$, and for some small positive constant $c$,
\[\P(Q_\infty>x)\le Ce^{-cx}.\]
Indeed, for $0<s<\frac1{2\max\lambda_\xi}$, $\E[e^{s\lambda E}]=(1-s\lambda)^{-1}$ for $E\sim\mathrm{Exp}(1)$, so by monotone convergence \Line{\E[e^{sQ_\infty}]=\prod_\xi e^{-s\lambda_\xi}(1-s\lambda_\xi)^{-1}\le e^{s^2\sigma^2}} (using $-\log(1-y)\le y+y^2$ for $0\le y\le\frac12$), and Markov's inequality gives the tail bound with $c=s$.\\

(c) Since $c_n>0$, \Line{|\P(c_nT_n\le x)-\P(Q_\infty\le x)|\le\dK(T_n,Q_\infty)+\Delta(x)} with \Line{\Delta(x)=|\P(Q_\infty\le x/c_n)-\P(Q_\infty\le x)|.} For $|x|\le C'\log n$, \Line{\Delta(x)\le\|f_\infty\|_\infty|x||c_n^{-1}-1|\le Cn^{-1}\log n.} For $x>C'\log n$ (with $C'\ge2/c$), both $\P(Q_\infty>x)$ and $\P(Q_\infty>x/c_n)$ are at most $Ce^{-cx}\le n^{-1}$, so $\Delta(x)\le n^{-1}$. For $x<-2\kappa(r)$, both probabilities vanish since $x/c_n<x<-\kappa(r)$. These ranges cover $\R$ for large $n$. Hence \[\dK(M_n,Q_\infty)\le Cn^{-1}\log n.\]

\medskip\noindent\textbf{Step 5 (from $M_n$ to $\mathcal F_n$).} By Lemma~\ref{lem:remmom} (Section~\ref{ts:lem:remmom}), \[\|\mathrm{Rem}\|_{L^q}\le Cq^{3/2}n^{-1/2}\] for even $q\le n^{1/3}$. Let $q=2\lceil\frac12\log n\rceil$ and $\varepsilon=e Cq^{3/2}n^{-1/2}$; Markov's inequality in $L^q$ gives \[\P(|\mathrm{Rem}|>\varepsilon)\le e^{-q}\le n^{-1}.\]
If $|\mathrm{Rem}|\le\varepsilon$, then \Line{\{M_n\le x-\varepsilon\}\subseteq\{\mathcal F_n-\E\mathcal F_n\le x\}\subseteq\{M_n\le x+\varepsilon\}.} Therefore,
\begin{align*}
  \P(\mathcal F_n-\E\mathcal F_n\le x)&= \P(\mathcal F_n-\E\mathcal F_n\le x,|\mathrm{Rem}|>\varepsilon) +\P(\mathcal F_n-\E\mathcal F_n\le x,|\mathrm{Rem}|\leq\varepsilon) \\
  &\leq n^{-1}+\P(M_n\le x+\varepsilon,|\mathrm{Rem}|\leq\varepsilon)\\
  &\leq n^{-1}+\P(M_n\le x+\varepsilon)
  \\
  &\leq n^{-1}+\P(Q_\infty\le x+\varepsilon)+\dK(M_n,Q_\infty),
\end{align*}
and
\begin{align*}
  \P(\mathcal F_n-\E\mathcal F_n\le x)&= \P(\mathcal F_n-\E\mathcal F_n\le x,|\mathrm{Rem}|>\varepsilon) +\P(\mathcal F_n-\E\mathcal F_n\le x,|\mathrm{Rem}|\leq\varepsilon) \\
  &\geq \P(M_n\le x-\varepsilon,|\mathrm{Rem}|\leq\varepsilon)\\
  &\geq \P(M_n\le x-\varepsilon)-\P(|\mathrm{Rem}|>\varepsilon)
  \\
  &\geq \P(Q_\infty\le x-\varepsilon)-\dK(M_n,Q_\infty)-n^{-1}.
\end{align*}
Since \Line{\P(Q_\infty\le x+\varepsilon)-\P(Q_\infty\le x)\le\|f_\infty\|_\infty\varepsilon} (and similarly for $x-\varepsilon$),
\[\begin{aligned}\dK(\mathcal{F}_n-\E[\mathcal{F}_n],Q_\infty)&\le\|f_\infty\|_\infty\varepsilon+\dK(M_n,Q_\infty)+\tfrac1n\\&\le C\big((\log n)^{3/2}n^{-1/2}+n^{-1}\log n\big)\le Cn^{-1/2}(\log n)^{3/2}.\end{aligned}\]
\hfill$\blacksquare$

\section{Proofs for the sparse regime}\label{sec:sparseproofs}

Throughout this section $c_{\min}\le nr_n^m\le c_{\max}$. Then \Line{(2r_n)^m=2^mr_n^m\in[2^mc_{\min}/n, 2^mc_{\max}/n],} $n(2r_n)^m\asymp1$, $\Dbar\asymp1$ and $\Dbar\ge2^{m-1}c_{\min}=:c_*$ for $n\ge2$, so Lemma~\ref{lem:mom} (Section~\ref{ts:lem:mom}) applies; and $r_n\le(c_{\max}/n)^{1/m}\le\frac1{64}$ for $n\ge n_0$. Constants may depend on $m,c_{\max},c_{\min}$.

The proof of Theorem~\ref{thm:sparse} has two parts. Section~\ref{sec:sp1} shows that $\Var[\mathcal F_n]\asymp1/n$. The upper bound follows from the Efron--Stein inequality (Theorem~\ref{tl:C} in Section~\ref{ts:C}), because moving one point changes $n\mathcal F_n$ by a quantity with bounded moments; the lower bound follows from a locality construction that finds order $n$ disjoint regions on which $n\mathcal F_n$ has independent, non-degenerate fluctuations. Section~\ref{sec:sp2} then proves the Gaussian limit with the rate $n^{-1/2}$ by applying Theorem~\ref{tl:H} (Section~\ref{ts:H}), a general normal approximation for sums of local functionals of i.i.d.\ points. This theorem is needed because $n\mathcal F_n$ is a sum of vertex scores that depend only on points within distance $2r_n$ (bounded range) and have bounded moments, while the variance of order $n$ from Section~\ref{sec:sp1} is exactly the lower bound that the theorem requires.

\subsection{Proof of Theorem~\ref{thm:sparse}(i)}\label{sec:sp1}

\paragraph{Upper bound $\Var[\mathcal{F}_n]\le C/n$.} We apply the Efron--Stein inequality (Theorem~\ref{tl:C} in Section~\ref{ts:C}) to $n\mathcal{F}_n$; it suffices to show that moving one point changes $n\mathcal F_n$ by an amount with bounded second moment.

\medskip\noindent\emph{The statistic.} By Proposition~\ref{prop:id} (Section~\ref{ts:prop:id}),
\Line{n\mathcal{F}_n=\sum_{i<j}A_{ij}T_{ij},\qquad T_{ij}=\frac{d_i}{d_j}+\frac{d_j}{d_i}-2\le d_i+d_j\quad(A_{ij}=1),}
because adjacent vertices have degrees $\ge1$. A unit change of one of the two degrees changes $T_{ij}$ by at most $1+\max(d_i,d_j)$, as in~\eqref{eq:crude}.

\medskip\noindent\emph{Moving one point.} Replace $X_k$ by an independent copy $X_k'$. This changes the edges at $k$, and the degrees of the neighbours of $k$ and of $k'$. The second endpoints of the affected edges can have degrees exceeding $\nu_k$, so, exactly as in~\eqref{eq:crude}, the count is in terms of $\tilde\nu=\tilde\nu_k+\tilde\nu_k'$, where
\Line{\tilde\nu_k=\#\{j\ne k:d(X_j,X_k)<3r_n\}}
(Lemma~\ref{lem:mom}(ii) in Section~\ref{ts:lem:mom}). Hence
\Line{|n\mathcal{F}_n-n\mathcal{F}_n^{(k)}|\le C(\tilde\nu+2)^3.}

\medskip\noindent\emph{Conclusion.} Since $nr_n^m\asymp 1$, every moment of $\tilde\nu$ is bounded (Lemma~\ref{lem:mom}(ii) in Section~\ref{ts:lem:mom}), so $\E[(n\mathcal{F}_n-n\mathcal{F}_n^{(k)})^2]\le C$. Theorem~\ref{tl:C} (Section~\ref{ts:C}) gives $\Var[n\mathcal{F}_n]\le\frac12nC$, that is, \[\Var[\mathcal{F}_n]\le \frac{C}{n}.\]

\paragraph{Lower bound $\Var[\mathcal{F}_n]\ge c/n$.} The idea is to find of order $n$ well-separated cubes in which $n\mathcal F_n$ has independent, non-degenerate fluctuations. The proof has three steps.

\medskip\noindent\textbf{Step 1: the cubes.} Choose an even integer $2K$ with
\Line{2K\in\big[\tfrac1{18r_n},\tfrac1{9r_n}\big]}
(the interval has length $\ge2$ because $r_n\le\frac1{64}$). Cut $\T^m$ into $(2K)^m$ congruent cubes of side
\Line{\varpi=\tfrac1{2K}\in[9r_n,18r_n].}
That is, cut each coordinate axis of $\T^m$ into $2K$ intervals $[a\varpi,(a+1)\varpi)$, $a\in\{0,1,\dots,2K-1\}$. The torus is then cut into $(2K)^m$ congruent cubes
\[
\mathcal C_{\mathbf a}=\prod_{b=1}^m[a_b\varpi,(a_b+1)\varpi),\qquad\mathbf a=(a_1,\dots,a_m)\in\{0,\dots,2K-1\}^m,
\]
each of side $\varpi$; $\mathbf a$ is the integer grid coordinate of the cube. Let $\mathcal Q_1,\dots,\mathcal Q_{K^m}$ be those cubes whose integer grid coordinates are all even. That is,
\[
\mathcal Q_1,\dots,\mathcal Q_{K^m}\in\{\mathcal C_{\mathbf a}:\ a_1,\dots,a_m\ \text{all even}\}.
\]
Each coordinate has $K$ even values in $\{0,\dots,2K-1\}$, so there are exactly $K^m$ such cubes, and
\Line{K^m\ge(36r_n)^{-m}\ge c n}
 using $r_n^{-m}\ge n/c_{\max}$. Two different $\mathcal Q_j$ are at distance $\ge\varpi>8r_n$. Let $\mathcal Q_j^{+}$ and $\mathcal Q_j^{++}$ be the sets of points at distance $<2r_n$, respectively $<4r_n$, from $\mathcal Q_j$. The sets $\mathcal Q_j^{++}$ are pairwise disjoint, because $\varpi>8r_n$, and each has volume
\Line{V=\mathrm{vol}(\mathcal Q_j^{++})\le(\varpi+8r_n)^m\le(26r_n)^m.}

\medskip\noindent\textbf{Step 2: conditional independence.} Let $L_i\in\{0,1,\dots,K^m\}$ be the label of the region containing $X_i$: $L_i=j$ if $X_i\in\mathcal Q_j$, and $L_i=0$ if $X_i\notin\bigcup_j\mathcal Q_j$. Let $\mathcal B$ be the $\sigma$-field generated by all labels $L_1,\dots,L_n$ and by the positions of the points with $L_i=0$. Given $\mathcal B$, the points with $L_i=j$ are independent and uniform on $\mathcal Q_j$, because the points are i.i.d.

The summand $\Delta_i$ of $n\mathcal F_n$ depends only on the points within distance $2r_n$ of $X_i$. Put
\Line{\Delta_{(j)}=\sum_{i: X_i\in\mathcal Q_j^+}\Delta_i,\qquad\Delta_{(0)}=\sum_{i: X_i\notin\bigcup_j\mathcal Q_j^+}\Delta_i,}
so that \[n\mathcal{F}_n=\Delta_{(0)}+\sum_{j=1}^{K^m}\Delta_{(j)}.\]
Then $\Delta_{(0)}$ depends only on points outside $\bigcup_j\mathcal Q_j$, hence is $\mathcal B$-measurable, and $\Delta_{(j)}$ is a function of the points in $\mathcal Q_j^{++}$, those outside $\mathcal Q_j$ being $\mathcal B$-measurable. Thus, given $\mathcal B$, the $\Delta_{(j)}$ are conditionally independent functions of the points inside $\mathcal Q_j$, and
\Line{\Var[n\mathcal{F}_n]\ge\E[\Var[n\mathcal{F}_n\mid\mathcal B]]=\sum_{j=1}^{K^m}\E[\Var[\Delta_{(j)}\mid\mathcal B]].}

\medskip\noindent\textbf{Step 3: a non-degenerate event in each cube.} Let $\mathcal E_j$ be the event that exactly three points lie in $\mathcal Q_j$ and no point lies in $\mathcal Q_j^{++}\setminus\mathcal Q_j$:
\[
\mathcal E_j=\{\text{exactly three points lie in }\mathcal Q_j\}\cap\{\text{no point lies in }\mathcal Q_j^{++}\setminus\mathcal Q_j\}.
\]
The first set is determined by the labels. A point in $\mathcal Q_j^{++}\setminus\mathcal Q_j$ does not lie  in any $\mathcal Q_j$, has label $0$, and its position is part of $\mathcal B$. So the event is $\mathcal B$-measurable and $\mathcal E_j\in\mathcal B$.

Fix a set $S\subset\{1,\dots,n\}$ with $|S|=3$ and consider the event $\mathcal E_j^S$ that $X_i\in\mathcal Q_j$ for $i\in S$ and $X_i\notin\mathcal Q_j^{++}$ for $i\notin S$. By independence and uniformity,
\[
\P(\mathcal E_j^S)=\prod_{i\in S}\mathrm{vol}(\mathcal Q_j)\prod_{i\notin S}\big(1-\mathrm{vol}(\mathcal Q_j^{++})\big)=\varpi^{3m}(1-V)^{n-3}.
\]
The events $\mathcal E_j^S$ for different $S$ are disjoint (exactly three points in $\mathcal Q_j$), and $\mathcal E_j=\bigcup_S\mathcal E_j^S$. There are $\binom n3$ sets $S$. Then
\Line{\P(\mathcal E_j)=\binom n3\varpi^{3m}(1-V)^{n-3}\ge c_*'>0,}
since $n\varpi^m\in[(9r_n)^mn,(18r_n)^mn]\subset[9^mc_{\min},18^mc_{\max}]$ and $nV\le26^mc_{\max}$.

On $\mathcal E_j$ the three points are i.i.d.\ uniform on $\mathcal Q_j$ given $\mathcal B$, and $\Delta_{(j)}$ is the sum of $\Delta_i$ over the three vertices and  $\Delta_i$ only depends on the three points. Let $X_1,X_2,X_3$ be the three points. If $1,2,3$ are isolated, then $\Delta_{(j)}=0$. If $1,2$ form an edge, and $3$ is isolated, then $\Delta_{(j)}=0$. If $1,2,3$ form a path, then $\Delta_{(j)}=1$. If $1,2,3$ form a triangle, then $\Delta_{(j)}=0$.
Both outcomes ($\Delta_{(j)}=0$ and $\Delta_{(j)}=1$) have conditional probability greater than some positive constant $c''$:
\begin{itemize}[leftmargin=*]
\item \emph{Two-path.} Suppose $X_1$ is the center of the path. Require $X_1$ to be at distance greater than $r_n$ from the faces of $\mathcal Q_j$ (probability $\ge(1-2r_n/\varpi)^m\ge(7/9)^m$); $X_2$ to have first coordinate in $(X_{1,1}+\frac{r_n}2,X_{1,1}+r_n)$ and the other coordinates within $\frac{r_n}2$ of those of $X_1$; and $X_3$ likewise with first coordinate in $(X_{1,1}-r_n,X_{1,1}-\frac{r_n}2)$ and the other coordinates within $\frac{r_n}2$ of those of $X_1$. Then $X_2,X_3$ are neighbours of $X_1$ but at distance greater than $r_n$ from each other, and each requirement has probability greater than $(r_n/(4\varpi))^m\ge72^{-m}$.
\item \emph{Empty graph.} Cut the first coordinate range of $\mathcal Q_j$ into five intervals of length $\varpi/5>r_n$ and require the three points to fall into the first, third and fifth. They are then pairwise at distance at least $\varpi/5>r_n$, with probability $5^{-3}$.
\end{itemize}
Let $c''=\min(5^{-3},(7/9)^m72^{-2m})\le\frac12$. On $\mathcal E_j$, $\Delta_{(j)}\in\{0,1\}$ with
\[
\P(\Delta_{(j)}=1\mid\mathcal B)\ge c'',\qquad\P(\Delta_{(j)}=0\mid\mathcal B)\ge c'',
\]
 $p:=\P(\Delta_{(j)}=1\mid\mathcal B)\in[c'',1-c'']$ and
\[
\Var[\Delta_{(j)}\mid\mathcal B]=p(1-p)\ge c''(1-c'')\ge\frac{c''}2 .
\]
Since $\mathcal E_j\in\mathcal B$ and $\P(\mathcal E_j)\ge c_*'$,
\[
\E\big[\Var[\Delta_{(j)}\mid\mathcal B]\big]\ge\E\big[\Var[\Delta_{(j)}\mid\mathcal B]\1_{\mathcal E_j}\big]\ge\frac{c''}2 \P(\mathcal E_j)\ge\frac{c''c_*'}2=:c_0 .
\]
Summing over the $K^m\ge cn$ cubes in Step~2 gives 
\[\Var[n\mathcal{F}_n]\ge c K^m\ge cn.\] \hfill$\square$

\subsection{Proof of Theorem~\ref{thm:sparse}(ii)}\label{sec:sp2}

\noindent  Sums of local functionals of i.i.d.\ points with exponential stabilisation, bounded moments and variance of order $n$ are normal at the rate $n^{-1/2}$ (Lachi\`eze-Rey--Schulte--Yukich \cite{LSY19}). We use this for Theorem~\ref{thm:sparse}(ii) (Section~\ref{ts:thm:sparse}), where each vertex score depends on points within distance $2r_n$.
\begin{tool}[normal approximation for local functionals; Lachi\`eze-Rey--Schulte--Yukich 2019 \cite{LSY19}, Thm.~2.1(b) with Cor.~2.2(b)]\label{tl:H}\tlsec{H}
Let $\mathcal X_n=\{X_1,\dots,X_n\}$ be i.i.d.\ with a probability law $\mathbb Q$ on a space with a measurable semi-metric satisfying the growth condition~(2.1) of \cite{LSY19} (which holds for the uniform law on $\T^m$ with the $L_\infty$ torus distance, with exponent $\gamma=m$), and let $H_n=\sum_{i=1}^n\xi_n(X_i;\mathcal X_n)$ where the scores satisfy, for every $n\ge9$, every $x$ and every set $\mathcal A$ of at most $7$ additional points, (a) \emph{exponential stabilisation}: $\xi_n(x;\mathcal X)$, with $\mathcal X=\mathcal X_{n-8}\cup\{x\}\cup \mathcal A$, is unchanged if $\mathcal X$ is cut to the ball $B(x,R_n)$ (condition (2.3) of \cite{LSY19}), and \Line{\P(R_n\ge r)\le C\exp(-c(n^{1/m}r)^\alpha);} (b) the moment bound \Line{\sup_{n\ge9}\sup_x\E[|\xi_n(x;\mathcal X_{n-8}\cup\{x\}\cup\mathcal A)|^{4+p}]<\infty} for some $p\in(0,1]$ (the decay condition~(2.9) of \cite{LSY19} is trivial when $K$ is the whole space). If $\Var [H_n]\ge cn$, then \[\dK\left(\frac{H_n-\E [H_n]}{\sqrt{\Var [H_n]}},\mathcal{N}\right)\le C n^{-1/2}.\]
\end{tool}
\noindent We use Cor.~2.2(b) of \cite{LSY19} with $K$ the whole space. We apply Theorem~\ref{tl:H} to $H_n=n\mathcal{F}_n$ with $\mathbb Q$ the uniform law on $\T^m$, the metric $d$ of Definition~\ref{def:model}, and $K=\T^m$; constants depend only on $m,c_{\max},c_{\min}$.

\medskip\noindent\textbf{Step 1: scores and stabilisation.} For a finite $\mathcal X\subset\T^m$ and $x\in\mathcal X$ let \Line{\begin{gathered}d_x(\mathcal X)=\#\{y\in\mathcal X\setminus\{x\}:d(x,y)<r_n\},\\ \xi_n(x;\mathcal X)=\frac1{d_x(\mathcal X)}\sum_{y: d(x,y)<r_n}d_y(\mathcal X)-d_x(\mathcal X)\end{gathered}} if $d_x(\mathcal X)\ge1$, $\xi_n=0$ otherwise. This is $\Delta_i$ of Definition~\ref{def:stat} for $x=X_i$, so \[H_n=n\mathcal{F}_n=\sum_{i=1}^n\xi_n(X_i;\mathcal X_n).\] Theorem~\ref{tl:H} is applied to configurations $\mathcal X=\mathcal X_{n-8}\cup\{x\}\cup\mathcal A$ with $|\mathcal A|\le7$. Here, $\mathcal X_{n-8}=\{X_1,X_2,\dots,X_{n-8}\}$, and $\mathcal A$ is a fixed (deterministic) set of at most 7 further points.  Every point entering $d_x$ or $d_y$ ($y$ a neighbour of $x$) lies within $2r_n$ of $x$ by the triangle inequality, then \Line{\xi_n(x;\mathcal X)=\xi_n(x;\mathcal X\cap B(x,2r_n)).}  The stabilisation radius is the deterministic $R_n=2r_n\le2(c_{\max}/n)^{1/m}$, and the tail condition (a) of Theorem~\ref{tl:H} holds trivially (e.g.\ with $\alpha=c=1$, $C=\exp(2c_{\max}^{1/m})$), since $\P(R_n\ge r)=0$ for $r>2(c_{\max}/n)^{1/m}$ and the right-hand side is $\ge1$ otherwise.

\medskip\noindent\textbf{Step 2: moments.} Let $\nu=\#(\mathcal X\cap B(x,2r_n))\setminus\{x\}$. As in Lemma~\ref{lem:mom}(iii) (Section~\ref{ts:lem:mom}), $d_x(\mathcal X)\le\nu$ and $d_y(\mathcal X)\le\nu+1$ for neighbours $y$, so $|\xi_n(x;\mathcal X)|\le\nu+1$. Moreover $\nu\le B+7$ with $B\sim\mathrm{Bin}(n-8,(4r_n)^m)$ and $\E B\le4^mc_{\max}$, and the moments of a binomial variable are bounded in terms of its mean (factorial moments are at most $(\E B)^j$). Hence $\E[|\xi_n(x;\mathcal X)|^5]\le C$ uniformly in $n,x,\mathcal A$, which is condition (b) of Theorem~\ref{tl:H} with $p=1$.

\medskip\noindent\textbf{Step 3: variance, growth and decay conditions.} By part (i), $\Var[H_n]=n^2\Var[\mathcal F_n]\ge cn$. The growth condition (2.1) of \cite{LSY19} holds with $\gamma=m$, $\kappa=2^m$, since $\mathbb Q(B(x,r))=\min(2r,1)^m$ has derivative $m2^mr^{m-1}$ for $r<\frac12$ and $0$ for $r>\frac12$. For $K=\T^m$ we have $d_n(x,K)=0$, so the decay condition (2.9) is trivial ($C_K=c_K=\alpha_K=1$), and $I_{K,n}=n$, so $I_{K,n}/\Var[H_n]\le1/c$.

\medskip\noindent\textbf{Step 4: conclusion.} Theorem~\ref{tl:H} gives
\[\dK\Bigg(\frac{H_n-\E [H_n]}{\sqrt{\Var [H_n]}},\mathcal{N}\Bigg)\le\frac{C}{\sqrt{Var[H_n]}}\le C n^{-1/2}.\]
 \hfill$\square$

\section{Proofs for the very sparse regime}\label{sec:vsparse}
\noindent Throughout $n(2r_n)^m\le\frac1{2e}$; then $(2r_n)^m\le\frac1{2en}$, so $r_n\le\frac12(2en)^{-1/m}<\frac16$ for $n\ge n_0(m)$ and all the results of Sections~\ref{sec:proofmean} and~\ref{sec:proofs} for $r_n<\frac16$ apply. $C$ depends on $m$ only.

In this regime the random graph consists mostly of isolated vertices, edges and small trees, so $n\mathcal F_n$ is essentially the number $N_\vee$ of induced two-paths. The section makes this precise and derives the limit laws in four steps.
\begin{itemize}[leftmargin=*]
\item Section~\ref{sec:vs1} collects the combinatorial inputs.  Lemma~\ref{lem:comp} shows that $n\mathcal F_n=N_\vee+\Xi$, where $\Xi$ comes only from components with at least four vertices, and Lemma~\ref{lem:csize} shows that such components are rare.
\item Section~\ref{sec:vs2} approximates $N_\vee$ by a Poisson variable with the Chen--Stein method (Theorem~\ref{tl:I} in Section~\ref{ts:I}). Lemma~\ref{lem:meanNvee} gives the exact mean $\lambda_n$ and Lemma~\ref{lem:b} the dependence bounds $b_1,b_2,b_3$ that the method requires.
\item Section~\ref{sec:vs3}, through Lemma~\ref{lem:VarD}, bounds the variance of $\Xi$ and so shows that the larger components are negligible.
\item Section~\ref{sec:vs4} combines these ingredients into the proof of Theorem~\ref{thm:vsparse}. For the Gaussian limit with a rate it also uses the Berry--Esseen theorem for i.i.d.\ sums (Theorem~\ref{tl:J}) to prove Lemma~\ref{lem:poisBE} (a Berry--Esseen bound for the Poisson law), and the perturbation lemma (Lemma~\ref{lem:pert} in Section~\ref{ts:lem:pert}) to absorb $\Xi$.
\end{itemize}

\subsection{Components of size at most three}\label{sec:vs1}

\begin{lemma}[Components]\label{lem:comp}\tlsec{lem:comp}
For a connected graph $\mathcal C$ let \[f(\mathcal C)=\sum_{\{i,j\}\in E(\mathcal C)}\frac{(d_i-d_j)^2}{d_id_j}\]
and $N_\vee(\mathcal C)$ the number of induced two-paths (triples $\{i,j,k\}$ with exactly the two edges $ij,jk$) in $\mathcal C$. Then
\begin{enumerate}[label=(\alph*),leftmargin=*]
\item $f(\mathcal C)=N_\vee(\mathcal C)$ if $|\mathcal C|\le3$;
\item $0\le f(\mathcal C),N_\vee(\mathcal C)\le|\mathcal C|^3$ for every $\mathcal C$.
\end{enumerate}
Consequently $n\mathcal{F}_n=\sum_{\mathcal C}f(\mathcal C)$ and $N_\vee=\sum_{\mathcal C}N_\vee(\mathcal C)$, and
\Line{\Xi:=n\mathcal{F}_n-N_\vee=\sum_{|\mathcal C|\ge4}\big(f(\mathcal C)-N_\vee(\mathcal C)\big),}
which satisfies
\Line{|\Xi|\le\sum_{v=1}^n|\mathcal C(v)|^2 \1\{|\mathcal C(v)|\ge4\},\qquad\mathcal C(v)=\text{component of }v .}
In particular, let $\Omega_4$ be the event that the graph has a connected component with at least four vertices; then $n\mathcal{F}_n=N_\vee$ on $\Omega_4^c$.
\end{lemma}
\noindent The proof of Lemma~\ref{lem:comp} is given in Section~\ref{pf:lem:comp} of the Supplement.

\begin{lemma}[Size of a component]\label{lem:csize}\tlsec{lem:csize}
For every integer $k\ge2$, \Line{\P(|\mathcal C(1)|\ge k)\le e (en(2r_n)^m)^{k-1}.} Consequently, for $p\ge0$, $n(2r_n)^m\le\frac1{2e}$, and $\Omega_4$ given in Lemma~\ref{lem:comp},
\Line{\E\big[|\mathcal C(1)|^p\1\{|\mathcal C(1)|\ge4\}\big]\le C_p\big(n(2r_n)^m\big)^{3},\qquad\P(\Omega_4)\le\frac23 n^{4}(2r_n)^{3m} .}
\end{lemma}
\noindent The proof of Lemma~\ref{lem:csize} is given in Section~\ref{pf:lem:csize} of the Supplement.

\subsection{The count of induced two-paths}\label{sec:vs2}

\noindent The Chen--Stein method bounds the total-variation distance between a sum of rare, locally dependent indicators and a Poisson variable by the quantities $b_1,b_2,b_3$ below. In the very sparse regime $n\mathcal F_n$ is essentially the number of induced two-paths, so this theorem gives the Poisson limit and, via Lemma~\ref{lem:poisBE} (Section~\ref{ts:lem:poisBE}), the Gaussian limit.
\begin{tool}[Poisson approximation by the Chen--Stein method; Arratia--Goldstein--Gordon 1989 \cite{AGG89}, Thm.~1; see also \cite{BHJ92}]\label{tl:I}\tlsec{I}
Let $(I_\alpha)_{\alpha\in \mathcal I}$ be a finite family of Bernoulli variables, $p_\alpha=\P(I_\alpha=1)$, $W=\sum_\alpha I_\alpha$, $\lambda=\E [W]$. For each $\alpha$ choose a set $B_\alpha\subseteq \mathcal I$ with $\alpha\in B_\alpha$ (the ``neighbourhood of dependence'') and put
\[\begin{aligned}&b_1=\sum_{\alpha\in \mathcal I}\sum_{\beta\in B_\alpha}p_\alpha p_\beta,\qquad b_2=\sum_{\alpha\in \mathcal I}\sum_{\beta\in B_\alpha\setminus\{\alpha\}}\E[I_\alpha I_\beta],\\
&b_3=\sum_{\alpha\in \mathcal I}\E\big[\Big|\E\big[I_\alpha-p_\alpha \big| \sigma(I_\beta:\beta\notin B_\alpha)\big]\Big|\big] .\end{aligned}\]
Then, for the Poisson law $\mathrm{Po}(\lambda)$,
\Line{d_{\mathrm{TV}}\big(W,\mathrm{Po}(\lambda)\big)\le2\Big[(b_1+b_2)\frac{1-e^{-\lambda}}{\lambda}+b_3 \min\big(1,1.4\lambda^{-1/2}\big)\Big]\le2 (b_1+b_2+b_3).}
In particular, if $I_\alpha$ is independent of $\sigma(I_\beta:\beta\notin B_\alpha)$ for all $\alpha$ then $b_3=0$.
\end{tool}

Index the induced two-paths by $\alpha=(j;\{i,k\})$ (centre $j$, ends $i,k$), \[I_\alpha=A_{ji}A_{jk}(1-A_{ik}),\]
vertex set $V_\alpha=\{i,j,k\}$, $\mathcal I$ the set of all such $\alpha$ ($|\mathcal I|=n\binom{n-1}2\le n^3/2$), $p_\alpha=\E [I_\alpha]$.

\begin{lemma}[Exact mean]\label{lem:meanNvee}\tlsec{lem:meanNvee}
For every $\alpha$, $p_\alpha=\kappa_m (2r_n)^{2m}$, hence \Line{\lambda_n=\E [N_\vee]=\frac{\kappa_m}2n(n-1)(n-2)(2r_n)^{2m}.}
\end{lemma}
\noindent The proof of Lemma~\ref{lem:meanNvee} is given in Section~\ref{pf:lem:meanNvee} of the Supplement.

\begin{lemma}[Neighbourhoods of dependence]\label{lem:b}\tlsec{lem:b}
Let $B_\alpha=\{\beta\in \mathcal I:V_\beta\cap V_\alpha\ne\emptyset\}$ (which contains $\alpha$), and let $b_1,b_2,b_3$ be the Chen--Stein quantities of Theorem~\ref{tl:I} for this choice of $B_\alpha$. 
Then $I_\alpha$ is independent of $(I_\beta)_{\beta\notin B_\alpha}$, so $b_3=0$, and for $n(2r_n)^m\le1$
\[b_1\le\tfrac52 n^{5}(2r_n)^{4m},\qquad b_2\le\tfrac52 n^{5}(2r_n)^{4m}+\tfrac92 n^{4}(2r_n)^{3m},\quad b_1+b_2\le10 n^{4}(2r_n)^{3m},\]
\[\big|\mathrm{Var}[N_\vee]-\lambda_n\big|\le b_1+b_2\le10 n^{4}(2r_n)^{3m} .\]
\end{lemma}
\noindent The proof of Lemma~\ref{lem:b} is given in Section~\ref{pf:lem:b} of the Supplement.

\subsection{The contribution of larger components}\label{sec:vs3}

\begin{lemma}[Variance of $\Xi$]\label{lem:VarD}\tlsec{lem:VarD}
For $n(2r_n)^m\le\frac1{2e}$, $\Var [\Xi]\le C n^{4}(2r_n)^{3m}$.
\end{lemma}
\noindent The proof of Lemma~\ref{lem:VarD} is given in Section~\ref{pf:lem:VarD} of the Supplement.

\subsection{Proof of Theorem~\ref{thm:vsparse}}\label{sec:vs4}
Throughout, $n(2r_n)^m\le\frac1{2e}$, $\kappa_m=1-(\frac34)^m\in[\frac14,1)$ and \Line{\lambda_n=\mathbb{E}[N_\vee]=\frac{\kappa_m}2n(n-1)(n-2)(2r_n)^{2m}.}

\begin{lemma}[Berry--Esseen for the Poisson law]\label{lem:poisBE}\tlsec{lem:poisBE}
There is an absolute constant $C_{\mathrm P}$ such that for all $\lambda>0$,
\Line{\dK\Big(\frac{\mathrm{Po}(\lambda)-\lambda}{\sqrt\lambda},\mathcal{N}\Big)\le C_{\mathrm P} \lambda^{-1/2}.}
\end{lemma}
\noindent The proof of Lemma~\ref{lem:poisBE} is based on the classical Berry--Esseen theorem for i.i.d.  sums and given in Section~\ref{pf:lem:poisBE} of the Supplement.

\medskip\noindent\emph{(i) Too sparse.} By Lemma \ref{lem:comp}, if $\mathcal{F}_n\ne0$, then $N_\vee\ge1$ or $\Omega_4$ occurs. By Markov's inequality, \[\P(N_\vee\ge1)\le\lambda_n\le\frac12n^3(2r_n)^{2m}.\]
Then by the union bound and Lemma \ref{lem:csize},
\Line{\P(\mathcal{F}_n\ne0)\le n^3(2r_n)^{2m}\big(\tfrac12+\tfrac23n(2r_n)^m\big)\le n^3(2r_n)^{2m}.}

\medskip\noindent\emph{(ii) Poisson regime.} Note that $N_\vee=\sum_\alpha I_\alpha$ over triples with a designated centre. By Lemma~\ref{lem:b} (Section~\ref{ts:lem:b}) and Theorem~\ref{tl:I} (Section~\ref{ts:I}) (with $b_3=0$), \Line{d_{\mathrm{TV}}(N_\vee,\mathrm{Po}(\lambda_n))\le2(b_1+b_2)\le20 n^{4}(2r_n)^{3m}.} Since $n\mathcal{F}_n=N_\vee$ off $\Omega_4$ by Lemma \ref{lem:comp}, for every set $\mathcal B$ the probabilities $\P(n\mathcal F_n\in\mathcal B)$ and $\P(N_\vee\in\mathcal B)$ agree on $\Omega_4^c$ and differ on $\Omega_4$ by at most $\P(\Omega_4)$; it then follows from Lemma \ref{lem:csize} that \Line{d_{\mathrm{TV}}(n\mathcal{F}_n,N_\vee)\le\P(\Omega_4)\le\frac23n^4(2r_n)^{3m}.} 
Then by triangle inequality, we have
\Line{d_{\mathrm{TV}}\big(n\mathcal{F}_n,\mathrm{Po}(\lambda_n)\big)\le21 n^{4}(2r_n)^{3m}.}

\medskip\noindent\emph{(iii) Variance and normal regime.} It is clear that $\lambda_n\asymp n^{3}(2r_n)^{2m}$ by Lemma~\ref{lem:meanNvee}.  Since \[n\mathcal{F}_n=N_\vee+\Xi\]
by Lemma \ref{lem:comp}, Cauchy--Schwarz inequality, Lemma~\ref{lem:b} (Section~\ref{ts:lem:b}) and Lemma~\ref{lem:VarD} (Section~\ref{ts:lem:VarD}) give
\[\big|\Var[n\mathcal{F}_n]-\lambda_n\big|\le|\Var [N_\vee]-\lambda_n|+2\sqrt{\Var [N_\vee]\Var [\Xi]}+\Var [\Xi]\le Cn\big(n(2r_n)^m\big)^{5/2}.\]
 Dividing by $\lambda_n$ gives \Line{\Var[n\mathcal{F}_n]=\lambda_n\big(1+O((n(2r_n)^m)^{1/2})\big).} 

For the normal approximation we may assume $n^3(2r_n)^{2m}\ge1$ (otherwise the bound exceeds $1\ge\dK$), so $\lambda_n\ge\frac1{16}$. Put \Line{T=n\mathcal{F}_n-\E[n\mathcal{F}_n]=(N_\vee-\lambda_n)+(\Xi-\E\Xi).}

\smallskip\noindent\emph{Step 1.} Theorem~\ref{tl:I} (Section~\ref{ts:I}) in its refined form carries the factor $1/\lambda_n$. By Lemma~\ref{lem:b},
\Line{d_{\mathrm{TV}}(N_\vee,\mathrm{Po}(\lambda_n))\le\frac{2(b_1+b_2)}{\lambda_n}\le320 n(2r_n)^m.} With $X=(N_\vee-\lambda_n)/\sqrt{\lambda_n}$ and Lemma~\ref{lem:poisBE} (Section~\ref{ts:lem:poisBE}),
\Line{\dK(X,\mathcal{N})\le Cn(2r_n)^m+C_{\mathrm P}\lambda_n^{-1/2}.}

\smallskip\noindent\emph{Step 2.} Let $Y=(\Xi-\E\Xi)/\sqrt{\lambda_n}$, so $T/\sqrt{\lambda_n}=X+Y$ and \Line{\E[Y^2]=\Var[\Xi]/\lambda_n\le Cn(2r_n)^m} by Lemma~\ref{lem:VarD} (Section~\ref{ts:lem:VarD}). Lemma~\ref{lem:pert} (Section~\ref{ts:lem:pert}) gives
\Line{\dK(X+Y,\mathcal{N})\le\dK(X,\mathcal{N})+C(\E [Y^2])^{1/3}\le C\big(\big(n(2r_n)^m\big)^{1/3}+(n^3(2r_n)^{2m})^{-1/2}\big),}
which is the bound in (iii). \hfill$\blacksquare$

\clearpage
\appendix
\begin{center}\textbf{\large Supplement to ``Limit laws for the friendship paradox in random geometric graphs''}\end{center}
\renewcommand{\thesection}{S}\setcounter{section}{1}\setcounter{subsection}{0}\setcounter{equation}{0}
\renewcommand{\thesubsection}{S\arabic{subsection}}\renewcommand{\theHsubsection}{S\arabic{subsection}}\renewcommand{\theHsection}{S}
\renewcommand{\theequation}{S.\arabic{equation}}

\begin{aug}
\author[A]{\fnms{Mingao}~\snm{Yuan}\ead[label=e1]{myuan2@utep.edu}}
\address[A]{Department of Mathematical Sciences, The University of Texas at El Paso, El Paso, Texas 79968, USA\printead[presep={,\ }]{e1}}
\end{aug}

\medskip

\noindent This supplement to the paper ``Limit laws for the friendship paradox in random geometric graphs'' contains the proofs of the lemmas and propositions stated in Sections~\ref{sec:proofmean}--\ref{sec:vsparse} of the paper, in the order in which they are stated; each lemma or proposition is copied before its proof. Numbers of theorems, lemmas, propositions, sections and equations that are not in this supplement refer to the paper; Theorems~A--K are the published results quoted there.

\subsection{Proof of Lemma~\ref{lem:kappa}}\label{pf:lem:kappa}
\mbox{}\par\begin{restated}{Lemma}{lem:kappa}{ (The triangle probability and $\kappa(r_n)$ for all $r_n<\frac12$)} For all $r_n<\frac12$,
\Line{\E[A_{12}A_{13}A_{23}]=\begin{cases} (3r_n^2)^m,&0<r_n\le\dfrac13,\\[2mm](12r_n^2-6r_n+1)^m ,&\dfrac13<r_n<\dfrac12.\end{cases}}
Then $\E[A_{12}A_{13}A_{23}]=(2r_n)^{2m}\rho(r_n)^m$,  
with $\rho(t)$ as in Theorem~\ref{thm:mean} (Section~\ref{ts:thm:mean}), so that $\kappa(r_n)=1-\rho(r_n)^m$. 
\end{restated}
\begin{proof}\smallskip\noindent\textbf{Step 1: reduction to the circle.} The cube indicator factorises, $\beta(x)=\prod_{b=1}^m\beta^{(1)}(x_b)$, where $\beta^{(1)}$ is the indicator of the arc $(-r_n,r_n)$ of the circle $\mathbb R/\mathbb Z$. The $L_\infty$ distance satisfies \Line{\{d(y,z)<r_n\}=\bigcap_b\{d_1(y_b,z_b)<r_n\},} where $d_1$ is the circular distance. Hence
\[
\E[A_{12}A_{13}A_{23}]=\iint\beta(y)\beta(z)\1\{d(y,z)<r_n\} dy dz=\vartheta^m,
\]
where
 \Line{\vartheta:=\big|\{(y,z)\in(-r_n,r_n)^2:\ d_1(y,z)<r_n\}\big|} is the area of a planar set.

\smallskip\noindent\textbf{Step 2: the area $\vartheta$.} The square $(-r_n,r_n)^2$ has area $4r_n^2$. Let $S=(-r_n,r_n)^2$. The area of the region outside the plannar set is
\begin{align*}
    \iint_{S}\1\{d_1(y,z)\geq r_n\}dydz&= \iint_{S}\1\{|y-z|\geq r_n\}\1\{1-|y-z|\geq r_n\}dydz.
\end{align*}
Let $y=y\in(-r_n,r_n)$, $u=y-z\in(-2r_n,2r_n)$. Note that $z=y-u\in(-r_n,r_n)$. Then $y\in (-r_n,r_n)\cap (u-r_n,u+r_n)$. The length of this interval is $(2r_n-|u|)$  (because $|u|<2r_n$). Then
\begin{align*}
    \iint_{S}\1\{d_1(y,z)\geq r_n\}dydz&= \int_{-2r_n}^{2r_n}\1\{|u|\geq r_n\}\1\{1-|u|\geq r_n\}\int_{(-r_n,r_n)\cap (u-r_n,u+r_n)}dydu\\
    &=\int_{-2r_n}^{2r_n}\1\{r_n\leq |u|\leq 1-r_n\}(2r_n-|u|)du.
\end{align*}
There are two cases.
\begin{itemize}[leftmargin=*,itemsep=1pt]
\item $r_n\le\frac13$: then $1-r_n\ge2r_n$, the excluded area is \[2\int_{r_n}^{2r_n}(2r_n-u) du=r_n^2,\]
hence $\vartheta=3r_n^2$.
\item $\frac13<r_n<\frac12$:  then $1-r_n<2r_n$, the excluded area is \Line{2\int_{r_n}^{1-r_n}(2r_n-u) du=6r_n-8r_n^2-1.}
\end{itemize}
Combining both cases yields the results of the lemma.
\end{proof}

\subsection{Proof of Lemma~\ref{lem:tree}}\label{pf:lem:tree}
\mbox{}\par\begin{restated}{Lemma}{lem:tree}{ (Trees)}
Let $\mathcal T$ be a tree whose vertex set is a set $V\subseteq\{1,\dots,n\}$ of $v$ vertices and whose $v-1$ edges are pairs $\{u,w\}\subset V$. Then, exactly,
\Line{\P\big(A_{uw}=1\ \text{for all edges }\{u,w\}\text{ of }\mathcal T\big)=(2r_n)^{m(v-1)}.}
\end{restated}
\begin{proof}Induction on $v$. For $v=1$ there is nothing to show (empty product, probability $1$). For $v\ge2$ the tree has a leaf $\ell$ with a unique neighbour $p$ in $\mathcal T$. The event for the smaller tree $\mathcal T'=\mathcal T-\ell$ depends only on $(X_u)_{u\in V\setminus\{\ell\}}$, and $X_\ell$ is independent of these. Conditioning on all $X_u$, $u\ne\ell$: \Line{\P(A_{\ell p}=1\mid(X_u)_{u\ne\ell})=\mathrm{vol} \mathcal A(X_p)=(2r_n)^m.} Therefore the probability is \Line{(2r_n)^m\cdot\P(\text{event for }\mathcal T')=(2r_n)^m\cdot (2r_n)^{m(v-2)}.}
\end{proof}

\subsection{Proof of Proposition~\ref{prop:id}}\label{pf:prop:id}
\mbox{}\par\begin{restated}{Proposition}{prop:id}{ (Identity)}
In every graph (with $\Delta_i=0$ for isolated $i$),\Line{\mathcal{F}_n=\frac1n\sum_{i<j}A_{ij}\frac{(d_i-d_j)^2}{d_id_j}\ge0 .}
\end{restated}
\begin{proof}If $d_i\ge1$, using $\sum_jA_{ij}=d_i$,
\Line{\Delta_i=\frac1{d_i}\sum_jA_{ij}d_j-d_i=\sum_jA_{ij}\Big(\frac{d_j}{d_i}-1\Big).}
If $d_i=0$ then $\Delta_i=0$ and the right-hand side is the empty sum $0$ as well, so the display holds for all $i$ (all terms with $A_{ij}=1$ have $d_i\ge1$). Summing over $i$ and grouping the two ordered pairs $(i,j)$ and $(j,i)$ of each edge $\{i,j\}$, $i<j$:
\Line{\sum_i\Delta_i=\sum_{i<j}A_{ij}\Big[\Big(\frac{d_j}{d_i}-1\Big)+\Big(\frac{d_i}{d_j}-1\Big)\Big]=\sum_{i<j}A_{ij} \frac{d_j^2+d_i^2-2d_id_j}{d_id_j}=\sum_{i<j}A_{ij}\frac{(d_i-d_j)^2}{d_id_j}.}
Divide by $n$. Every term is $\ge0$ (for an edge both degrees are $\ge1$).
\end{proof}

\subsection{Proof of Lemma~\ref{lem:S}}\label{pf:lem:S}
\mbox{}\par\begin{restated}{Lemma}{lem:S}{}
$\mathcal{S}=\sum_id_i(d_i-1)-6\mathcal T$, where $\mathcal T$ is the number of triangles of the graph.
\end{restated}
\begin{proof}For an edge $\{i,j\}$, \Line{(A_{ik}-A_{jk})^2=A_{ik}+A_{jk}-2A_{ik}A_{jk}} (check the four cases of $(A_{ik},A_{jk})$). Summing over $k\ne i,j$ and writing $\gamma_{ij}=\sum_{k\ne i,j}A_{ik}A_{jk}$ for the number of common neighbours,
\Line{\sum_{k\ne i,j}(A_{ik}-A_{jk})^2=(d_i-1)+(d_j-1)-2\gamma_{ij}\qquad(A_{ij}=1).}
Therefore \Line{\mathcal{S}=\sum_{i<j}A_{ij}(d_i+d_j)-2\sum_{i<j}A_{ij}-2\sum_{i<j}A_{ij}\gamma_{ij}.} Now (a) $\sum_{i<j}A_{ij}(d_i+d_j)=\sum_id_i\cdot d_i$ because vertex $i$ occurs in exactly $d_i$ edges; (b) $\sum_{i<j}A_{ij}=\frac12\sum_id_i$; (c) $\sum_{i<j}A_{ij}\gamma_{ij}=3\mathcal T$ because each triangle is counted once for each of its three edges. So $\mathcal{S}=\sum_id_i^2-\sum_id_i-6\mathcal T$.
\end{proof}

\subsection{Proof of Lemma~\ref{lem:ECn}}\label{pf:lem:ECn}
\mbox{}\par\begin{restated}{Lemma}{lem:ECn}{}
$\E [\Upsilon_n]=0$, and every summand \Line{H_{ij;kl}=A_{ij}(A_{ik}-A_{jk})(A_{il}-A_{jl})} of $\Upsilon_n$ satisfies $\E[H_{ij;kl}\mid X_U]=0$ whenever $k\notin U$ or $l\notin U$, $U\subseteq\{i,j,k,l\}$.
\end{restated}
\begin{proof}\smallskip\noindent\textbf{Step 1: the case $k\notin U$.} Let $U'=\{i,j,l\}\cup U$; it does not contain $k$, because $k\notin U$ and $k\notin\{i,j,l\}$. Since $A_{ik}$ is a function of $(X_i,X_k)$ and $X_k$ is independent of $(X_i,X_j,X_l)$,
\[
\E[A_{ik}\mid X_i,X_j,X_l]=\E[A_{ik}\mid X_i]=(2r_n)^m,\qquad\text{likewise for }A_{jk},
\]
so $\E[A_{ik}-A_{jk}\mid X_i,X_j,X_l]=0$.

\smallskip\noindent\textbf{Step 2: the summand.} $A_{ij}$ and $A_{il}-A_{jl}$ are functions of $(X_i,X_j,X_l)$, so they factor out:
\[
\E[H_{ij;kl}\mid X_i,X_j,X_l]=A_{ij}(A_{il}-A_{jl}) \E[A_{ik}-A_{jk}\mid X_i,X_j,X_l]=0 .
\]

\smallskip\noindent\textbf{Step 3: conclusion.} By the tower property, \Line{\E[H_{ij;kl}\mid X_U]=\E\big[\E[H_{ij;kl}\mid X_i,X_j,X_l,X_U]\mid X_U\big]=0;} conditioning first on the larger set $\{i,j,l\}\cup U$, which still misses $k$, gives the same zero. The case $l\notin U$ is identical by symmetry, and $U=\emptyset$ gives $\E[H]=0$.
\end{proof}

\subsection{Proof of Proposition~\ref{prop:alpha}}\label{pf:prop:alpha}
\mbox{}\par\begin{restated}{Proposition}{prop:alpha}{ (The constants $\alpha_1,\alpha_2,\alpha_3$)}
Let $0<r<\frac12$, $r_n=r$, and let $X_k,X_l$ be independent and uniform on $\T^m$, with $P_2=P_2(X_k,X_l)$ and $P_3=P_3(X_k,X_l)$. Let $\alpha_1,\alpha_2,\alpha_3$ be the piecewise polynomials~\eqref{eq:alpha1}--\eqref{eq:alpha3}.
\[\E[P_2^2]=\alpha_1(r)^m,\ \ \ \ \ \ \ \  \E[P_2P_3]=\alpha_2(r)^m,\ \ \ \ \ \ \ \   \E[P_3^2]=\alpha_3(r)^m,\] 
and
\begin{align*}
\E[g(X_k,X_l)^2]=4\big[(2r)^{2m}\alpha_1(r)^m-2(2r)^m\alpha_2(r)^m+\alpha_3(r)^m\big].
\end{align*}
\end{restated}
\begin{proof}\emph{Factorisation.}
Since $\beta(x)=\prod_{b=1}^m\1_{(-r,r)}(x_b)$ (the arc $(-r,r)$ of $\mathbb R/\mathbb Z$), the integrals defining $P_2,P_3$ factorise over the coordinates: \Line{\begin{gathered}P_2(x)=P_2(x,0)=\prod_b\PtD{x_b},\\ P_3(x)=P_3(x,0)=\prod_b\PthD{x_b},\end{gathered}} where $\E_1$ is the expectation of the one-dimensional model: $m=1$, points $X_1, X_2, \dots,X_n$ uniform on the circle $\mathbb R/\mathbb Z$, $A_{uv}=\1\{d(X_u,X_v)<r\}$. Write \Line{\begin{gathered}\tilde\alpha_1=\E_1[(\PtD{U})^2],\\ \tilde\alpha_2=\E_1[\PtD{U}\PthD{U}],\\ \tilde\alpha_3=\E_1[(\PthD{U})^2]\end{gathered}} with $U$ uniform on $\mathbb R/\mathbb Z$.

Steps~1--3 below show that $\tilde\alpha_i=\alpha_i$ (the polynomials~\eqref{eq:alpha1}--\eqref{eq:alpha3}). Since $X_k,X_l$ are independent with $X_k-X_l$ uniform, \Line{\E[P_2^2]=\prod_b\E_1[(\PtD{U_b})^2]=\tilde\alpha_1^m=\alpha_1(r)^m,} and likewise \[\E[P_2P_3]=\tilde\alpha_2^m=\alpha_2(r)^m,\ \ \ \ \ \ \ \ \ \E[P_3^2]=\tilde\alpha_3^m=\alpha_3(r)^m,\]
where $U_1,\dots,U_m$ are i.i.d.\ uniform on $[0,1)$. 

\emph{Density of sum of uniform variables.} Fix $0<r<\frac12$. If $U_1,\dots,U_n$ are i.i.d.\ uniform on $(0,1)$, then $U_1+\dots+U_n$ has the Irwin--Hall density \Line{f_n(y)=\frac1{(n-1)!}\sum_{i=0}^n(-1)^i\binom ni(y-i)_+^{n-1},\ \ \ \ 0<y<n.} 
Let $V_1,\dots,V_n$ be i.i.d.\ uniform on $(-r,r)\subset\mathbb R$, let $S_n=V_1+\dots+V_n$ with density $p_n(x)$ on $\mathbb R$. Since $V_k=2rU_k-r$ in law, \[p_n(x)=\frac1{2r}f_n\big(\frac{x+nr}{2r}\big)=\frac1{(2r)^n(n-1)!}\sum_{i=0}^n(-1)^i\binom ni(x+(n-2i)r)_+^{n-1},\] and $p_n(x)=0$ for $|x|\ge nr$. 
 The density of $S_n$ modulo $1$ on $\mathbb R/\mathbb Z$ is
\begin{equation}\label{exprqnx}
q_n(x)=\sum_{j\in\mathbb Z}p_n(x+j).  
\end{equation}

\emph{Step 1: reduction to $q_n(0)$.} Let $I=(-r,r)$ and $\1_I(x)=\1\{x\in I\ (\text{mod}\ 1)\}$ be the 1-periodic indicator function of $I$. Clearly, $\1_I(x)=\1_I(-x)$. Let $X_1,\dots,X_n$ be i.i.d. uniform on $\R/\mathbb{Z}\cong[0,1)$.  Recall that
\[
A_{uv}=\1\{d(X_u,X_v)<r\},\qquad d(a,b)=\min\{|a-b|,1-|a-b|\}=\min_{k\in\mathbb{Z}}|a-b+k|.
\]
Since $d(X_u,X_v)<r$ is equivalent to $X_u-X_v\in I\ (\text{mod}\ 1)$, then 
\[A_{uv}=\1_I(X_u-X_v).\]
It then follows that
\[\begin{aligned}\PtD{x}=\E\big[\1_{I}(U_1)\1_{I}(x-U_1)\big]=(\1_I*\1_I)(x). \end{aligned}\]
Here, $\1_I*\1_I$ is the convolution of the indicator function $\1_I(x)$.
Note that the density of $V_1+V_2$ meodulo $1$ at $x$ is $(2r)^{-2}\1_I*\1_I(x)$. It is equal to $q_2(x)$ in \eqref{exprqnx}. Then,
\Line{\PtD{x}=(2r)^2q_2(x).}
Similarly, for $P_3$,
\[\PthD{x}=(2r)^3q_3(x).\]

Let $U$ be uniform on $\mathbb R/\mathbb Z$. Since $q_2$ and $q_3$ are even, this gives
\[\begin{aligned}\alpha_1&=\E\big[\big(\PtD{U}\big)^2\big]\\&=(2r)^4 \E[q_2(U)q_2(-U)]\\
&=(2r)^4q_2*q_2(0)\\
&=(2r)^4q_4(0), 
\end{aligned}
\]
\[
\begin{aligned}
\alpha_2&=\E\big[\PtD{U} \PthD{U}\big]=(2r)^5q_5(0),\\
\alpha_3&=\E\big[\big(\PthD{U}\big)^2\big]=(2r)^6q_6(0).\end{aligned}\]

\emph{Step 2: an explicit formula.} With $\alpha_1,\alpha_2,\alpha_3$ corresponding to $n=4,5,6$,  it follows from \eqref{exprqnx} that
\LineE{eq:alphaformula}{(2r)^nq_n(0)=(2r)^n\sum_{j\in\mathbb Z}p_n(j)=\frac1{(n-1)!}\sum_{|j|<nr}\ \sum_{i=0}^n(-1)^i\binom ni\big(j+(n-2i)r\big)_+^{n-1}.}

\emph{Step 3: evaluation.} The term with indices $(j,i)$ in~\eqref{eq:alphaformula} is present iff $|j|<nr$ and $j+(n-2i)r>0$. The set of present terms changes only when $r=|j|/|n-2i|$ or $r=|j|/n$. In $(0,\frac12)$ these values are $\frac14$ for $n=4$; $\frac15,\frac13,\frac25$ for $n=5$; $\frac16,\frac14,\frac13$ for $n=6$. In each open interval between consecutive breakpoints the set of present terms is constant, so~\eqref{eq:alphaformula} is a polynomial in $r$. For each $j$ the present indices are $i=0,\dots,i_{\max}(j)$, where $i_{\max}(j)$ is as follows ($-$ means that no term with this $j$ is present):
\begin{center}\small\renewcommand{\arraystretch}{1.15}
\begin{adjustbox}{max width=\linewidth}\begin{tabular}{@{}cccccccc@{}}
\toprule
$n$ & interval & $j=-2$ & $j=-1$ & $j=0$ & $j=1$ & $j=2$ & polynomial\\
\midrule
$4$ & $(0,\frac14]$ & $-$ & $-$ & $1$ & $-$ & $-$ & $\frac{16}3r^3$\\
$4$ & $[\frac14,\frac12)$ & $-$ & $0$ & $1$ & $3$ & $-$ & $-\frac13+4r-16r^2+\frac{80}3r^3$\\
$5$ & $(0,\frac15]$ & $-$ & $-$ & $2$ & $-$ & $-$ & $\frac{115}{12}r^4$\\
$5$ & $[\frac15,\frac13]$ & $-$ & $0$ & $2$ & $4$ & $-$ & $\frac1{12}-\frac53r+\frac{25}2r^2-\frac{125}3r^3+\frac{185}3r^4$\\
$5$ & $[\frac13,\frac25]$ & $-$ & $1$ & $2$ & $3$ & $-$ & $-\frac13+\frac{10}3r-10r^2+\frac{10}3r^3+\frac{335}{12}r^4$\\
$5$ & $[\frac25,\frac12)$ & $0$ & $1$ & $2$ & $3$ & $4$ & $1-10r+40r^2-80r^3+80r^4$\\
$6$ & $(0,\frac16]$ & $-$ & $-$ & $2$ & $-$ & $-$ & $\frac{88}5r^5$\\
$6$ & $[\frac16,\frac14]$ & $-$ & $0$ & $2$ & $5$ & $-$ & $-\frac1{60}+\frac r2-6r^2+36r^3-108r^4+\frac{736}5r^5$\\
$6$ & $[\frac14,\frac13]$ & $-$ & $1$ & $2$ & $4$ & $-$ & $\frac1{12}-\frac32r+10r^2-28r^3+20r^4+\frac{224}5r^5$\\
$6$ & $[\frac13,\frac12)$ & $0$ & $1$ & $2$ & $4$ & $5$ & $-\frac9{20}+\frac{13}2r-38r^2+116r^3-196r^4+\frac{872}5r^5$\\
\bottomrule
\end{tabular}\end{adjustbox}
\end{center}
\emph{Expansion in every interval.} On an interval with present indices $i\le i_{\max}(j)$ put \Line{S_j(r)=\sum_{i=0}^{i_{\max}(j)}(-1)^i\binom ni\big(j+(n-2i)r\big)^{n-1},} so that~\eqref{eq:alphaformula} equals $\frac1{(n-1)!}\sum_jS_j(r)$. The subtotals and the resulting polynomials are the following (each $S_j$ is a polynomial in $r$ of degree $n-1$):
\medskip

\noindent\emph{$n=4$, $r\in(0,\frac{1}{4}]$.}\[\begin{aligned}S_{0}&=32r^{3},\hskip 1cm \tfrac1{6}\textstyle\sum_jS_j=\frac{16}{3}r^{3}.\end{aligned}\]
\noindent\emph{$n=4$, $r\in[\frac{1}{4},\frac{1}{2})$.}\[\begin{aligned}S_{-1}&=-1 + 12r - 48r^{2} + 64r^{3},\hskip 1cm S_{0} =32r^{3},\hskip 1cmS_{1} =-1 + 12r - 48r^{2} + 64r^{3},\end{aligned}\]
\[ \tfrac1{6}\textstyle\sum_jS_j=-\frac{1}{3} + 4r - 16r^{2} + \frac{80}{3}r^{3}.\]
\noindent\emph{$n=5$, $r\in(0,\frac{1}{5}]$.}\[\begin{aligned}S_{0}=230r^{4},\hskip 1cm \tfrac1{24}\textstyle\sum_jS_j=\frac{115}{12}r^{4}.\end{aligned}\]
\noindent\emph{$n=5$, $r\in[\frac{1}{5},\frac{1}{3}]$.}\[\begin{aligned}S_{-1}&=1 - 20r + 150r^{2} - 500r^{3} + 625r^{4},\hskip 1cm S_{0}=230r^{4},\\S_{1}&=1 - 20r + 150r^{2} - 500r^{3} + 625r^{4},\\\tfrac1{24}\textstyle\sum_jS_j&=\frac{1}{12} - \frac{5}{3}r + \frac{25}{2}r^{2} - \frac{125}{3}r^{3} + \frac{185}{3}r^{4}.\end{aligned}\]
\noindent\emph{$n=5$, $r\in[\frac{1}{3},\frac{2}{5}]$.}\[\begin{aligned}S_{-1}&=-4 + 40r - 120r^{2} + 40r^{3} + 220r^{4},\hskip 1cm S_{0}=230r^{4},\\S_{1}&=-4 + 40r - 120r^{2} + 40r^{3} + 220r^{4},\\\tfrac1{24}\textstyle\sum_jS_j&=-\frac{1}{3} + \frac{10}{3}r - 10r^{2} + \frac{10}{3}r^{3} + \frac{335}{12}r^{4}.\end{aligned}\]
\noindent\emph{$n=5$, $r\in[\frac{2}{5},\frac{1}{2})$.}\[\begin{aligned}S_{-2}&=16 - 160r + 600r^{2} - 1000r^{3} + 625r^{4},\\S_{-1}&=-4 + 40r - 120r^{2} + 40r^{3} + 220r^{4},\hskip 1cm S_{0}=230r^{4},\\S_{1}&=-4 + 40r - 120r^{2} + 40r^{3} + 220r^{4},\\S_{2}&=16 - 160r + 600r^{2} - 1000r^{3} + 625r^{4},\\\tfrac1{24}\textstyle\sum_jS_j&=1 - 10r + 40r^{2} - 80r^{3} + 80r^{4}.\end{aligned}\]
\noindent\emph{$n=6$, $r\in(0,\frac{1}{6}]$.}\[\begin{aligned}S_{0}&=2112r^{5},\hskip 1cm \tfrac1{120}\textstyle\sum_jS_j=\frac{88}{5}r^{5}.\end{aligned}\]
\noindent\emph{$n=6$, $r\in[\frac{1}{6},\frac{1}{4}]$.}\[\begin{aligned}S_{-1}&=-1 + 30r - 360r^{2} + 2160r^{3} - 6480r^{4} + 7776r^{5},\\S_{0}&=2112r^{5},\\S_{1}&=-1 + 30r - 360r^{2} + 2160r^{3} - 6480r^{4} + 7776r^{5},\\\tfrac1{120}\textstyle\sum_jS_j&=-\frac{1}{60} + \frac{1}{2}r - 6r^{2} + 36r^{3} - 108r^{4} + \frac{736}{5}r^{5}.\end{aligned}\]
\noindent\emph{$n=6$, $r\in[\frac{1}{4},\frac{1}{3}]$.}\[\begin{aligned}S_{-1}&=5 - 90r + 600r^{2} - 1680r^{3} + 1200r^{4} + 1632r^{5},\\S_{0}&=2112r^{5},\\S_{1}&=5 - 90r + 600r^{2} - 1680r^{3} + 1200r^{4} + 1632r^{5},\\\tfrac1{120}\textstyle\sum_jS_j&=\frac{1}{12} - \frac{3}{2}r + 10r^{2} - 28r^{3} + 20r^{4} + \frac{224}{5}r^{5}.\end{aligned}\]
\noindent\emph{$n=6$, $r\in[\frac{1}{3},\frac{1}{2})$.}\[\begin{aligned}S_{-2}&=-32 + 480r - 2880r^{2} + 8640r^{3} - 12960r^{4} + 7776r^{5},\\S_{-1}&=5 - 90r + 600r^{2} - 1680r^{3} + 1200r^{4} + 1632r^{5},\\S_{0}&=2112r^{5},\\S_{1}&=5 - 90r + 600r^{2} - 1680r^{3} + 1200r^{4} + 1632r^{5},\\S_{2}&=-32 + 480r - 2880r^{2} + 8640r^{3} - 12960r^{4} + 7776r^{5},\\\tfrac1{120}\textstyle\sum_jS_j&=-\frac{9}{20} + \frac{13}{2}r - 38r^{2} + 116r^{3} - 196r^{4} + \frac{872}{5}r^{5}.\end{aligned}\]

Finally, it is easy to verify that
\begin{align*}
\E[g(X_k,X_l)^2]&=4\big((2r_n)^{2m}P_2^2-2(2r_n)^{m}P_2P_3+P_3^2\big)\\
&= 4\big[(2r)^{2m}\alpha_1(r)^m-2(2r)^m\alpha_2(r)^m+\alpha_3(r)^m\big].
\end{align*}
\end{proof}

\subsection{Proof of Lemma~\ref{lem:mom}}\label{pf:lem:mom}
\mbox{}\par\begin{restated}{Lemma}{lem:mom}{ (Moment bounds)}
Let \Line{\nu_k=\#\{j\ne k:d(X_j,X_k)<2r_n\}} be the number of points within distance $2r_n$ of $X_k$. Assume $\Dbar\ge c_*$ for a constant $c_*>0$ (this holds when $r_n<\frac16$ and $nr_n^m\ge1$, and in the sparse regime of Theorem~\ref{thm:sparse} in Section~\ref{ts:thm:sparse}; $C_p$ below depends on $m,p,c_*$). For every $p\ge1$:
\begin{enumerate}[label=(\roman*),leftmargin=*]
\item $\|d_i-\Dbar\|_{L^p}\le C_p\Dbar^{1/2}$ and $\|d_i\|_{L^p}\le C_p\Dbar$;
\item $\|\nu_k\|_{L^p}\le C_p\Dbar$, and the same bound holds for \Line{\tilde\nu_k=\#\{j\ne k:d(X_j,X_k)<3r_n\};}
\item $d_k\le \nu_k$, and every neighbour $l$ of $k$ has $d_l\le \nu_k+1$; hence $\sum_{l:A_{kl}=1}d_l\le \nu_k(\nu_k+1)$;
\item the bounds in (i)--(ii) remain true, with the same constants, for the conditional laws of $d_v$ and $\nu_k$ given the positions of any fixed finite set of other points, provided the number of conditioned points is bounded by a constant (the bounds then hold with $C_p$ depending also on that number).
\end{enumerate}
\end{restated}
\begin{proof}We use the Rosenthal inequality: for $p\ge2$ there is $C_p$ such that for independent centred $\xi_1,\dots,\xi_N$ with finite $p$-th moments,
\Line{\E\big[\Big|\sum_{i=1}^N\xi_i\Big|^p\big]\le C_p\Big[\Big(\sum_{i=1}^N\E[\xi_i^2]\Big)^{p/2}+\sum_{i=1}^N\E[|\xi_i|^p]\Big].}
In particular, if $B$ is a sum of independent Bernoulli variables with parameters $\le p_0$ and $\mu=Np_0\ge1$, then $\|B-\E [B]\|_{L^p}\le C_p\sqrt\mu$ and $\|B\|_{L^p}\le C_p\mu$ (here $\|Y\|_{L^p}=(\E[|Y|^p])^{1/p}$).

\smallskip\noindent\textbf{(i), (ii): moment bounds.} By~\eqref{eq:degbin}, $d_i\sim\mathrm{Bin}(n-1,(2r_n)^m)$. Similarly
\begin{itemize}[leftmargin=*,itemsep=1pt]
\item \Line{\nu_k\sim\mathrm{Bin}\big(n-1,\min(1,(4r_n)^m)\big),} and $(4r_n)^m=2^m(2r_n)^m$;
\item \Line{\tilde\nu_k\sim\mathrm{Bin}\big(n-1,\min(1,(6r_n)^m)\big),} and $\min(1,(6r_n)^m)\le3^m(2r_n)^m$ (for $r_n\ge\frac16$ use $1\le3^m(2r_n)^m$).
\end{itemize}
Each of these is a sum of independent Bernoulli variables with parameters at most $C(2r_n)^m$ and mean at most $C\Dbar$. The Rosenthal inequality gives the claimed bounds. It is stated for means $\mu\ge1$; if $\mu<1$, use $\mu\vee1\le\max(1,1/c_*) \mu$ (valid because $\Dbar\ge c_*$), which only changes the constants.

\smallskip\noindent\textbf{(iii): comparison of degrees.} Let $l$ be a neighbour of $k$, so $d(X_k,X_l)<r_n$.
\begin{itemize}[leftmargin=*,itemsep=1pt]
\item Every neighbour $u\ne k$ of $l$ satisfies \Line{d(X_u,X_k)\le d(X_u,X_l)+d(X_l,X_k)<2r_n,} so $u$ is counted by $\nu_k$. Together with $k$ itself this gives $d_l\le\nu_k+1$.
\item Every neighbour of $k$ lies within $r_n<2r_n$ of $X_k$, so $d_k\le\nu_k$.
\end{itemize}
Summing $d_l\le\nu_k+1$ over the $d_k\le\nu_k$ neighbours of $k$ gives $\sum_{l:A_{kl}=1}d_l\le\nu_k(\nu_k+1)$.

\smallskip\noindent\textbf{(iv): conditioning.} Given the positions of finitely many points, every remaining point is still independent and uniform, and falls in the relevant region with probability at most $C(2r_n)^m$. So each count is again a sum of at most $n$ independent Bernoulli variables with such parameters, plus a bounded deterministic contribution from the conditioned points. The argument of (i)--(ii) applies.
\end{proof}

\subsection{Proof of Proposition~\ref{prop:kolm}}\label{pf:prop:kolm}
Before proving Proposition~\ref{prop:kolm}, we present a lemma.

\begin{lemma} \label{prelemma}
Assume $r_n<\frac16$. Let $g$ be as in~\eqref{eq:gdef}, \Line{\begin{gathered}G(X_k,X_l)=\E[g(X_k,X_i) g(X_i,X_l)\mid X_k,X_l],\\ \Theta(X_k,X_l)=\E[g(X_k,X_i)^2 g(X_i,X_l)\mid X_k,X_l]\end{gathered}} ($i\ne k,l$) and $\omega=\E[g(X_k,X_l)^2]$. Then 
\begin{enumerate}[label=(\alph*),leftmargin=*]
\item $g$ satisfies
\[\E[g(X_k,X_l)]=\E[H_{ij;kl}]=0, \ \ \ \ \ \mathrm{supp} g(x,0)\subset\{\|x\|_\infty\le3r_n\}\]
\[\|g\|_\infty\le4(2r_n)^{2m},\ \ \ \ \ \ \ \ \ \ \E |g(X_k,X_l)|\le4(2r_n)^{3m}.\]
\item $G$ is symmetric, and
\Line{\E[G(X_k,X_l)]=0,\qquad\|G\|_\infty=\omega\quad(\text{attained when }X_k=X_l),}
\Line{\E|G(X_k,X_l)|\le16(2r_n)^{6m},\qquad\E[G(X_k,X_l)^2]\le16(2r_n)^{6m}\omega;}
\item \Line{\E[g(X_k,X_l)^4]\le4(2r_n)^{4m}\omega,\qquad\E[g(X_k,X_l)^6]\le16(2r_n)^{8m}\omega,}
\Line{\E[\Theta(X_k,X_l)]=0,\qquad\E[\Theta(X_k,X_l)^2]\le\E[g(X_k,X_l)^4] \big(\E|g(X_k,X_l)|\big)^2\le64(2r_n)^{10m}\omega.}
\end{enumerate}
\end{lemma}
\begin{proof}Throughout, $g=g(X_k,X_l)$, and for a third index $i\ne k,l$ we write $g_{ki}=g(X_k,X_i)$, $g_{il}=g(X_i,X_l)$.

\medskip\noindent\emph{Proof of (a).} Note that
\[   P_2(X_k,X_l)=\E[A_{ki}A_{il}|X_k,X_l]\le \sqrt{\E[A_{ki}|X_k,X_l]\E[A_{il}|X_k,X_l]}=(2r_n)^m,\]
\[   P_3(X_k,X_l)=\E[A_{ki}A_{ij}A_{jl}|X_k,X_l]\le \E[A_{ki}A_{ij}|X_k,X_l]=(2r_n)^{2m}.\]
Therefore,
\[|g(X_k,X_l)|\le2 \big[(2r_n)^mP_2(X_k,X_l)+P_3(X_k,X_l)\big].\]
Then it is easy to get the (a).

\medskip\noindent\emph{Proof of (b).} By definition $G(X_k,X_l)=\E[g_{ki}g_{il}\mid X_k,X_l]$, and $G$ is symmetric because $g$ is.
\begin{itemize}[leftmargin=*]
\item \emph{$\|G\|_\infty=\omega$.} Given $X_k,X_l$, the point $X_i$ is uniform, so \Line{\E[g_{ki}^2\mid X_k,X_l]=\E[g_{il}^2\mid X_k,X_l]=\omega} by translation invariance. The conditional Cauchy--Schwarz inequality gives $|G(X_k,X_l)|\le\omega$, with equality when $X_k=X_l$ (then $g_{ki}=g_{il}$).
\item \emph{$\E|G|\le16(2r_n)^{6m}$.} \Line{|G(X_k,X_l)|\le\E[|g_{ki}||g_{il}|\mid X_k,X_l].} Given $X_i$, the variables $g_{ki}$ and $g_{il}$ are independent (functions of $X_k$, resp.\ $X_l$), and \Line{\E[|g_{ki}|\mid X_i]=\E|g|=\E[|g_{il}|\mid X_i].} Hence, by (a), \Line{\E|G(X_k,X_l)|\le(\E|g|)^2\le16(2r_n)^{6m}} by (a).
\item \emph{$\E[G]=0$.} \Line{\E[G(X_k,X_l)]=\E[\E[g_{ki}\mid X_i] \E[g_{il}\mid X_i]]=0,} because $\E[g_{ki}\mid X_i]=0$ (for every $x$, $\E[g(X_k,x)]=\E[g]=0$ by translation invariance).
\item \emph{$\E[G^2]\le16(2r_n)^{6m}\omega$.} \Line{\E[G^2]\le\|G\|_\infty \E|G|\le\omega\cdot16(2r_n)^{6m}.}
\end{itemize}

\medskip\noindent\emph{Proof of (c).} By (a), \Line{\begin{gathered}\E[g^4]\le\|g\|_\infty^2 \E[g^2]\le4(2r_n)^{4m}\omega,\\ \E[g^6]\le\|g\|_\infty^4 \E[g^2]\le16(2r_n)^{8m}\omega.\end{gathered}} For $\Theta(X_k,X_l)=\E[g_{ki}^2g_{il}\mid X_k,X_l]$: \Line{\E[\Theta]=\E[\E[g_{ki}^2\mid X_i] \E[g_{il}\mid X_i]]=\omega\cdot0=0.} Moreover, by the conditional Cauchy--Schwarz inequality applied to the factorisation $g_{ki}^2g_{il}$ $=(g_{ki}^2\sqrt{|g_{il}|})\cdot(\mathrm{sgn}(g_{il})\sqrt{|g_{il}|})$,
\[\Theta(X_k,X_l)^2\le\E\big[g_{ki}^4|g_{il}|\mid X_k,X_l\big]\cdot\E\big[|g_{il}|\mid X_k,X_l\big]=\E\big[g_{ki}^4|g_{il}|\mid X_k,X_l\big] \E|g|,\]
because \Line{\E[|g_{il}|\mid X_k,X_l]=\E[|g_{il}|\mid X_l]=\E|g|} by translation invariance. Taking expectations, and using that given $X_i$ the variables $g_{ki}$ and $g_{il}$ are independent with $\E[g_{ki}^4\mid X_i]=\E[g^4]$ and $\E[|g_{il}|\mid X_i]=\E|g|$,
\[\begin{aligned}\E[\Theta^2]&\le\E|g|\cdot\E\big[g_{ki}^4|g_{il}|\big]\\&=\E|g|\cdot\E\big[\E[g_{ki}^4\mid X_i] \E[|g_{il}|\mid X_i]\big]=\E[g^4] (\E|g|)^2,\end{aligned}\]
and by (a) and the bound on $\E[g^4]$ above,
\[\E[\Theta(X_k,X_l)^2]\le\E[g^4] (\E|g|)^2\le4(2r_n)^{4m}\omega\cdot16(2r_n)^{6m}=64(2r_n)^{10m}\omega.\]
\end{proof}

\begin{restated}{Proposition}{prop:kolm}{ (Kolmogorov bound for the pair statistic)}
If $r_n<\frac16$ and $nr_n^m\ge1$, then
\Line{\dK\big(W_n/\varsigma_n,\mathcal{N}\big)\le C ((2r_n)^{m/2}+n^{-1/2})\le C r_n^{m/2},}
where $\varsigma_n^2=\Var [W_n]$.
\end{restated}

\begin{proof} The proof is based on Theorem \ref{tl:K}. Let $I$ be uniform on $\{1,\dots,n\}$ and $X_I'$ an independent copy of $X_I$, all independent of $X=(X_1,\dots,X_n)$. Let $W=W_n/\varsigma_n$ and let $W'$ be $W$ computed with $X_I$ replaced by $X_I'$; $(W,W')$ is exchangeable. Put $M_I=\sum_{l\ne I}g_{Il}$ and $M_I'=\sum_{l\ne I}g(X_I',X_l)$, so $\Delta=W-W'=(M_I-M_I')/\varsigma_n$. Let $\E'$ denote expectation over $X_I'$ alone.

\smallskip\noindent\emph{Step A: linear regression.} \Line{\E'[M_I']=\sum_{l\ne I}\E'[g(X_I',X_l)]=0,} so
\Line{\E[\Delta\mid X]=\frac1{n\varsigma_n}\sum_IM_I=\frac1{n\varsigma_n}\cdot2\sum_{k<l}g_{kl}=\frac2nW,}
and since $W$ is a function of $X$, \[\E[\Delta\mid W]=\E[\E[\Delta\mid X]\mid W]=\frac2nW.\]
This is the hypothesis of \TS{K} with $\lambda=2/n$, $R=0$; also $\E[|W|]\le(\E [W^2])^{1/2}=1$.

{\relpenalty=0 \binoppenalty=700 \smallskip\noindent\emph{Step B: first term of \TS{K}.} Since
\Line{\E'[M_I'^2]=\sum_{l,l'\ne I}\E'[g(X_I',X_l)g(X_I',X_{l'})]=\sum_{l,l'\ne I}G(X_l,X_{l'})=:Q_I}
and $\E'[M_IM_I']=0$,
\Line{\frac1{2\lambda}\E[\Delta^2\mid X]=\mathcal R:=\frac{\mathcal S}{4\varsigma_n^2},\qquad \mathcal S=\sum_I\big(M_I^2+Q_I\big).}}
Expanding (all sums over indices $\ne I$): $M_I^2=\sum_lg_{Il}^2+\sum_{l\ne l'}g_{Il}g_{Il'}$ and \Line{Q_I=(n-1)\omega+\sum_{l\ne l'}G(X_l,X_{l'}).} Summing over $I$, $\mathcal S=\mathcal S_1+\mathcal S_2+\mathcal S_3+n(n-1)\omega$ with
\Line{\mathcal S_1=2\sum_{k<l}g_{kl}^2,\qquad \mathcal S_2=\sum_{I\ne l\ne l'}g_{Il}g_{Il'},\qquad \mathcal S_3=2(n-2)\sum_{k<l}G(X_k,X_l).}
Thus $\E [\mathcal S]=2n(n-1)\omega$ by Lemma \ref{prelemma}, and $\E [\mathcal R]=1$ (as $\varsigma_n^2=\binom n2\omega$). 

Note that $\mathcal S_2=\mathcal S_3+\mathcal S_4$ where \Line{\mathcal S_4=\sum_{I\ne l\ne l'}\big(g_{Il}g_{Il'}-G(X_l,X_{l'})\big).} Then \Line{\mathcal S-\E [\mathcal S]=(\mathcal S_1-\E [\mathcal S_1])+\mathcal S_4+2\mathcal S_3.} Variances: (i) $g^2$ has constant marginal $\omega$, so the terms of $\mathcal S_1$ are uncorrelated.  By Lemma \ref{prelemma}, we have
\Line{\Var [\mathcal S_1]\le4\binom n2\E [g^4]\le8n^2(2r_n)^{4m}\omega;} (ii) $\E[G(X_k,X_l)]=0$, so  by Lemma \ref{prelemma}, \Line{\Var [\mathcal S_3]=4(n-2)^2\binom n2\E[G(X_k,X_l)^2]\le32n^4(2r_n)^{6m}\omega;} (iii) The variance of $\mathcal S_4$ is bounded by
\[\Var [\mathcal S_4]=\sum_{I\ne l\ne l'}\E[\big(g_{Il}g_{Il'}-G(X_l,X_{l'})\big)^2]\le\binom n3\cdot36 \omega^2\le6n^3\omega^2.\] Since $\varsigma_n^4\ge n^4\omega^2/16$,
\[\begin{aligned}\Var [\mathcal R]&=\frac{\Var [\mathcal S]}{16\varsigma_n^4}\le\frac{3(8n^2(2r_n)^{4m}\omega+6n^3\omega^2+128n^4(2r_n)^{6m}\omega)}{n^4\omega^2}\\
&\le C\Big(\frac{(2r_n)^{4m}}{n^2\omega}+\frac1n+\frac{(2r_n)^{6m}}\omega\Big)\le C\Big(\frac1{n^2(2r_n)^m}+\frac1n+(2r_n)^m\Big)\le C(2r_n)^m,\end{aligned}\]
using $\omega\ge c_g(2r_n)^{5m}$ and $n(2r_n)^m\ge c$. By the Cauchy-Schwarz inequality,
\Line{\E[|1-\frac1{2\lambda}\E[\Delta^2|W]|]\le\E[|1-\mathcal R|]\le\sqrt{\Var [\mathcal R]}\le C\sqrt{(2r_n)^m}.}

\smallskip\noindent\emph{Step C: third-moment term.} Note that
\[\E[\Delta^3|W]=\E[\E[\Delta^3|X]\mid W].\]
Then \[\E[|\E[\Delta^3|W]|]\le\E[|\E[\Delta^3|X]|]. \]
Expanding the cube and using $\E'[M_I']=0$,
\[\begin{aligned}\E[\Delta^3\mid X]&=\frac{\mathcal K}{n\varsigma_n^3},\qquad \mathcal K=\sum_I\big(M_I^3+3M_IQ_I-\Gamma_I\big),\\
\Gamma_I&=\E'[M_I'^3]=\sum_{l_1,l_2,l_3\ne I}\Gamma(X_{l_1},X_{l_2},X_{l_3}),\end{aligned}\]
where \Line{\Gamma(X_{l_1},X_{l_2},X_{l_3})=\E[g(X_{l_1},X_I') g(X_{l_2},X_I') g(X_{l_3},X_I')\mid X_{l_1},X_{l_2},X_{l_3}].} Split according to which of the summation indices coincide (all indices $\ne I$; $\Theta(X_l,X_{l'})$ is defined in Lemma~\ref{prelemma}), so 
\begin{align*}
\mathcal K&=\mathcal K_1+\mathcal K_2+\mathcal K_3+\mathcal K_4+\mathcal K_5+\mathcal K_6,\\
\mathcal K_1&=\sum_{I\ne l}\big(g_{Il}^3-\E[g^3]\big),\quad \mathcal K_2=3(n-1)\omega\sum_IM_I,\quad \mathcal K_3=3\sum_{I\ne l\ne l'}\big(g_{Il}^2g_{Il'}-\Theta(X_l,X_{l'})\big),\\
\mathcal K_4&=6\sum_{I\ne l\ne l'}g_{Il}G(X_l,X_{l'}),\quad \mathcal K_5=3\sum_{I, l'',l,l'\text{ distinct}}g_{Il''}G(X_l,X_{l'}),\\\mathcal K_6&=\sum_{I,l_1,l_2,l_3\text{ distinct}}\big(g_{Il_1}g_{Il_2}g_{Il_3}-\Gamma(X_{l_1},X_{l_2},X_{l_3})\big).
\end{align*}
We bound the variance of each piece by Lemma~\ref{prelemma} and Hoeffding projections (\TS{D}).
\begin{itemize}[leftmargin=*]
\item $\mathcal K_1$: pair terms with constant marginal $\E[g^3]$, so uncorrelated: \Line{\Var [\mathcal K_1]\le4\binom n2\E [g^6]\le32n^2(2r_n)^{8m}\omega.}
\item $\mathcal K_2=6(n-1)\omega\sum_{k<l}g_{kl}$: \Line{\Var [\mathcal K_2]\le36n^2\omega^2\binom n2\omega\le18n^4\omega^3.}
\item $\mathcal K_3$: with $f=g_{Il}^2g_{Il'}-\Theta(X_l,X_{l'})$ one checks $\E[f|X_I,X_{l'}]=\omega g_{Il'}$ (as $\E[g_{Il}^2\mid X_I]=\omega$) while the conditional expectations of $f$ given $(X_I,X_l)$, $(X_l,X_{l'})$ or a single coordinate vanish (as $\E[g_{Il'}|X_I]=0$ and $\E\Theta=0$). Hence $f-\omega g_{Il'}$ is canonical on $\{I,l,l'\}$, and \Line{3\sum_{I\ne l\ne l'}\omega g_{Il'}=6\omega(n-2)\sum_{k<l}g_{kl}.} So \Line{\mathcal K_3=6\omega(n-2)\sum_{k<l}g_{kl}+3\sum_{I\ne l\ne l'}\big(g_{Il}^2g_{Il'}-\Theta(X_l,X_{l'})-\omega g_{Il'}\big);} the first part has variance $\le Cn^4\omega^3$, the canonical part $\le Cn^3\E [f^2]$, and \Line{\E[f^2]\le2\E[g^4] \omega+2\E[\Theta^2]\le C(2r_n)^{4m}\omega^2.} Then
\[\Var [\mathcal K_3]\le Cn^4\omega^3 +Cn^3(2r_n)^{4m}\omega^2\]
\item $\mathcal K_4$: with $f=g_{Il}G(X_l,X_{l'})$, $\E[f|X_I,X_{l'}]=G_3(X_I,X_{l'})$ and the other projections vanish; by conditional Cauchy--Schwarz, \Line{\E[G_3^2]\le(\E|g|)^2\|G\|_\infty\E|G|\le256(2r_n)^{12m}2\omega.} So
\Line{\Var [\mathcal K_4]\le Cn^4 \E[G_3^2]+Cn^3\E [f^2]\le Cn^4(2r_n)^{12m}2\omega+Cn^3(2r_n)^{6m}\omega^2,}
using $\E [f^2]=\omega \E[G(X_l,X_{l'})^2]$.
\item $\mathcal K_5$: $f=g_{Il''}G(X_l,X_{l'})$ is canonical on the $4$-set: \Line{\Var [\mathcal K_5]\le Cn^4\omega \E[G(X_l,X_{l'})^2]\le Cn^4(2r_n)^{6m}\omega^2.}
\item $\mathcal K_6$: canonical on the $4$-set (subtracting $\Gamma$ removes the projection on $\{l_1,l_2,l_3\}$, projections containing $I$ vanish since $\E[g]=0$). There are at most $n^4$ four-sets with $24$ terms each, and given $X_I$ the three factors are independent, so $\E[f^2]\le\E[(g_{Il_1}g_{Il_2}g_{Il_3})^2]=\omega^3$; hence \[\Var[\mathcal K_6]\le24^2n^4\omega^3.\]
\end{itemize}
Dividing by $n^4\omega^3$ and using $\omega\ge c_g(2r_n)^{5m}$, $n(2r_n)^m\ge c$, each ratio is $\le C$, therefore,
\[\Var [\mathcal K]\le Cn^4\omega^3.\]  
With $\varsigma_n^3\ge n^3\omega^{3/2}/8$,
\[\frac{\E[|\E[\Delta^3|W]|]}\lambda\le\frac n2\cdot\frac{\E[|\mathcal K|]}{n\varsigma_n^3}\le\frac{C}{n}.\]
\emph{Conclusion.} \TS{K} gives \Line{\dK(W,\mathcal{N})\le C\sqrt{(2r_n)^m}+Cn^{-1/2},} and $n^{-1/2}\le C\sqrt{(2r_n)^m}$ because $n(2r_n)^m\ge c$.
\end{proof}

\subsection{Proof of Lemma~\ref{lem:pert}}\label{pf:lem:pert}
\mbox{}\par\begin{restated}{Lemma}{lem:pert}{ (Perturbation)}
Let $X,Y$ be random variables and $\varepsilon>0$. Then
\Line{\dK(X+Y,\mathcal{N})\le\dK(X,\mathcal{N})+\frac{\varepsilon}{\sqrt{2\pi}}+\P(|Y|>\varepsilon).}
In particular \Line{\dK(X+Y,\mathcal{N})\le\dK(X,\mathcal{N})+C(\E [Y^2])^{1/3}.} Moreover for $c\in[\frac12,2]$: \Line{\dK(cX,\mathcal{N})\le\dK(X,\mathcal{N})+C|c-1|,} and for every $x_0\in\R$: \Line{\dK(X+x_0,\mathcal{N})\le\dK(X,\mathcal{N})+|x_0|/\sqrt{2\pi}.}
\end{restated}
\begin{proof}On $\{|Y|\le\varepsilon\}$ we have $X+Y\le x\Rightarrow X\le x+\varepsilon$. Hence
\begin{align*}
\P(X+Y\le x)&\le\P(X\le x+\varepsilon)+\P(|Y|>\varepsilon)\\
&\le\Phi(x+\varepsilon)+\dK(X,\mathcal{N})+\P(|Y|>\varepsilon)\\
&\le\Phi(x)+\frac{\varepsilon}{\sqrt{2\pi}}+\dK(X,\mathcal{N})+\P(|Y|>\varepsilon),
\end{align*}
using $\sup_x\Phi'(x)=1/\sqrt{2\pi}$. Symmetrically, since \Line{\{X\le x-\varepsilon\}\subseteq\{X+Y\le x\}\cup\{|Y|>\varepsilon\},}
\begin{align*}
\P(X+Y\le x)&\ge\P(X\le x-\varepsilon)-\P(|Y|>\varepsilon)\\
&\ge\Phi(x)-\frac\varepsilon{\sqrt{2\pi}}-\dK(X,\mathcal{N})-\P(|Y|>\varepsilon).
\end{align*}
For the second claim take $\varepsilon=(\E [Y^2])^{1/3}$ and use Chebyshev's inequality \Line{\P(|Y|>\varepsilon)\le\E [Y^2]/\varepsilon^2=(\E [Y^2])^{1/3}.} For the scaling, $\P(cX\le x)=\P(X\le x/c)$ and \Line{\sup_x|\Phi(x/c)-\Phi(x)|\le C|c-1|} because $\frac{d}{dc}\Phi(x/c)=-\frac xc\phi(x/c)$ is bounded uniformly in $x$ and $c\in[\frac12,2]$ ($\phi$ is the normal density). The shift claim is the mean value theorem for $\Phi$.
\end{proof}

\subsection{Proof of Lemma~\ref{lem:degmom}}\label{pf:lem:degmom}
\mbox{}\par\begin{restated}{Lemma}{lem:degmom}{ (Moments of degenerate sums)}
Let $X_1,\dots,X_n$ be independent, $t\ge1$ and $s\ge1$ integers with $n\ge 2st$, and for each $t$-element set $T\subseteq\{1,\dots,n\}$ let $\chi_T(X_T)$ be a random variable with $|\chi_T|\le B$ and $\E[\chi_T\mid X_{T\setminus\{i\}}]=0$ for every $i\in T$ (complete degeneracy). Then, with $q=2s$,
\Line{\Big\|\sum_{|T|=t}\chi_T(X_T)\Big\|_{L^q}\le C_t q^{t/2} n^{t/2}B .}
\end{restated}
\begin{proof}\smallskip\noindent\emph{Step 1 (expansion).} Write
\Line{\E\Big[\Big(\sum_T\chi_T\Big)^{2s}\Big]=\sum_{(T_1,\dots,T_{2s})}\E\big[\chi_{T_1}\cdots \chi_{T_{2s}}\big],}
the sum being over all $2s$-tuples of $t$-subsets of $\{1,\dots,n\}$.

\smallskip\noindent\emph{Step 2 (most terms vanish).} Suppose an index $i$ belongs to exactly one set of the tuple, say $T_a$. Condition on all variables $X_j$, $j\ne i$: the factors $\chi_{T_b}$, $b\ne a$, do not depend on $X_i$ and come out of the conditional expectation, while $\E[\chi_{T_a}\mid X_{T_a\setminus\{i\}}]=0$. So the term vanishes. Hence only tuples in which every index of $V=T_1\cup\dots\cup T_{2s}$ occurs in at least two of the sets survive.

\smallskip\noindent\emph{Step 3 (few distinct indices).} 
Write each set as a list of its $t$ elements, $T_a=\{i_{a,1},\dots,i_{a,t}\}$ (distinct elements). A \emph{slot} is a position $(a,r)$ with $1\le a\le2s$ and $1\le r\le t$; the slot $(a,r)$ \emph{holds} the index $i_{a,r}$. Thus, the number of slots is $J=2s\cdot t=2st$ (independent of the tuple); every slot holds exactly one index, so the slots are partitioned according to the index they hold. Let $u=|V|$. For $i\in V$, let $m_i\ge1$ be the number of slots holding $i$. Then
\begin{equation*} 
\sum_{i\in V}m_i=J,
\end{equation*}
because each of the $J$ slots holds exactly one index of $V$. Hence $m_i\ge2$. Then
\[
J=\sum_{i\in V}m_i\ \ge\ \sum_{i\in V}2=2u.
\]
Hence, $u\le J/2=st$.

\smallskip\noindent\emph{Step 4 (counting).} We count the tuples with $|V|=u$ in two stages. First choose $V$, a $u$-element subset of $\{1,\dots,n\}$: there are
\Line{\binom nu\le\frac{n^u}{u!}\le\Big(\frac{en}{u}\Big)^u}
choices, because $u!\ge(u/e)^u$. Then choose each $T_a$ as a $t$-subset of $V$: there are $\binom ut\le u^t$ choices for each of the $2s$ sets, hence at most $u^{2st}=u^J$ for the whole tuple. A tuple is determined by $V$ and its sets, so
\Line{\#\{\text{tuples with }|V|=u\}\le\Big(\frac{en}u\Big)^uu^J=e^u n^u u^{J-u}.}

\smallskip\noindent\emph{Step 5 (summation).} Each term is at most $B^{2s}$ in absolute value. Since $u\le J/2\le n$ (here $n\ge J$ is used), \Line{n^uu^{J-u}=n^{J/2}u^{J/2}(u/n)^{J/2-u}\le n^{J/2}(J/2)^{J/2},} and $e^u\le e^{J/2}$. There are at most $J/2$ values of $u$, so
\Line{\E\Big[\Big(\sum_T\chi_T\Big)^{2s}\Big]\le\frac J2 e^{J/2}\Big(\frac J2\Big)^{J/2}n^{J/2}B^{2s}=st e^{st}(st)^{st}n^{st}B^{2s}.}
Taking the $2s$-th root gives $(st)^{1/(2s)}e^{t/2}(st)^{t/2}n^{t/2}B$. Here \Line{(st)^{1/(2s)}=t^{1/(2s)}s^{1/(2s)}\le t^{1/2}e^{1/(2e)},} because $s^{1/s}\le e^{1/e}$ for $s\ge1$, and $(st)^{t/2}=(t/2)^{t/2}q^{t/2}$ since $s=q/2$. Hence \Line{\big\|\sum_T\chi_T\big\|_{L^q}=\big(\E[(\sum_T\chi_T)^{2s}]\big)^{1/q}} is at most $C_tq^{t/2}n^{t/2}B$ with $C_t=t^{1/2}e^{1/(2e)}e^{t/2}(t/2)^{t/2}$.
\end{proof}

\subsection{Proof of Lemma~\ref{lem:remmom}}\label{pf:lem:remmom}
\mbox{}\par\begin{restated}{Lemma}{lem:remmom}{}
For constant $r\in(0,\frac12)$ there are $C,n_0$ such that for $n\ge n_0$ and every even integer $q$ with $2\le q\le n^{1/3}$, \[\|\mathrm{Rem}\|_{L^q}\le C q^{3/2} n^{-1/2}.\]
\end{restated}
\begin{proof}We use Lemma~\ref{lem:degmom} (Section~\ref{ts:lem:degmom}) and \TS{L} to prove Lemma~\ref{lem:remmom} (Section~\ref{ts:lem:remmom}). Recall $n\Dbar^2\ge cn^3$ and \Line{\mathrm{Rem}=\mathrm{Rem}^{(1)}+\mathrm{Rem}^{(2)}+\mathrm{Rem}^{(3)}} from~\eqref{eq:split}.

\emph{$\mathrm{Rem}^{(1)}$.} By \TS{D}, \[\Upsilon_n^{(3)}=\sum_{|T|=3}\pi_T\Upsilon_n,\ \ \ \ \ \Upsilon_n^{(4)}=\sum_{|U|=4}\pi_U\Upsilon_n,\]
and $\pi_S\Upsilon_n$ is a function of $X_S$ with $\E[\pi_S\Upsilon_n|X_{S\setminus\{i\}}]=0$ for $i\in S$ (\TS{D}(c)), i.e.\ completely degenerate. Each summand $H_{ij;kl}$ of $\Upsilon_n$ satisfies $|H_{ij;kl}|\le1$, and $\pi_SH_{ij;kl}$ is an alternating sum of $2^{|S|}$ conditional expectations of $H_{ij;kl}$, so $|\pi_SH_{ij;kl}|\le2^{|S|}$; a given triple lies in the index sets of $12(n-3)$ summands and a given $4$-set carries $12$. Thus $|\pi_T\Upsilon_n|\le96n$ and $|\pi_U\Upsilon_n|\le192$. Lemma~\ref{lem:degmom} (Section~\ref{ts:lem:degmom}) with $t=3$, $B=96n$ and with $t=4$, $B=192$ gives \[\|\Upsilon_n^{(3)}\|_{L^q}\le Cq^{3/2}n^{5/2},\ \ \ \ \ \ \ \ \|\Upsilon_n^{(4)}\|_{L^q}\le Cq^2n^2\] (note $q\le n^{1/3}$ implies $n\ge4q$ for large $n$). Dividing by $n\Dbar^2\ge cn^3$ and using $q^2n^{-1}\le q^{3/2}n^{-1/2}$: \[\|\mathrm{Rem}^{(1)}\|_{L^q}\le Cq^{3/2}n^{-1/2}.\]

\emph{$\mathrm{Rem}^{(2)}$.} By Lemma~\ref{lem:S} (Section~\ref{ts:lem:S}), $\mathcal S=\sum_{|T|=3}s(X_T)$ with, for $T=\{i,j,k\}$, \Line{s(X_T)=2(A_{ij}A_{ik}+A_{ji}A_{jk}+A_{ki}A_{kj})-6A_{ij}A_{jk}A_{ki}} (each vertex $i$ contributes $d_i(d_i-1)=2\sum_{\{j,k\}}A_{ij}A_{ik}$ and each triangle is counted once), so $|s|\le12$. By translation invariance $\E[s(X_T)|X_i]$ does not depend on $X_i$, so $\E[\mathcal  S|X_i]=\E \mathcal  S$ and $\mathcal  S$ has no one-point component. The two-point components are \Line{\pi_{\{k,l\}}\mathcal  S=\sum_{T\supset\{k,l\}}\pi_{\{k,l\}}s(X_T),} a sum of $n-2$ terms bounded by $4\cdot12$, hence $\le48n$, and the three-point components are $\pi_Ts$, bounded by $96$. Lemma~\ref{lem:degmom} (Section~\ref{ts:lem:degmom}) ($t=2$, $B=48n$; $t=3$, $B=96$) gives \Line{\|\mathcal  S-\E \mathcal S\|_{L^q}\le C(q n^2+q^{3/2}n^{3/2}),} so \[\|\mathrm{Rem}^{(2)}\|_{L^q}\le Cq^{3/2}n^{-1}.\]

\emph{$\mathrm{Rem}^{(3)}=R_n-\E R_n$.} By \TS{L}, namely~\eqref{eq:BBLM}, it suffices to show 
\[\|n(R_n-R_n^{(k)})\|_{L^q}\le Cq\]
for all $k$, where $R_n^{(k)}$ is $R_n$ with $X_k$ replaced by an independent copy $X_k'$: then \Line{\|R_n-\E R_n\|_{L^q}\le2\sqrt{\varkappa q} (n C^2q^2/n^2)^{1/2}=Cq^{3/2}n^{-1/2}.} We use the notation of Step~5 of Section~\ref{sec:proofs} ($\nu$, $\mathcal D_k$, $\mathcal E_k$, $N(k)$, $N'(k)$), now with constant radius, so that $\Dbar\asymp n$, $\nu\le\tilde\nu\le2n$ and $|\mathcal D_k|\le n$ deterministically; take $n_0$ so large that $\Dbar>16$. Each degree is binomial, $d_v\sim\mathrm{Bin}(n-1,(2r)^m)$, and Bernstein's inequality integrated against the tail gives \Line{\|d_v-\Dbar\|_{L^p}\le C(\sqrt{pn}+p)\le C\sqrt{pn}} for $1\le p\le n$.

\emph{On $\mathcal E_k$.} The bound derived in Step~5 reads
\[\begin{aligned}|n(R_n-R_n^{(k)})|\1_{\mathcal E_k}&\le\frac C{\Dbar^3}\Big[\sum_{l\in N(k)}Y_{kl}^3+\sum_{l\in N'(k)}(Y_{kl}')^3+(\tilde\nu+2)\sum_{v\in\mathcal D_k}\big((d_v-\Dbar)^2+1\big)\Big],\\
Y_{kl}&=|d_k-\Dbar|+|d_l-\Dbar|,\qquad Y_{kl}'=|d_k'-\Dbar|+|d_l'-\Dbar| .\end{aligned}\]
Note that
\Line{Y_{kl}^3\le 4(|d_k-\Dbar|^3+|d_l-\Dbar|^3)} and likewise for $Y_{kl}'$, the number of terms is at most $2n$, and the second sum has at most $n$ terms, each $\le(d_v-\Dbar)^2+1$. It is clear that
\Line{\begin{gathered}\big\| |d_v-\Dbar|^3\big\|_{L^q}=\|d_v-\Dbar\|_{L^{3q}}^3\le C(qn)^{3/2},\\ \big\|(d_v-\Dbar)^2+1\big\|_{L^q}\le C(qn+1),\end{gathered}} using $q\le n$. By Minkowski's inequality and Hölder's inequality 
\[\begin{aligned}\Big\|\sum_{l}Y_{kl}^3\Big\|_{L^q}&\le Cn (qn)^{3/2},\\\Big\|(\tilde\nu+2)\sum_{v\in\mathcal D_k}\big((d_v-\Dbar)^2+1\big)\Big\|_{L^q}&\le Cn\cdot n (qn+1)\le Cqn^3 .\end{aligned}\]
Dividing by $\Dbar^3\ge cn^3$ gives \Line{\big\|n(R_n-R_n^{(k)})\1_{\mathcal E_k}\big\|_{L^q}\le C(q^{3/2}n^{-1/2}+q)\le Cq.}

\emph{On $\mathcal E_k^c$.} By~\eqref{eq:crude}, \Line{\begin{gathered}|n(R_n-R_n^{(k)})|\le C(\tilde\nu+2)^3\le Cn^3,\\ \P(\mathcal E_k^c)\le C\Dbar e^{-\Dbar/60}\le Cne^{-cn}\end{gathered}} (Step~5), so \Line{\big\|n(R_n-R_n^{(k)})\1_{\mathcal E_k^c}\big\|_{L^q}\le Cn^3\big(Cne^{-cn}\big)^{1/q}\le\exp\big(\log(Cn^3)+(\log(Cn)-cn)/q\big)\le1,} because $q\le n^{1/3}$.

Adding the two cases ($\mathcal E_k$ and $\mathcal E_k^c$) with Minkowski's inequality gives \Line{\|n(R_n-R_n^{(k)})\|_{L^q}\le Cq+1\le C'q,} which is what was needed for $\mathrm{Rem}^{(3)}$. Combining the three pieces with Minkowski's inequality in $L^q$, \Line{\|\mathrm{Rem}\|_{L^q}\le\sum_{i=1}^3\|\mathrm{Rem}^{(i)}\|_{L^q}\le C q^{3/2}n^{-1/2}} (note $q^{3/2}n^{-1}\le q^{3/2}n^{-1/2}$), which is the claim of the lemma.
\end{proof}

\subsection{Proof of Lemma~\ref{lem:comp}}\label{pf:lem:comp}
\mbox{}\par\begin{restated}{Lemma}{lem:comp}{ (Components)}
For a connected graph $\mathcal C$ let $f(\mathcal C)=\sum_{\{i,j\}\in E(\mathcal C)}\frac{(d_i-d_j)^2}{d_id_j}$ and $N_\vee(\mathcal C)$ the number of induced two-paths (triples $\{i,j,k\}$ with exactly the two edges $ij,jk$) in $\mathcal C$. Then
\begin{enumerate}[label=(\alph*),leftmargin=*]
\item $f(\mathcal C)=N_\vee(\mathcal C)$ if $|\mathcal C|\le3$;
\item $0\le f(\mathcal C),N_\vee(\mathcal C)\le|\mathcal C|^3$ for every $\mathcal C$.
\end{enumerate}
Consequently $n\mathcal{F}_n=\sum_{\mathcal C}f(\mathcal C)$ and $N_\vee=\sum_{\mathcal C}N_\vee(\mathcal C)$ (sums over connected components of the random graph), and
\Line{\Xi:=n\mathcal{F}_n-N_\vee=\sum_{|\mathcal C|\ge4}\big(f(\mathcal C)-N_\vee(\mathcal C)\big),}
which satisfies
\Line{|\Xi|\le\sum_{v=1}^n|\mathcal C(v)|^2 \1\{|\mathcal C(v)|\ge4\},\qquad\mathcal C(v)=\text{component of }v .}
In particular, let $\Omega_4$ be the event that the graph has a connected component with at least four vertices; then $n\mathcal{F}_n=N_\vee$ on $\Omega_4^c$.
\end{restated}
\begin{proof}We use Proposition~\ref{prop:id} (Section~\ref{ts:prop:id}) to prove Lemma~\ref{lem:comp} (Section~\ref{ts:lem:comp}).
Throughout, $|\mathcal C|$ is the number of vertices of $\mathcal C$; in a component the degrees $d_i$ are the degrees in the whole graph, because all neighbours of a vertex lie in its component.

\medskip\noindent\textbf{(a) Components with at most three vertices.} A connected graph on at most three vertices is a single vertex, an edge, a path on three vertices, or a triangle. We check that $f(\mathcal C)=N_\vee(\mathcal C)$ in each case. A single vertex has no edge, so $f=0$, and it contains no triple, so $N_\vee=0$. An edge has degrees $(1,1)$, so its only term is $\frac{(1-1)^2}{1\cdot1}=0$, and it has only two vertices, so $N_\vee=0$. A path on three vertices has degrees $(1,2,1)$ and two edges, each contributing $\frac{(1-2)^2}{1\cdot2}=\frac12$, so $f=1$; its only triple is the path itself, which has exactly two edges and is therefore an induced two-path, so $N_\vee=1$. A triangle has degrees $(2,2,2)$, so every edge term vanishes and $f=0$; its only triple has three edges, not two, so it is not an induced two-path and $N_\vee=0$. Hence $f(\mathcal C)=N_\vee(\mathcal C)$ in all four cases, which is~(a).

\medskip\noindent\textbf{(b) The bounds $0\le f,N_\vee\le|\mathcal C|^3$.} Both quantities are non-negative (each edge term is a square divided by a positive number; $N_\vee$ is a count).
\begin{itemize}[leftmargin=*,itemsep=2pt]
\item \emph{$f(\mathcal C)\le|\mathcal C|^3$.} For an edge $\{i,j\}$ both degrees are at least $1$ and at most $|\mathcal C|-1$, so, for $x=d_i$, $y=d_j$,
\[
\frac{(x-y)^2}{xy}=\frac xy+\frac yx-2\le x+y\le2|\mathcal C| ,
\]
because $\frac xy\le x$ and $\frac yx\le y$ when $x,y\ge1$. The component has at most $\binom{|\mathcal C|}2\le\frac{|\mathcal C|^2}2$ edges, hence $f(\mathcal C)\le\frac{|\mathcal C|^2}2\cdot2|\mathcal C|=|\mathcal C|^3$.
\item \emph{$N_\vee(\mathcal C)\le|\mathcal C|^3/2$.} An induced two-path is determined by its centre ($|\mathcal C|$ choices) and the unordered pair of its two ends (at most $\binom{|\mathcal C|}2\le\frac{|\mathcal C|^2}2$ choices), so $N_\vee(\mathcal C)\le|\mathcal C|\cdot\frac{|\mathcal C|^2}2=\frac{|\mathcal C|^3}2$.
\end{itemize}

\medskip\noindent\textbf{Splitting over components.} By Proposition~\ref{prop:id} (Section~\ref{ts:prop:id}), $n\mathcal F_n=\sum_{\{i,j\}\in E}\frac{(d_i-d_j)^2}{d_id_j}$, a sum over the edges of the graph. Every edge lies in exactly one component, and the degrees of its endpoints are the same in the component, so $n\mathcal F_n=\sum_{\mathcal C}f(\mathcal C)$. Likewise, a triple that induces a two-path is connected, so it lies inside one component, and being induced does not depend on whether we look at the component or at the whole graph; hence $N_\vee=\sum_{\mathcal C}N_\vee(\mathcal C)$. Subtracting, and using~(a) for the components with at most three vertices (which contribute $0$),
\[
\Xi=n\mathcal F_n-N_\vee=\sum_{|\mathcal C|\ge4}\big(f(\mathcal C)-N_\vee(\mathcal C)\big).
\]

\medskip\noindent\textbf{The bound for $\Xi$.} Since $f(\mathcal C)$ and $N_\vee(\mathcal C)$ are both non-negative, \Line{|f(\mathcal C)-N_\vee(\mathcal C)|\le\max\big(f(\mathcal C),N_\vee(\mathcal C)\big)\le|\mathcal C|^3} by~(b). Therefore $|\Xi|\le\sum_{|\mathcal C|\ge4}|\mathcal C|^3$. Finally, a component with $k$ vertices has $k$ vertices $v$, each with $|\mathcal C(v)|^2=k^2$, so
\[
\sum_{|\mathcal C|\ge4}|\mathcal C|^3=\sum_{v=1}^n|\mathcal C(v)|^2 \1\{|\mathcal C(v)|\ge4\},
\]
which is the stated bound.
\end{proof}

\subsection{Proof of Lemma~\ref{lem:csize}}\label{pf:lem:csize}
\mbox{}\par\begin{restated}{Lemma}{lem:csize}{ (Size of a component)}
For every $k\ge2$, \Line{\P(|\mathcal C(1)|\ge k)\le e (en(2r_n)^m)^{k-1}.} Consequently, for $p\ge0$ and $n(2r_n)^m\le\frac1{2e}$, where $\Omega_4$ is the event of Lemma~\ref{lem:comp} (Section~\ref{ts:lem:comp}),
\Line{\E\big[|\mathcal C(1)|^p\1\{|\mathcal C(1)|\ge4\}\big]\le C_p\big(n(2r_n)^m\big)^{3},\qquad\P(\Omega_4)\le\frac23 n^{4}(2r_n)^{3m} .}
\end{restated}
\begin{proof}We use Lemma~\ref{lem:tree} (Section~\ref{ts:lem:tree}) and Cayley's formula to prove Lemma~\ref{lem:csize} (Section~\ref{ts:lem:csize}). Throughout, $\mathcal C(1)$ denotes the connected component of vertex $1$ in the graph $\mathcal G_n(r_n)$.

\medskip\noindent\textbf{Step 1: a component with $\ge k$ vertices contains a tree on $k$ vertices through $1$.} Suppose $|\mathcal C(1)|\ge k$. There are a $k$-set $V\ni1$ and a tree $\mathcal T$ on $V$ (a labelled tree) such that $A_{uw}=1$ for every edge $\{u,w\}$ of $\mathcal T$. 

\medskip\noindent\textbf{Step 2: union bound.} The number of sets $V\subset\{1,\dots,n\}$ with $|V|=k$ and $1\in V$ is $\binom{n-1}{k-1}\le\frac{n^{k-1}}{(k-1)!}$; for a given $V$ the number of labelled trees on $V$ is $k^{k-2}$ (Cayley's formula); and, by Lemma~\ref{lem:tree} (Section~\ref{ts:lem:tree}), each fixed tree has all its edges present with probability exactly $(2r_n)^{m(k-1)}$. Hence, by Step~1 and the union bound, the following is a bound for the \emph{tail} event $\{|\mathcal C(1)|\ge k\}$:
\Line{\P\big(|\mathcal C(1)|\ge k\big)\le\frac{n^{k-1}}{(k-1)!} k^{k-2} (2r_n)^{m(k-1)}=\frac{k^{k-2}}{(k-1)!} \big(n(2r_n)^m\big)^{k-1}=\frac{k^{k-1}}{k!} \big(n(2r_n)^m\big)^{k-1}.}

\medskip\noindent\textbf{Step 3: the elementary estimate.} Since $e^k=\sum_{j\ge0}k^j/j!\ge k^k/k!$ we have $k!\ge(k/e)^k$, so $\frac{k^{k-1}}{k!}\le\frac{k^{k-1}e^k}{k^k}=\frac{e^k}{k}$ and
\Line{\P\big(|\mathcal C(1)|\ge k\big)\le\frac{e^k}{k}\big(n(2r_n)^m\big)^{k-1}\le e (en(2r_n)^m)^{k-1}\qquad(k\ge2),}
which is the first claim. (No restriction on $n(2r_n)^m$ was used so far; the bound is only informative for $en(2r_n)^m<1$.)

\medskip\noindent\textbf{Step 4: the moment bound.} Put $x=en(2r_n)^m$ and assume $n(2r_n)^m\le\frac1{2e}$, i.e.\ $x\le\frac12$. For a non-negative integer variable $Y=|\mathcal C(1)|$ and $p\ge0$, using $\P(Y=k)\le\P(Y\ge k)$ and Step~3,
\Line{\E\big[Y^p\1\{Y\ge4\}\big]=\sum_{k\ge4}k^p \P(Y=k)\le\sum_{k\ge4}k^p \P(Y\ge k)\le e\sum_{k\ge4}k^px^{k-1}=e x^3\sum_{j\ge0}(j+4)^px^j,}
where we substituted $k=j+4$ and pulled out $x^3=x^{(j+4)-1}x^{-j}$. Because $x\le\frac12$, the last series is at most $S_p:=\sum_{j\ge0}(j+4)^p2^{-j}<\infty$, so
\[
\E\big[|\mathcal C(1)|^p\1\{|\mathcal C(1)|\ge4\}\big]\le e S_p x^3=e^4S_p \big(n(2r_n)^m\big)^{3}=:C_p\big(n(2r_n)^m\big)^{3} ,
\]
and $C_p$ depends only on $p$.

\medskip\noindent\textbf{Step 5: the probability of $\Omega_4$.} Recall that $\Omega_4$ is the event that \emph{some} component of the graph has at least four vertices. On $\Omega_4$, there exists  a $4$-set $V$ carrying a tree $\mathcal T$ on $V$. So $\Omega_4$ is contained in the union, over all $4$-subsets $V\subset\{1,\dots,n\}$ and all labelled trees $\mathcal T$ on $V$, of the events \Line{\{A_{uw}=1\text{ for the three edges of }\mathcal T\}.} There are $\binom n4\le\frac{n^4}{24}$ four-sets, $4^{4-2}=16$ labelled trees on each, and each tree is present with probability $(2r_n)^{3m}$ (Lemma~\ref{lem:tree} in Section~\ref{ts:lem:tree}). Hence
\Line{\P(\Omega_4)\le\frac{n^4}{24}\cdot16\cdot (2r_n)^{3m}=\frac23 n^4(2r_n)^{3m}=\frac23 n^{4}(2r_n)^{3m},}
which is the second claim of the lemma. This step needs no assumption on $n(2r_n)^m$.
\end{proof}

\subsection{Proof of Lemma~\ref{lem:meanNvee}}\label{pf:lem:meanNvee}
\mbox{}\par\begin{restated}{Lemma}{lem:meanNvee}{ (Exact mean)}
$p_\alpha=\kappa_m (2r_n)^{2m}$ for every $\alpha$, hence \Line{\lambda_n=\E [N_\vee]=\frac{\kappa_m}2n(n-1)(n-2)(2r_n)^{2m}.}
\end{restated}
\begin{proof}We use Lemma~\ref{lem:kappa} (Section~\ref{ts:lem:kappa}) to prove Lemma~\ref{lem:meanNvee} (Section~\ref{ts:lem:meanNvee}). Since 
\[I_\alpha=A_{ji}A_{jk}-A_{ji}A_{jk}A_{ik},\]
then
\[ 
p_\alpha=\E[A_{ji}A_{jk}]-\E[A_{ji}A_{jk}A_{ki}]=(2r_n)^{2m}-(2r_n)^{2m}(1-\kappa_m)=(2r_n)^{2m}\kappa_m.
\]
(for $r_n<\frac16$, $\kappa(r_n)=\kappa_m$).
Summing over the $|\mathcal I|=n(n-1)(n-2)/2$ indices gives $\lambda_n$.
\end{proof}

\subsection{Proof of Lemma~\ref{lem:b}}\label{pf:lem:b}
\mbox{}\par\begin{restated}{Lemma}{lem:b}{ (Neighbourhoods of dependence)}
Let $B_\alpha=\{\beta\in \mathcal I:V_\beta\cap V_\alpha\ne\emptyset\}$ (which contains $\alpha$), and let $b_1,b_2,b_3$ be the Chen--Stein quantities of \TS{I} for this choice of $B_\alpha$, namely
\[\begin{aligned}&b_1=\sum_{\alpha\in \mathcal I}\sum_{\beta\in B_\alpha}p_\alpha p_\beta,\qquad b_2=\sum_{\alpha\in \mathcal I}\sum_{\beta\in B_\alpha\setminus\{\alpha\}}\E[I_\alpha I_\beta],\\
&b_3=\sum_{\alpha\in \mathcal I}\E\Big[\Big|\E\big[I_\alpha-p_\alpha \big| \sigma(I_\beta:\beta\notin B_\alpha)\big]\Big|\Big] .\end{aligned}\]
Then $I_\alpha$ is independent of $(I_\beta)_{\beta\notin B_\alpha}$, so $b_3=0$, and for $n(2r_n)^m\le1$
\[b_1\le\tfrac52 n^{5}(2r_n)^{4m},\qquad b_2\le\tfrac52 n^{5}(2r_n)^{4m}+\tfrac92 n^{4}(2r_n)^{3m},\]
\[\big|\mathrm{Var}[N_\vee]-\lambda_n\big|\le b_1+b_2\le10 n^{4}(2r_n)^{3m} .\]
\end{restated}
\begin{proof}\emph{Independence.} $I_\alpha$ is a function of $(X_u)_{u\in V_\alpha}$ and $(I_\beta)_{\beta\notin B_\alpha}$ is a function of $(X_u)_{u\notin V_\alpha}$ (their vertex sets are disjoint from $V_\alpha$); these two families of points are independent.

\emph{Counting.} Given $\alpha$, the number of $\beta$ with $V_\beta\cap V_\alpha\ne\emptyset$ is at most $3\binom{n-1}2\cdot3\le\frac92n^2\le5n^2$ (choose a shared vertex, two more vertices, a centre), and the number with $|V_\beta\cap V_\alpha|\ge2$ is at most $3(n-2)\cdot3\le9n$ (choose which two vertices are shared, the third vertex and the centre).

\emph{$b_1$.}
\Line{b_1=\sum_\alpha\sum_{\beta\in B_\alpha}p_\alpha p_\beta\le\frac{n^3}2\cdot5n^2\cdot (2r_n)^{4m}=\frac52n^5(2r_n)^{4m}=\frac52n^{5}(2r_n)^{4m}}
(as $p\le (2r_n)^{2m}$).

\emph{$b_2$.} Let $\beta\in B_\alpha\setminus\{\alpha\}$. (1) If $|V_\alpha\cap V_\beta|=1$: $\{I_\alpha=I_\beta=1\}$ forces the two edges of $\alpha$ and the two edges of $\beta$ to be present; these four edges are distinct and form a connected graph on $|V_\alpha\cup V_\beta|=5$ vertices, i.e.\ a tree, so by Lemma~\ref{lem:tree} (Section~\ref{ts:lem:tree}) the probability is at most $(2r_n)^{4m}$. There are at most $\frac{n^3}2\cdot5n^2$ such ordered pairs, contributing $\le\frac52n^5(2r_n)^{4m}=\frac52n^{5}(2r_n)^{4m}$. (2) If $|V_\alpha\cap V_\beta|=2$: the union has $4$ vertices and is connected (the two triples overlap), so it contains a spanning tree with three present edges; probability $\le (2r_n)^{3m}$. At most $\frac{n^3}2\cdot9n$ pairs: contribution $\le\frac92n^4(2r_n)^{3m}=\frac92n^{4}(2r_n)^{3m}$. (3) If $V_\beta=V_\alpha$ and $\beta\ne\alpha$, then $\beta$ has a different centre, say $\alpha=(j;\{i,k\})$ and $\beta=(i;\{j,k\})$; $I_\alpha=1$ requires $A_{ik}=0$ while $I_\beta=1$ requires $A_{ik}=1$ (the same for the third centre), so $\E [I_\alpha] I_\beta=0$.

\emph{Variance.} \smallskip\noindent\emph{(a) Expansion.} Write $N_\vee=\sum_{\alpha\in\mathcal I}I_\alpha$, so that \Line{\Var [N_\vee]=\sum_{\alpha}\sum_{\beta\in\mathcal I}\mathrm{Cov}(I_\alpha,I_\beta).} If $\beta\notin B_\alpha$, then $I_\beta$ is one of the variables $(I_{\beta'})_{\beta'\notin B_\alpha}$, and $I_\alpha$ is independent of this family, so $\mathrm{Cov}(I_\alpha,I_\beta)=0$. The remaining terms are those with $\beta\in B_\alpha$, and $B_\alpha$ contains $\alpha$ itself. For $\beta=\alpha$ the covariance is $\Var [I_\alpha]=p_\alpha(1-p_\alpha)$ because $I_\alpha$ is a Bernoulli variable. Hence
\Line{\Var [N_\vee]=\sum_\alpha p_\alpha(1-p_\alpha)+\sum_\alpha\sum_{\beta\in B_\alpha\setminus\{\alpha\}}\mathrm{Cov}(I_\alpha,I_\beta).}
Since $\lambda_n=\E [N_\vee]=\sum_\alpha p_\alpha$, subtracting $\lambda_n$ gives the exact identity
\LineE{eq:varid}{\Var [N_\vee]-\lambda_n=-\sum_\alpha p_\alpha^2+\sum_\alpha\sum_{\beta\in B_\alpha\setminus\{\alpha\}}\mathrm{Cov}(I_\alpha,I_\beta).}

\smallskip\noindent\emph{(b) Summation.} Since $\mathrm{Cov}(I_\alpha,I_\beta)=\E[I_\alpha I_\beta]-p_\alpha p_\beta$ with both terms non-negative, \Line{|\mathrm{Cov}(I_\alpha,I_\beta)|\le\E[I_\alpha I_\beta]+p_\alpha p_\beta.} By the definitions of \TS{I}, $b_2=\sum_\alpha\sum_{\beta\in B_\alpha\setminus\{\alpha\}}\E[I_\alpha I_\beta]$, and $b_1=\sum_\alpha\sum_{\beta\in B_\alpha}p_\alpha p_\beta$ contains the diagonal terms $\beta=\alpha$, so that
\[
b_1=\sum_\alpha p_\alpha^2+\sum_\alpha\sum_{\beta\in B_\alpha\setminus\{\alpha\}}p_\alpha p_\beta,\qquad\text{in particular}\qquad\sum_\alpha p_\alpha^2\le b_1 .
\]
Applying the triangle inequality to~\eqref{eq:varid} and this bound to each covariance,
\[\begin{aligned}\big|\Var [N_\vee]-\lambda_n\big|&\le\sum_\alpha p_\alpha^2+\sum_\alpha\sum_{\beta\in B_\alpha\setminus\{\alpha\}}\Big(\E[I_\alpha I_\beta]+p_\alpha p_\beta\Big)\\
&=\sum_\alpha p_\alpha^2+b_2+\Big(b_1-\sum_\alpha p_\alpha^2\Big)=b_1+b_2 .\end{aligned}\]
Finally, $b_1+b_2\le10 n^{4}(2r_n)^{3m}$ was shown above, so \Line{|\Var[N_\vee]-\lambda_n|\le b_1+b_2\le10 n^{4}(2r_n)^{3m},} which is the last claim of the lemma.
\end{proof}

\subsection{Proof of Lemma~\ref{lem:VarD}}\label{pf:lem:VarD}
\mbox{}\par\begin{restated}{Lemma}{lem:VarD}{ (Variance of $\Xi$)}
For $n(2r_n)^m\le\frac1{2e}$, $\Var [\Xi]\le C n^{4}(2r_n)^{3m}$.
\end{restated}
\begin{proof}We use the Efron--Stein inequality (\TS{C}) and Lemma~\ref{lem:csize} (Section~\ref{ts:lem:csize}) to prove Lemma~\ref{lem:VarD} (Section~\ref{ts:lem:VarD}). Let $\Xi^{(k)}$ be $\Xi$ with $X_k$ replaced by an independent copy $X_k'$; it suffices to show \[\E[(\Xi-\Xi^{(k)})^2]\le C(n(2r_n)^m)^3,\] since then $\Var[\Xi]\le\frac12n C(n(2r_n)^m)^3$.

Write $\Xi(\mathsf H)=\sum_{\mathcal C}\mathrm{ex}(\mathcal C)$ over components of a graph $\mathsf H$, with $\mathrm{ex}(\mathcal C)=f(\mathcal C)-N_\vee(\mathcal C)$ computed from $\mathcal C$ alone. By Lemma~\ref{lem:comp} (Section~\ref{ts:lem:comp}), $\mathrm{ex}(\mathcal C)=0$ if $|\mathcal C|\le3$ and $|\mathrm{ex}(\mathcal C)|\le|\mathcal C|^3$ (both $f,N_\vee\ge0$).

Let $\mathsf G^-$ be the graph without $X_k$ and $\mathcal C_k$ the component of $k$ in $\mathsf G$. The components of $\mathsf G$ not containing $k$ are components of $\mathsf G^-$, and the remaining components of $\mathsf G^-$ are the pieces $\mathcal P_1,\dots,\mathcal P_s$ of $\mathcal C_k\setminus\{k\}$, with $\sum|\mathcal P_i|=|\mathcal C_k|-1$. Hence \Line{\Xi(\mathsf G)-\Xi(\mathsf G^-)=\mathrm{ex}(\mathcal C_k)-\sum_i\mathrm{ex}(\mathcal P_i).} If $|\mathcal C_k|\le3$ this is $0$ (all $|\mathcal P_i|\le2$); otherwise, using $\sum a_i^3\le(\sum a_i)^3$, it is at most $|\mathcal C_k|^3+(|\mathcal C_k|-1)^3\le2|\mathcal C_k|^3$. So \Line{|\Xi(\mathsf G)-\Xi(\mathsf G^-)|\le2|\mathcal C_k|^3\1\{|\mathcal C_k|\ge4\}.}

The resampled graph $\mathsf G^{(k)}$ has the same $\mathsf G^-$, so the same bound holds with the component $\mathcal C_k'$ of the new vertex. Write \Line{\Xi-\Xi^{(k)}=[\Xi(\mathsf G)-\Xi(\mathsf G^-)]-[\Xi(\mathsf G^{(k)})-\Xi(\mathsf G^-)].}
Then
\Line{\E\big[(\Xi-\Xi^{(k)})^2\big]\le8 \E\big[|\mathcal C_k|^6\1\{|\mathcal C_k|\ge4\}\big]+8 \E\big[|\mathcal C_k'|^6\1\{|\mathcal C_k'|\ge4\}\big]\le16 C_6 \big(n(2r_n)^m\big)^3,}
since $\mathsf G^{(k)}\overset d=\mathsf G$, the vertices are exchangeable, and Lemma~\ref{lem:csize} (Section~\ref{ts:lem:csize}) with $p=6$ applies for $n(2r_n)^m\le\frac1{2e}$. Hence $\Var[\Xi]\le8C_6 n^4(2r_n)^{3m}$.
\end{proof}

\subsection{Proof of Lemma~\ref{lem:poisBE}}\label{pf:lem:poisBE}
\mbox{}\par

\noindent The classical Berry--Esseen theorem for i.i.d.\ sums is used  to prove Lemma~\ref{lem:poisBE} (Section~\ref{ts:lem:poisBE}), by writing a Poisson variable as a sum of i.i.d.\ Poisson variables.

\begin{tool}[Berry--Esseen for i.i.d.\ sums; Berry 1941, Esseen 1942]\label{tl:J}\tlsec{J}
There is an absolute constant $C_{\mathrm{BE}}$ such that for i.i.d.\ $Y_1,\dots,Y_N$ with mean $\mu$, variance $\varsigma^2>0$ and $\rho=\E[|Y_1-\mu|^3]<\infty$,
\Line{\sup_x\Big|\P\Big(\frac{\sum_{i=1}^N(Y_i-\mu)}{\varsigma\sqrt N}\le x\Big)-\Phi(x)\Big|\le\frac{C_{\mathrm{BE}} \rho}{\varsigma^3\sqrt N}.}
\end{tool}

\begin{restated}{Lemma}{lem:poisBE}{ (Berry--Esseen for the Poisson law)}
There is an absolute constant $C_{\mathrm P}$ such that for all $\lambda>0$,
\Line{\dK\Big(\frac{\mathrm{Po}(\lambda)-\lambda}{\sqrt\lambda},\mathcal{N}\Big)\le C_{\mathrm P} \lambda^{-1/2}.}
\end{restated}
\begin{proof}We use \TS{J} (the Berry--Esseen theorem) to prove Lemma~\ref{lem:poisBE} (Section~\ref{ts:lem:poisBE}). For any integer $N\ge1$, $\mathrm{Po}(\lambda)$ has the same law as a sum of $N$ i.i.d.\ $\mathrm{Po}(\theta)$ variables $Y_i$, $\theta=\lambda/N$, with $\mu=\theta$ and $\varsigma^2=\theta$. For the third absolute central moment, for $\theta\le1$,
\begin{align*}
\rho=\E\big[|Y_1-\theta|^3\big]&=e^{-\theta}\Big[\theta^3+(1-\theta)^3\theta+\sum_{k\ge2}(k-\theta)^3\frac{\theta^k}{k!}\Big]\\
&\le\theta^3+\theta+\theta^2\sum_{k\ge2}\frac{k^3}{k!}\\
&\le\theta+(1+5e) \theta^2\\
&\le2\theta\qquad\text{if }\theta\le\tfrac1{15}.
\end{align*}
(For $k\ge2$ we used $|k-\theta|\le k$, $\theta^k\le\theta^2$ and $\sum_{k\ge0}k^3/k!=5e$.) Choose $N\ge15\lambda$. Theorem J gives
\Line{\dK\le C_{\mathrm{BE}} 2\theta/(\theta^{3/2}\sqrt N)=2C_{\mathrm{BE}}(\theta N)^{-1/2}=2C_{\mathrm{BE}}\lambda^{-1/2}.}
\end{proof}

\end{document}